%11MAR2004
\input amstex
\documentstyle{amsppt}
\document

\magnification 1100

\def\gen{{\frak{g}}}

\def\aen{\frak{a}}
\def\ben{\frak{b}}
\def\uen{\frak{u}}
\def\ven{\frak{v}}
\def\hen{\frak{h}}

\def\len{\frak{l}}
\def\nen{\frak{n}}
\def\pen{\frak{p}}
\def\qen{\frak{q}}
\def\cen{\frak{c}}
\def\sen{\frak{s}}

\def\zen{\frak{z}}

\def\slen{\sen\len}

\def\a{{\alpha}}
\def\o{{\omega}}
\def\O{{\Omega}}

\def\l{{\lambda}}
\def\g{{\gamma}}

\def\b{{\beta}}
\def\eps{{\varepsilon}}

\def\1b{{\bold 1}}

\def\bb{{\bold b}}

\def\cb{{\bold c}}
\def\db{{\bold d}}

\def\fb{{\bold f}}
\def\hb{{\bold h}}

\def\kb{{\bold k}}

\def\ob{{\bold o}}

\def\Ab{{\bold A}}
\def\Bb{{\bold B}}
\def\Cb{{\bold C}}
\def\Db{{\bold D}}

\def\Fb{{\bold F}}
\def\Gb{{\bold G}}
\def\Ib{{\bold I}}
\def\Jb{{\bold J}}
\def\Kb{{\bold K}}
\def\Mb{{\bold M}}
\def\Nb{{\bold N}}
\def\Hb{{\bold H}}

\def\Qb{{\bold Q}}
\def\Rb{{\bold R}}

\def\Tb{{\bold T}}
\def\Lb{{\bold L}}
\def\bbu{{\underline\bb}}
\def\dbu{{\underline\db}}
\def\fbu{{\underline\fb}}
\def\Bbu{{\underline\Bb}}
\def\Fbu{{\underline\Fb}}
\def\Lbu{{\underline\Lb}}
\def\Ub{{\bold U}}

\def\Modb{{\bold{Mod}}}
\def\Modbu{{\underline\Modb}}

\def\Qcohb{{\Qb\cb\ob\hb}}
\def\Indb{{\bold{Ind}}}
\def\Prob{{\bold{Pro}}}
\def\Extb{{\bold{Ext}}}
\def\Homb{{\bold{Hom}}}
\def\Endb{{\bold{End}}}

\def\Setsb{{\bold{Sets}}}

\def\benu{{\underline\ben}}
\def\cenu{{\underline\cen}}
\def\genu{{\underline\gen}}
\def\henu{{\underline\hen}}

\def\senu{{\underline\sen}}

\def\chiu{{\underline\chi}}
\def\phiu{{\underline\phi}}

\def\du{{\underline d}}
\def\eu{{\underline e}}
\def\fu{{\underline f}}
\def\hu{{\underline h}}

\def\su{{\underline s}}

\def\xu{{\underline x}}

\def\Gu{{\underline G}}

\def\Mu{{\underline M}}
\def\Hu{{\underline H}}
\def\Gbu{{\underline\Gb}}
\def\Hbu{{\underline\Hb}}

\def\Mbu{{\underline\Mb}}
\def\Tbu{{\underline\Tb}}

\def\Iu{{\underline I}}

\def\Wu{{\underline W}}
\def\Xu{{\underline X}}
\def\Yu{{\underline Y}}

\def\Bcu{{\underline\Bc}}
\def\Ccu{{\underline\Cc}}
\def\Zcu{{\underline\Zc}}

\def\Ncu{{\underline\Nc}}
\def\Scu{{\underline\Sc}}
\def\Phiu{{\underline\Phi}}
\def\Psiu{{\underline\Psi}}
\def\Piu{{\underline\Pi}}

\def\h{{\roman h}}

\def\t{{\roman t}}

\def\R{{\roman R}}

\def\top{{\text{top}}}
\def\ad{\roman{ad}}

\def\ch{\roman{ch}}
\def\con{\roman{con}}

\def\Spec{\roman{Spec}}

\def\Id{\roman{Id}}

\def\Rep{\text{Rep}\,}

\def\simto{\,{\buildrel\sim\over\to}\,}

\def\End{\roman{End}}

\def\IM{\roman{IM}}

\def\Ker{\text{Ker}\,}

\def\Lie{\text{Lie}}
\def\GL{{\text{GL}}}

\def\SL{{\text{SL}}}
\def\SP{{\text{SP}}}

\def\mod{\text{mod}\,}
\def\Oplus{\ts\bigoplus}

\def\Sum{\ts\sum}

\def\CC{{\Bbb C}}
\def\DD{{\Bbb D}}

\def\QQ{{\Bbb Q}}
\def\RR{{\Bbb R}}

\def\ZZ{{\Bbb Z}}

\def\Bc{{\Cal B}}
\def\Cc{{\Cal C}}
\def\Dc{{\Cal D}}
\def\Ec{{\Cal E}}
\def\Fc{{\Cal F}}
\def\Gc{{\Cal G}}

\def\Lc{{\Cal L}}

\def\Nc{{\Cal N}}
\def\Oc{{\Cal O}}
\def\Pc{{\Cal P}}
\def\Qc{{\Cal Q}}

\def\Sc{{\Cal S}}
\def\Tc{{\Cal T}}
\def\Uc{{\Cal U}}
\def\Vc{{\Cal V}}

\def\Xc{{\Cal X}}
\def\Yc{{\Cal Y}}
\def\Zc{{\Cal Z}}

\def\and{{\ \text{and}\ }}

\def\ts{\textstyle}                
\def\ss{\scriptstyle}                
                
\def\qed{\hfill $\sqcap \hskip-6.5pt \sqcup$}        % White box                
\overfullrule=0pt                                    % No black boxes

\def\lra{{{\longrightarrow}}}
\def\bda{{{\big\downarrow}}}
\def\bua{{{\big\uparrow}}}
\def\ind{{\lim\limits_{\longrightarrow}}}
\def\pro{{\lim\limits_{\longleftarrow}}}

\def\u1{{\underline 1}}

\def\la{{\langle}}
\def\ra{{\rangle}}
\def\lla{{\longleftarrow}}

\newdimen\Squaresize\Squaresize=14pt
\newdimen\Thickness\Thickness=0.5pt
\def\Square#1{\hbox{\vrule width\Thickness
	      \alphaox to \Squaresize{\hrule height \Thickness\vss
	      \hbox to \Squaresize{\hss#1\hss}
	      \vss\hrule height\Thickness}
	      \unskip\vrule width \Thickness}
	      \kern-\Thickness}
\def\Vsquare#1{\alphaox{\Square{$#1$}}\kern-\Thickness}

%%%%%%%%%%%%%%%%%%%%%%%%%%%%%%%%%%%%%%%%%%%%%%%%%%%%%%%%%%%%%%%%%%%%%%%%
\topmatter
\title Induced and simple modules of double affine Hecke algebras\endtitle
\rightheadtext{}
\author Eric Vasserot\endauthor
\address D\'epartement de Math\'ematiques,
Universit\'e de Cergy-Pontoise, 2 av. A. Chauvin,
BP 222, 95302 Cergy-Pontoise Cedex, France\endaddress
\email eric.vasserot\@math.u-cergy.fr\endemail
\thanks
2000{\it Mathematics Subject Classification.}
Primary 17B37; Secondary 17B67, 14M15, 16E20.
\endthanks
\abstract
We classify the simple integrable modules 
of double affine Hecke algebras via perverse sheaves.
We get also some estimate for the Jordan-H\"older multiplicities 
of induced modules.
\endabstract
\endtopmatter
\document

\head Contents\endhead

1. Introduction

2. Reminder on flag manifolds 

3. Reminder on double affine Hecke algebras

4. The convolution algebra

5. Induced modules

6. The regular case

7. Induction of sheaves and Fourier transform

8. The type A case

9. Other examples

A. Appendix : K-theory

B. Appendix : infinite dimensional vector spaces 

\head 1. Introduction\endhead
\subhead 1.0\endsubhead
The present work is partially motivated by some recent progress in the 
representation theory of rational Cherednik algebras, see \cite{BEG2}
and the references therein.
The main tool, there, is an induction functor from modules over 
a Weyl group to modules over the rational algebra.
It is desirable to have a geometric model for
the induction functor from the affine Hecke algebra to the
double affine Hecke algebra, and for the Jordan-H\"older multiplicities. 
To do so the first step is 
to construct geometrically a class of simple modules. 
We do this via equivariant K-theory and perverse sheaves.
In the classification of simple modules two main cases appear.
In the first one we deal with equivariant perverse sheaves on a
variety with a finite number of orbits.
Then the construction is quite parallel to the affine case.
In the second case there is an infinite number of orbits and the picture
is less clear. 
Fortunately a Fourier transform exchanges both cases, 
yielding a complete classification of the simple integrable modules.

\subhead 1.1\endsubhead
The proof uses the Kashiwara affine flag manifold $\Xc$.
It is a non Noetherian scheme which is locally
pro-smooth and of pro-finite type.
The idea to generalize K-theoretic technics from the affine Hecke
algebra to the double affine Hecke algebra is not new.
It was mentioned to us by Ginzburg about 10 years ago,
and was probably known before that.
The idea to use Kashiwara-Tanisaki's tools 
for double affine Hecke algebras 
appears first in \cite{GG}, where a geometric construction 
of the regular representation is sketched. 
Our construction is different from {\it loc. cit.}
We construct a ring and not only a faithfull representation.
Moreover it is essential for us to
do a reduction to the fixed point subset
respectively to some semisimple element whose
`rotation' component is not a root of unity.
The case of roots of unity is still open.

Let $\Hb$ be the double affine Hecke algebra.
One Important step in our work is Theorem 4.9
which yields a ring homomorphism from $\Hb$ to 
a ring defined via the equivariant K-theory of 
an analogue of the Steinberg variety.
Let us sketch briefly this construction.

We first define an affine analogue of the Steinberg variety,
which we denote by $\Zc$.
It is an ind-scheme of ind-infinite type.
It comes with a filtration by subsets $\Zc_{\leq y}$, 
with $y$ in the affine Weyl group $W$. 
The subsets $\Zc_{\leq y}$ are reduced separated schemes of infinite type,
and the natural inclusions $\Zc_{\leq y'}\subseteq\Zc_{\leq y}$,
with $y'\leq y$ are closed immersions.
The set $\Zc$ is endowed with an action of a torus $A$
wich preserves each term of the filtration.

For a well-chosen element $a\in A$,
the fixed point set $\Zc^a\subset\Zc$ 
is a scheme locally of finite type. 
Hence we have a convolution ring $\Kb^A(\Zc^a)$ : 
it is the inductive limit of the system of $\Rb_A$-modules 
$\Kb^A((\Zc_{\leq y})^a)$ with $y\in W$.
Each $\Kb^A((\Zc_{\leq y})^a)$ is itself the projective limit
of the system of $\Rb_A$-modules $\Kb^A(({}^w\Zc_{\leq y})^a)$ with $w\in W$ 
and ${}^w\Zc_{\leq y}\subset\Zc_{\leq y}$ an open subscheme such that 
$({}^w\Zc_{\leq y})^a$ is of finite type.
See 4.4.

Although the convolution product on $\Kb(\Zc^a)$ is easily defined, it is not
so easy to construct a ring homomorphism $\Hb\to\Kb^A(\Zc^a)$
without using the whole $\Zc$.
Let $\Hbu$ be the affine Hecke algebra. 
Recall that in the geometric construction of $\Hbu$ in \cite{CG} or \cite{KL1}, 
one constructs a ring isomorphism from $\Hbu$ to the equivariant
K-theory group of the Steinberg variety $\Zcu$.
However it is more delicate to construct directly
a ring homomorphism $\Hbu\to\Kb(\Zcu^a)$.
More precisely, although it is quite possible to assign elements
in $\Kb(\Zcu^a)$ to the standard generators of $\Hbu$, it is
delicate to check the required relations without using the whole 
variety $\Zcu$. 
This explains several technical difficulties in our work.

We do not define a convolution ring $\Kb^A(\Zc)$.
This is not necessary.
To construct a ring homomorphism 
$\Psi_a\,:\,\Hb\to\Kb^A(\Zc^a)_a$ 
we proceed as follows :  

{\sl Step 1.} 
For any non empty proper set $J$ of simple roots,
let $W_J\subsetneq W$ be the corresponding parabolic subgroup.
It is well-known that $W_J$ is a Weyl group of finite type.
Let $G_J$ be the connected, simply connected, simple group
whose Weyl group is $W_J$. 
Let $\Zc_J$ be the Steinberg variety of $G_J$.
Let $\Hb_J$ be the affine Hecke algebra associated to $G_J$.
There is a ring isomorphism 
$\Psi_J\,:\,\Hb_J\to\Kb^{G_J\times\CC^\times}(\Zc_J)$.

{\sl Step 2.} 
Let $w_J\in W_J$ be the longuest element.
The Grothendieck group $\Kb^A({}^{w}\Zc_{\leq w_J})$ of $A$-equivariant
coherent sheaves on ${}^w\Zc_{\leq w_J}$ makes sense.
It is the inductive limit of 
$\Kb^A({}^{w}\Zc^{k\ell}_{\leq w_J})$, where 
${}^{w}\Zc^{k\ell}_{\leq w_J}$ is a Noetherian quotient scheme 
of ${}^{w}\Zc_{\leq w_J}$.
The groups $\Kb^A({}^{w}\Zc_{\leq w_J})$ 
are endowed with a convolution product.
Remark that the natural embedding 
${}^w\Zc^{k\ell}_{\leq y'}\subset{}^w\Zc^{k\ell}_{\leq y}$
may be not closed, see Remark 4.2.
This explains a technical complication
in the definition of the product on $\Kb^A({}^{w}\Zc_{\leq w_J})$.
See 4.5.
For a similar reason we have to consider the Thomason
concentration isomorphism on non separated schemes.
A reminder is given in the appendix.
Induction of sheaves yields a map
$\Kb^{G_J\times\CC^\times}(\Zc_J)\to\Kb^A({}^{w}\Zc_{\leq w_J})$.
Composing it with the bivariant concentration map and the natural map
$\Kb^A(({}^w\Zc_{\leq w_J})^a)_a\to\Kb^A(\Zc^a)_a$ 
yields a ring homomorphism
$\Theta_J\,:\,\Kb^{G_J\times\CC^\times}(\Zc_J)\to\Kb^A(\Zc^a)_a$.

%$\Psi_{J,a}\,:\,\Hb_J\to\Kb^A((\Zc_J)^a)$.
{\sl Step 3.} 
The maps $\Theta_J$ are compatible in the following sense.
If $J_1\subset J_2$ then $\Hb_{J_1}$ embeds in $\Hb_{J_2}$ and
the composition of the chain of maps
$$\Hb_{J_1}\to\Hb_{J_2}\to\Kb^{G_{J_2}\times\CC^\times}(\Zc_{J_2})
\to\Kb^A(\Zc^a)_a$$
equals $\Theta_{J_1}\circ\Psi_{J_1}$.
Now each algebra $\Hb_J$ embeds as a parabolic subalgebra of $\Hb$,
and $\Hb$ is obviously generated by $\bigcup_J\Hb_J$.
The compatibility condition above insures that the morphisms 
$\Theta_{J}\circ\Psi_J$
glue together in a morphism 
$\Psi_a\,:\,\Hb\to\Kb^A(\Zc^a)_a$.
It is uniquely determined. 

The map $\Psi_a$ becomes surjective 
after a suitable completion of $\Hb$.
It is certainly not injective.
Using $\Psi_a$, a standard sheaf-theoretic construction,
due to Ginzburg in the $\Hbu$-case,
gives us a collection of simple $\Hb$-modules.
They are precisely the simple integrable modules.
Remark that it is important for us to construct
the convolution product on $\Kb^A(\Zc^a)_a$ and 
the ring homomorphism $\Psi_a$, rather than the regular representation
of $\Hb$ as in \cite{CG}, in order to use Ginzburg's machinery.
Otherwise, we can not get a geometric interpretation
for the Jordan-H\"older multiplicities of induced modules
as in 7.5.

We conclude this subsection by two remarks :

First, as far as fixed point subsets are concerned, 
we may consider only the ind-scheme $\Bc$, see 2.7 for the notation,
rather than Kashiwara's scheme $\Xc$ : 
although $\Xc$ contains $\Bc$ as a proper subset,
the fixed point subsets $\Bc^a$, $\Xc^a$ coincide,
see Lemma 2.13.$(ii)$.
However, we use $\Kb^A({}^{w}\Zc_{\leq w_J})$
in the definition of $\Theta_{J}$ and to check
the compatibilities in Step 3, see Lemma 4.7.
Recall that ${}^w\Zc_{\leq w_J}$ is a scheme of infinite type.
We do not know how to make sense of the Grothendieck group
of coherent sheaves with infinite dimensional support on $\Zc_{\leq w_J}$,
i.e. $\Kb^A({}^{w}\Zc_{\leq w_J})$,
without using Kashiwara's affine flag manifold.
For instance, viewing $\Zc_{\leq w_J}$ as an ind-scheme rather than a scheme
would only give us coherent sheaves with finite dimensional support 
on $\Zc_{\leq w_J}$.
See Section 4, in particular the preliminary subsection, for more details.

Second, our construction is similar to Lusztig's action
of the affine Weyl group on the homology of
the affine Springer fibers in \cite{L5}. 
Indeed, one can prove that the actions of $\Hb$ and $W$
are compatible in some sense, see 9.4 for details.
Our work is more technical than \cite{L5},
because Lusztig does not use fixed points sets relatively to $a$,
and because he does not construct any simple module.

\subhead 1.2\endsubhead
Let us say a few words on the structure of this paper.
It contains nine sections.

Section 2 contains basic results on reductive groups
and flag manifolds.
The definition of the Kashiwara flag manifold 
and some basic property are given in 2.3-5. 
The scheme $\Xc$ is locally of infinite type.
We construct a scheme $\Tc$ which is a substitute to the
cotangent bundle of $\Xc$.
The definition of $\Tc$
and its basic properties are given in 2.6-7. 
The fixed point subsets of $\Xc$, $\Tc$ are studied in 2.13.

Section 3 contains the definition of the double affine 
Hecke algebra $\Hb$. 

Section 4 contains the definition of $\Zc$. 
Basic properties of $\Zc$ are given in 4.1-2.
The convolution product on the equivariant K-groups 
of $\Zc^a$ is given in 4.4.
The convolution algebra is related to $\Hb$ in 4.6-9.

Section 5 contains a geometric construction of some induced $\Hb$-modules.
It is an affine analogue of \cite{L3}.

In Sections 6,7 we state and prove the main theorem. 
The simple and induced modules in the regular case are studied in Section 6.
The classification of the simple modules is achieved in 7.4, 7.6,
the study of the induced modules in 7.5.

The type $A$ is described in Section 8.

Section 9 contains examples.

Appendix A contains recollections on equivariant K-theory.
In particular we consider the case of non-separated schemes.
Appendix B contains recollections on endomorphism algebras 
of infinite dimensional vector spaces.

\subhead 1.3\endsubhead
Let us insist on a few points which are important 
for the representation theory of the double affine Hecke algebra.

First, when the parameters are generic the
representations induced from irreducible representations
of the affine Hecke algebra remain irreducible, 
due to an argument of Cherednik via intertwiners.
Geometrically, it means that the nilpotent orbits
which contains the Jordan-H\"older multiplicities of induced modules
are the same in the affine and the double affine cases,
see 9.1 for details.
Thus the main problem is in the regular case, 
see Definition 2.14.

Second, in type $A$ case Cherednik has recently given a combinatorial
approach to the representation of double affine Hecke algebras.
It is coherent with our description, see 8.2-3.
Moreover, in this case the Jordan-H\"older multiplicities of the induced modules
are related to the canonical basis of the cyclic quiver.
It is quite remarkable that, for type $A$ , the multiplicities are given by 
Kazhdan-Lusztig polynomials of type $A^{(1)}$.
The appearance of Kazhdan-Lusztig polynomials was already
suggested in a less precise way in \cite{AST}.
Theorem 8.5 gives an algorithm to compute 
the Jordan-H\"older multiplicities of induced modules. 
It is inspired from \cite{A}, \cite{VV}.

Finally, finite dimensional representations are still unknown.
However a special class of finite dimensional modules, which where introduced
previously by Cherednik (see \cite{C2, Theorem 8.5} for instance),
can be nicely described in geometrical terms.
There are related to some regular semi-simple nil-elliptic elements
introduced by Lusztig a few years ago, see 9.3.

\subhead 1.4. General conventions\endsubhead
By scheme we always mean scheme over $\CC$.
By variety we mean a reduced separated scheme of finite type.
Let $G$ be a linear group.
Set $G^\circ$ equal to the identity component of $G$.

Let $\Xc$ be any Noetherian $G$-scheme 
(i.e. a Noetherian scheme with a regular action of the group $G$).
Set $\Kb^G(\Xc)$ equal to the complexified Grothendieck group 
of the Abelian category of $G$-equivariant coherent sheaves on $\Xc$.
Set $\Kb_G(\Xc)$ equal to the complexified Grothendieck group 
of the exact category of $G$-equivariant locally free sheaves on $\Xc$.
Let $\Kb^G_i(\Xc)$, $\Kb_G^i(\Xc)$, $i\geq 1$,
denote the higher complexified K-groups.
We may write $\Kb^G_0(\Xc)$, $\Kb_G^0(\Xc)$
for $\Kb^G(\Xc)$, $\Kb_G(\Xc)$.

We write $\Rb_G$ for $\Kb_G(point)$.
If $G$ is a diagonalizable group and $g\in G$,
let $\Jb_g\subset\Rb_G$ be the maximal ideal associated to $g$,
and $\CC_g=\Rb_G/\Jb_g$.
For any $\Rb_G$-module $\Mb$ we set $\Mb_g=\Mb\otimes_{\Rb_G}\CC_g$.
See the appendix for other notations relative to K-theory.

For any subset $S\subset G$ we write 
$\Xc^S\subset\Xc$ for the fixed point subset.
Let $\delta_\Xc$ denote the diagonal immersion $\Xc\to\Xc\times\Xc$.
Let $q_{ij}\,:\,\Xc_1\times\Xc_2\times\Xc_3\to\Xc_i\times\Xc_j$, 
$q_i\,:\,\Xc_1\times\Xc_2\to\Xc_i$ 
denote the obvious projections.

Given a one-parameter subgroup $\rho\,:\,\CC^\times\to G$,
we write $\CC^\times_\rho$ for $\rho(\CC^\times)$.
We write $\la g\ra\subset G$ for the closed subgroup generated by 
the element $g\in G$, and $N_G(H)$ for the normalizer of 
a subgroup $H\subset G$.
%Given $g\in G$, $x\in\Lie(G)$,
%we may write $g(x)$ for the element $\ad_gx$.

Assume now that $\Xc$ is a variety. 
Let $\Hb_*(\Xc)$ denote the complex Borel-Moore homology of $\Xc$.
Occasionally we consider the equivariant topological
K-homology groups $\Kb_{i,\top}^G(\Xc)$, with $i=0,1$.
As above we may write $\Kb_\top^G(\Xc)$ for $\Kb_{0,\top}^G(\Xc)$.
Let also $\Db(\Xc)$ be the derived category
of bounded complexes of $\CC$-sheaves on $\Xc$
whose cohomology groups are constructible sheaves,
and $\Db_G(\Xc)$ be the triangulated category
of $G$-equivariant complexes of constructible sheaves on 
$\Xc$ constructed in \cite{BL}.
Let $\CC_\Xc$ be the constant sheaf,
and $\DD_\Xc$ be the dualizing complex.
For any $G$-equivariant local system $\chi$ 
on a locally closed set $\Oc\subset\Xc$, 
let $IC_\chi$ be the corresponding intersection cohomology complex.
It is a simple perverse sheaf if $\chi$ is an irreducible local system.
If $\Oc$ is a $G$-orbit and $e\in\Oc$,
we may identify a representation of the component group of
the isotropy subgroup of $e$ in $G$ with the corresponding 
$G$-equivariant local system on $\Oc$.

Given graded vector spaces $V,W$ we write $V\dot{=}W$
for a linear isomorphism that does not necessarily preserve the gradings.
We will also use the notation $\dot{=}$ to denote quasi-isomorphisms 
that only hold up to a shift in the derived category. 

If $V$ is a finite dimensional vector space, let $V^*$ be the
dual space and $\Fb_V\,:\,\Db_\con(V)\to\Db_\con(V^*)$ 
be the Fourier transform.
Here $\Db_\con(V)$ is the full subcategory of
$\Db(V)$ of conic sheaves, 
i.e. complexes of sheaves whose cohomology is
constant on the $\RR^+$-orbits.
Set $\CC^\times_V=\CC^\times\cdot\Id_V\subset\GL(V)$.
The $\CC^\times_V$-equivariant perverse sheaves are conic.
If $E$ is a $G\times\CC^\times_V$-equivariant 
perverse sheaf then $\Fb_V(E)$
is a $G\times\CC^\times_{V^*}$-equivariant perverse sheaf on $V^*$,
and $\Fb_{V^*}\circ\Fb_V(E)\simeq E$.
If moreover $E$ is simple then $\Fb_V(E)$ is also simple.

Given an algebra $\Ab$, let $\Modb_\Ab$ 
denote the category of left $\Ab$-modules.
A morphism of algebras $c\,:\,\Ab\to\Bb$ 
induces a pull-back map $\Modb_\Bb\to\Modb_\Ab$,
still denoted by $c^\bullet$.

\subhead Acknowledgements\endsubhead
I wish to thank G. Laumon and V. Kac for their help during the
preparation of this paper, and I. Cherednik for his interest.
While writing this paper I enjoy the hospitality of MSRI and IHES.
The author was also partially supported by EEC grant
no. ERB FMRX-CT97-0100.

\head 2. Reminder on flag manifolds\endhead

Sections 2.0 to 2.13 contain notations and basic facts on flag manifolds.
Proposition 2.14 will be used in the sheaf theoretic approach to the 
convolution algebra in Section 6.

\subhead 2.0\endsubhead
Let $\Gu$ be a simple connected simply-connected linear algebraic group.
Let $\Gu$ denote also the $\CC$-points in $\Gu$. 
Let $\genu$ be the Lie algebra of $\Gu$.
Let $\henu\subset\genu$ be a Cartan subalgebra, and let
$\benu\subset\genu$ be a Borel subalgebra containing $\henu$.
Let $\Hu\subset\Gu$ be the Cartan subgroup corresponding to $\henu$.
Let $\Phiu$ be the root system of $\genu$.
For any $\a\in\Phiu$ let $\genu_\a\subset\genu$ 
be the corresponding root subspace.
Let $\Phiu^\vee$ be the root system dual to $\Phiu$.
Set $\Xu,\Yu\subset\henu^*$ 
(resp. $\Xu^\vee,\Yu^\vee\subset\henu$)
equal to the weight and the root lattices of $\Phiu$
(resp. of $\Phiu^\vee$).
Recall that, $\Gu$ being simply connected,
$\Xu$ is isomorphic to the group of characters of $\Hu$,
while $\Yu^\vee$ is isomorphic to the group of the
one-parameter subgroups of $\Hu$.
Let $\a_i,\o_i\in\henu^*$, $\a_i^\vee,\o^\vee_i\in\henu$, $i\in\Iu$, 
be the simple roots, the fundamental weights,
the simple coroots and the fundamental coweights.
Hence $\Yu\subset\Xu$, $\Yu^\vee\subset\Xu^\vee$
and $\{\a_i\}$, $\{\a^\vee_i\}$,
$\{\o_i\}$, $\{\o^\vee_i\}$ are bases of $\Yu$, $\Yu^\vee$, 
$\Xu$, $\Xu^\vee$ respectively.
Let $\theta\in\Phiu$ be the highest root, 
and $\theta^\vee$ be the corresponding coroot.
Let $\Wu$, $W$ be the Weyl group and the affine Weyl group of $\Gu$.
Let $\ell(w)$ denote the length of the element $w\in W$.

Set $I=\Iu\sqcup\{0\}$.
Consider the lattices
$Y=\bigoplus_{i\in I}\ZZ\a_i\subset 
X=\ZZ\delta\oplus\bigoplus_{i\in I}\ZZ\o_i$,
$Y^\vee=\bigoplus_{i\in I}\ZZ\a^\vee_i\subset 
X^\vee=\ZZ c\oplus\bigoplus_{i\in I}\ZZ\o_i^\vee$,
where $\a_0$, $\a_0^\vee$, $\o_0$, $\o_0^\vee$, $\delta$, $c$
are new variables such that 
$\a_0=\delta-\theta,$ and $\a^\vee_0=c-\theta^\vee.$
There is a unique pairing $X\times Y^\vee\to\ZZ$ such that
$(\o_i:\a_j^\vee)=\delta_{ij}$ and 
$(\delta:\a_j^\vee)=0$.

We identify $\Iu$ (resp. $I$) with the set 
of simple reflexions in $\Wu$ (resp. in $W$).
Let $s_i\in W$ be the simple reflexion corresponding to $i\in I$.
For all $i,j\in I$ let $m_{ij}$ 
be the order of the element $s_is_j$ in $W$.
Let also $\O$ denote the Abelian group $\Xu^\vee/\Yu^\vee$.

\vskip3mm

\subhead 2.1\endsubhead
Set $K=\CC((\eps))$, and $K_d=\CC((\eps^{1/d}))$ for any $d\in\ZZ_{>0}$.
Let $\hen$ be the Abelian Lie algebra on the vector space
$X^\vee\otimes_\ZZ\CC=\henu\oplus\CC c\oplus\CC\o_0^\vee$.
The sum $\gen=(\genu\otimes K)+\hen$ is endowed with the usual Lie bracket,
such that $c$ is central and $\o_0^\vee$ acts by derivations.

Let $\Phi\subset Y$ be the set of roots of $\gen$,
$\Phi^+\subset\Phi$ the set of positive roots, $\Phi^-=-\Phi^+.$
For any $\a\in\Phi$ let $\gen_\a\subset\gen$ 
be the corresponding weight subspace.
Given $\Theta\subset\Phi^+$
such that $(\Theta+\Theta)\cap\Phi^+\subset\Theta$,
we set $\uen(\pm\Theta)=\prod_{\a\in\pm\Theta}\gen_\a$.
Note that $\uen(\pm\Theta)$ is a
Lie algebra in the category of schemes.

Given a decreasing sequence of subsets
$\Theta_\ell\subseteq\Phi^+$, with $\ell\in\ZZ_{\geq 0}$, such that
$(\Theta_\ell+\Phi^+)\cap\Phi^+\subset\Theta_\ell,$
$|\Phi^+\setminus\Theta_\ell|<\infty,$
and $\bigcap_{\ell\geq 0}\Theta_\ell=\emptyset,$
we set $\uen_\ell=\uen(\Theta_\ell)$ and
$\uen^-_\ell=\uen(-\Theta_\ell)$.
We write $\uen$ and $\uen^-$ for $\uen(\Phi^+)$ and $\uen(\Phi^-)$.
Then $\uen_\ell\subset\uen$ and $\uen^-_\ell\subset\uen^-$
are closed Lie ideal.

\subhead 2.2\endsubhead
Let $G'$ be the central extension of $\Gu(K)$ by $\CC^\times$, see \cite{G}.
The affine Kac-Moody group $G_{KM}\subset G'$ 
is the restriction of the central extension to the subgroup
$\Gu(\CC[\eps^{\pm 1}])\subset\Gu(K),$ see \cite{Kc, 2.8}.
The groups $\Gu(K)$, $G'$ or $G_{KM}$ share the similar properties to the
group $\Gu$, such as Bruhat decomposition.
See \cite{K}, \cite{Kc} and the references therein.
The group $\CC^\times$ acts on $\Gu(K_d)$ by `rotating the loops', 
i.e. $\tau\cdot y(\eps^{1/d})=y(\tau\eps^{1/d})$. 
If $d=1$ this action extends uniquely to $G'$.
Let $G=G'\rtimes\CC^\times$ and $G_d=\Gu(K_d)\rtimes\CC^\times$
be the corresponding semidirect products.
We write $\CC^\times_\delta$ for the subgroup 
$\{1\}\times\CC^\times\subset G, G_d$.
The coroots $\a_i^\vee$, hence $c$, are viewed as 
one-parameter subgroups of $G'$.
We write $\CC^\times_{\o_0}$ for $\la c\ra\subset G'$.

Set also 
$H'=\prod_{i\in I}\la\a_i^\vee\ra\subset G'$,
and $H=H'\times\CC^\times_\delta\subset G$.
Hence $H'=\Spec(\CC[X'])=\Hu\times\CC^\times_{\o_0}$ 
with $X'=X/\ZZ\delta$, and $H=\Spec(\CC[X])$.
The group $W$ acts on $X$ by 
$s_i(\l)=\l-(\l:\a_i^\vee)\a_i$ for all $i\in I$, $\l\in X$. 
This action induces a $W$-action on $H$.
For any $s\in H$ we set $W_s=\{w\in W\,;\,w(s)=s\}$.

Given a finite-dimensional nilpotent Lie algebra $\aen$, 
let $\exp(\aen)$ denote the unipotent 
algebraic group with Lie algebra $\aen$.
We consider the group schemes
$U=\pro_\ell\exp(\uen/\uen_\ell)$ and
$U^-=\pro_\ell\exp(\uen^-/\uen^-_\ell)$.
Let $U(\Theta)\subset U$ and
$U(-\Theta)\subset U^-$ be the closed subgroups
with Lie algebras $\uen(\Theta)$ and $\uen(-\Theta)$.
We write $U_\ell$ and $U_\ell^-$ for $U(\Theta_\ell)$ and $U(-\Theta_\ell)$.
For any real root $\a\in\Phi^+$ we also write
$U_{\pm\a}$ for $U(\{\pm\a\})$.
Let $B=H\ltimes U$, $B^-=H\ltimes U^-$.
The groups $U$, $B$ are identified with closed subgroups of $G$
(for the topology on $G$ induced by the natural topology on $K$).

\subhead 2.3\endsubhead
Let $\Xc$ be the flag manifold of the triple 
$(\gen,\hen,X)$, see \cite{K}. 
Recall that $\Xc$ is a reduced separated scheme, locally of infinite type, 
equal to the quotient of a separated scheme $\Gc$ by
the right action of $B$.
The group $B^-$ acts on $\Gc$ on the left.
The actions of $B$, $B^-$ on $\Gc$ are locally free.
See \cite{K, 5.7} for the definition of locally freeness.
Let $\pi\,:\,\Gc\to\Xc$ be the canonical map.

Let $\dot W$ be the braid group of $W$.
It acts on $\Gc$ from the left and the right.
The left actions of $B^-,$ $\dot W$ 
commute to the right actions of $B$, $\dot W$.
The left $B$-action and right $B^-$-action on $\Gc$ are denoted 
by `a dot', as well as the $\dot W$-action (left or right) on $\Gc$.
The left actions of $B^-,$ $\dot W$  on $\Xc$ is also denoted by `a dot'.

Let $\ad$ denote the adjoint action of
$H$, $\dot W$, $G$ on $\gen$.
We embed $W$ in $\dot W$ as usual : 
let $\dot w\in\dot W$ denote the image of $w$.
For any $H$-stable subset $\Yc\subseteq\Xc$ (resp. $\Yc\subset\gen$)
the set $\dot w\cdot\Yc$ (resp. $\ad_{\dot w}\Yc$) depends only on $w$,
hence we denote it by $w\cdot\Yc$ (resp. $\ad_w\Yc$).
The $\dot W$-action on $\gen$ gives an isomorphism of group schemes 
$U(\Phi^-\cap w^{-1}\Phi^+)\to U(\Phi^+\cap w\Phi^-)$,
$u\mapsto\ad_{\dot w}u$.

Let $1_\Gc\in\Gc$ be the 'identity element', 
see \cite{K, 5.1}, and $1_\Xc=\pi(1_\Gc)$.
We set 
$\Xc_y=yU(\Phi^-\cap y^{-1}\Phi^+)\cdot 1_\Xc$,
$\Xc^w=U^-w\cdot 1_\Xc$,
$\Uc^w=wU^-\cdot 1_\Xc$,
$\Xc_{\leq y}=\bigcup_{y'\leq y}\Xc_{y'}$,
${}^w\Xc=\bigcup_{w'\leq w}\Uc^{w'}$.
The maps
$U(\Phi^-\cap w^{-1}\Phi^-)\to\Xc^w,$ $u\mapsto wu\cdot 1_\Xc$, and
$U^-\to\Uc^w$, $u\mapsto wu\cdot 1_\Xc$ are isomorphisms.
Put also ${}^w\Gc=\pi^{-1}({}^w\Xc)$.
Given $y,y',y''\in W$ we write $y''\succeq(y,y')$
if $y''\geq ww'$ for any $w\leq y,$ $w'\leq y'$.

\proclaim{Lemma} 
Fix $y,y',y'',w,w'\in W$.

(i)
$\dot w u\cdot 1_\Gc=1_\Gc\cdot(\ad_{\dot w}u)w$
for any $u\in U(\Phi^-\cap w^{-1}\Phi^+)$.

(ii) 
If $w'\succeq(w,y)$ then
$wU^-\cdot\Xc_{\leq y}\subset{}^{w'}\Xc$. 

(iii)
If $y''\succeq(y,y')$ then
$yU(\Phi^-\cap y^{-1}\Phi^+)\cdot\Xc_{\leq y'}\subset\Xc_{\leq y''}.$

(iv)
$\Xc_{\leq w}\subset{}^w\Xc$,
${}^w\Xc=\bigcup_{w'\leq w}\Xc^{w'}$,
and $\Xc=\bigcup_{w\in W}\Xc^w$.

(v)
$\Uc^w, {}^w\Xc$ are open in $\Xc$;
$\Xc_y,\Xc^w$ are locally closed in $\Xc$; 
$\Xc_{\leq y}$ is projective.
\endproclaim

\noindent{\sl Proof:}
$(i)$ is proved in
\cite{K, Theorem 5.1.10.$(iii)$} if $\ell(w)=1$.
If $\ell(w)=\ell>1$ we fix a reduced expression 
$w=s_{i_1}s_{i_2}\cdots s_{i_\ell}.$
The element $u\in U(\Phi^-\cap w^{-1}\Phi^+)$ can be decomposed 
as a product
$u=
(\ad_{\dot s_{i_\ell}}^{-1}\cdots\ad_{\dot s_{i_2}}^{-1}u_1)
\cdots
(\ad_{\dot s_{i_\ell}}^{-1}u_{\ell-1})
u_\ell,
$
with $u_k\in U_{-\a_{i_k}}$.
Then
$\dot w u\cdot 1_\Gc=
\dot s_{i_1}u_1
\dot s_{i_2}u_2
\cdots
\dot s_{i_\ell}u_\ell\cdot 1_\Gc.$
Claim $(i)$ follows by induction on $\ell(w)$.
If $y=s_{i_1}s_{i_2}\cdots s_{i_\ell}$
is a reduced expression, then
$s_{i_1}U_{-\a_{i_1}}
s_{i_2}U_{-\a_{i_2}}
\cdots s_{i_\ell}U_{-\a_{i_\ell}}\cdot 1_\Xc=\Xc_y.$
This yields $(iii)$.
By $(iii)$ we have
$wU^-\cdot\Xc_{\leq y}\subset 
U^-wU(\Phi^-\cap w^{-1}\Phi^+)\cdot\Xc_{\leq y}
\subset U^-\cdot\Xc_{\leq w'}.$
We have also $U^-\cdot\Xc_{\leq w'}\subset{}^{w'}\Xc$ by $(iv)$.
This yields $(ii)$.
The first part of Claim $(iv)$ is obvious.
For the second, \cite{K, Lemma 4.5.3} implies that
${}^w\Xc\subset\bigcup_{w'\le w}\Xc^{w'}$ and the opposite
inclusion follows from \cite{K, Lemma 4.5.7}.
The last part of $(iv)$ is proved in \cite{K, Proposition 4.5.9}.
Claim $(v)$ is proved in \cite{KT2, Propositions 1.3.1, 1.3.2}.
\qed

\subhead 2.4\endsubhead
Set ${}^w\Xc^k=U_k^-\setminus{}^w\Xc.$
Let $p^{k_2k_1}\,:\,{}^w\Xc^{k_2}\to{}^w\Xc^{k_1}$, for $k_2\geq k_1$, and 
$p^k\,:\,{}^w\Xc\to{}^w\Xc^{k}$ 
be the natural projections.
The $U_k^-$-action on $\Uc^w$
is free if $w^{-1}(\Theta_k)\subset\Phi^+$.
If $k,k_2,k_1\gg 0$ 
then ${}^w\Xc^k$ is a smooth Noetherian scheme,
$p^{k_2k_1}$ is smooth and affine,
and $p^k$ restricts to a closed immersion
$\Xc_{\leq w}\to{}^w\Xc^{k}$, see \cite{KT2}.
In particular the scheme ${}^w\Xc$ represents the 
pro-object $({}^w\Xc^k,p^k)$.
Note that ${}^w\Xc^k$ may be not separated.

\subhead 2.5\endsubhead
The group $B$ acts on $\uen$ by conjugation. 
The quotient $\Tc=\Gc\times_B\uen$ is a reduced and separated scheme.
Let $\rho\,:\,\Tc\to\Xc$ be the canonical map.
Set $\Vc^w=\rho^{-1}(\Uc^w)$, ${}^w\Tc=\rho^{-1}({}^w\Xc)$,
and $\Tc_y=\rho^{-1}(\Xc_y)$.
Let $[g:x]$ denote the image 
of $(g,x)$ by the projection 
$\Gc\times\uen\to\Tc$.
There is a left action of $\dot W$, $B^-$ on 
$\Tc$ such that $u\cdot[g:x]=[u\cdot g:x]$.

\subhead 2.6\endsubhead
We have $\ad_B\uen_\ell=\uen_\ell$ because 
$\uen_\ell\subset\uen$ is a Lie ideal. Set
${}^w\Tc^{k\ell}=(U_k^-\setminus {}^w\Gc)\times_B(\uen/\uen_\ell).$
If $k\gg 0$ the schemes ${}^w\Tc^{k\ell}$ are Noetherian and smooth.
They may be not separated.
If $k_2\geq k_1\gg 0$ and $\ell_2\geq\ell_1\geq 0$ 
we have the Cartesian square with smooth affine maps
$$\matrix
{}^w\Tc^{k_2\ell_2}
&{\buildrel p^{k_2k_1}_{\ell_2}\over\lra}&
{}^w\Tc^{k_1\ell_2}\cr
{\ss p_{\ell_2\ell_1}^{k_2}}\bda\qquad&
&\qquad\bda{\ss p^{k_1}_{\ell_2\ell_1}}\cr
{}^w\Tc^{k_2\ell_1}
&{\buildrel p^{k_2k_1}_{\ell_1}\over\lra}&
{}^w\Tc^{k_1\ell_1}.
\endmatrix
$$
The scheme ${}^w\Tc$ represents the pro-object
$({}^w\Tc^{k\ell}, p^k_{\ell_2\ell_1},p^{k_2k_1}_\ell),$
because $p^k_{\ell_2\ell_1},p^{k_2k_1}_\ell$ are affine morphisms.
Recall that two pro-objects as above representing the same scheme
${}^w\Tc$ are necessarily isomorphic, see A.5.
We set also
$\Tc^\ell=\Gc\times_B(\uen/\uen_\ell)$.
The map $p_\ell\,:\,\Tc\to\Tc^\ell$, $[g:x]\mapsto[g:x+\uen_\ell]$,
is a vector bundle.
Let $p^k_\ell\,:\,{}^w\Tc\to{}^w\Tc^{k\ell}$ be the obvious projection.

\subhead 2.7\endsubhead
Let $\ben\subset\genu\otimes K$ be the preimage of 
$\benu\subset\genu$ by the natural projection 
$\genu\otimes\CC[[\eps]]\to\genu$.
An Iwahori subalgebra of $\genu\otimes K$ is
a Lie subalgebra of the form $\ad_g\ben$ with $g\in G$.
Recall that, if $g\in\Gu(K)$, then $\ad_g\ben=\ben$ 
if and only if $g\in B$, see \cite{IM, Section 2} for instance.
Let $\Bc$ be the set of all Iwahori subalgebras in $\genu\otimes K$.
Set $\ben_w=\ad_w\ben$.
The group $U(\Phi^+\cap w\Phi^-)$ is identified
with a closed subgroup of $G$.
The map 
$\ad_{U(\Phi^+\cap w\Phi^-)}\ben_w\to\Xc,$
$g\mapsto w\cdot(\ad_{\dot w^{-1}}g)\cdot 1_\Xc$,
is an immersion onto $\Xc_w$.
Hence $\Bc$ is identified with the subset
$\bigcup_{w\in W}\Xc_w\subset\Xc$.

Observe that $\Xc$ is a scheme, hence it is endowed with the corresponding
topology, while $\Bc$ is only a set, which we dentify with a subset in $\Xc$.

\subhead 2.8\endsubhead
Let $\dot\gen=\{(x,\pen)\in\gen\times\Bc\,;\, x\in\pen\}$,
and $p\,:\,\dot\gen\to\gen$ be the 1-st projection.
We set $\Bc_x=p^{-1}(x)$.
Let $\Nc\subset\genu\otimes K$ be the set of {\sl nil-elements}, i.e.
$x\in\Nc$ if $(\ad_x)^n\to 0$ when $n\to\infty$ 
for the topology induced by $K$. 
Recall that $\Bc_x\neq\emptyset$ if $x\in\Nc$,
where $\Nc=\ad_G\uen$, and $\Nc\cap\ben=\uen$ by \cite{KL2, Lemma 2.1}. 
Put $\dot\Nc=p^{-1}(\Nc)$.
Hence
$\dot\Nc\cap(\gen\times\Xc_y)=
\{(\ad_{g\dot y}x,\ad_g\ben_y)\,;\,
g\in U(\Phi^+\cap y\Phi^-), x\in\uen\}.$
Thus the map 
$$\dot\Nc\cap(\gen\times\Xc_y)\to\Tc,\ 
(\ad_{g\dot y}x,\ad_g\ben_y)\mapsto
[g\dot y\cdot 1_\Gc:x]$$ 
is an immersion onto $\Tc_y$.
Therefore $\dot\Nc=\rho^{-1}(\Bc).$ 

Once again, observe that $\Tc$ is a scheme, 
while $\Bc_x$, $\Nc$, $\dot\Nc$, and $\dot\gen$ are only sets,
viewed as subsets in $\Tc$.

\subhead 2.9\endsubhead
Fix a subset $J\subsetneq I$.
Let $\Phi_J\subset\Phi$ be the root system generated
by $\{\a_i\,;\,i\in J\}$, and $\Phi_J^\pm=\Phi_J\cap\Phi^\pm.$ 
Set $U_J=U(\Phi^+\setminus\Phi_J^+)$ and
$\uen_J=\uen(\Phi^+\setminus\Phi_J^+)$.
Let $G_J\subset G$ be the subgroup generated by
$H$, $U(\Phi^+_J)$, and $U(\Phi^-_J)$.
Note that $J$ is a finite root system.
Thus $G_J$ is reductive in the usual sense.
It is also connected and
its derived subgroup is simply-connected by 2.2.
Put $B_J=G_J\ltimes U_J$.
Let $\ben_J\subset\gen$ be the Lie algebra of $B_J$.

The Weyl group of $G_J$ is identified with the subgroup
$W_J\subset W$ generated by $\{s_i\,;\,i\in J\}$.
Let $w_J\in W_J$ be the longest element.
Let ${}^J W$ (resp. $W^J$) be the set of all $y\in W$ such that
$y$ has minimal length in the coset $W_Jy$ (resp. $yW_J$).

One proves as in \cite{K} that
the group $B_J$ acts locally freely on $\Gc$ on the right,
and that the quotient $\Xc^J=\Gc/B_J$ is a separated scheme,
see also \cite{KT3}.

If $w\in W^Jw_J$ then ${}^w\Gc$ is $B_J$-stable and 
we put ${}^w\Xc^J={}^w\Gc/B_J$, ${}^w\Xc^{J,k}=U_k^-\setminus{}^w\Xc^J$. 
Note that $\Xc^J=\bigcup_w{}^w\Xc^J$, where each 
${}^w\Xc^J$ is open, pro-smooth, and of infinite type.
See A.5 for generalities on pro-objects.

Set $\Bc^J$ equal the set of $G$-conjugate of $\ben_J$.
There is a natural embedding $\Bc^J\subset\Xc^J$, see 2.7.
The natural projection $\pi_J\,:\,\Xc\to\Xc^J$ takes
the Lie algebra $\pen\in\Bc$ to the unique element in $\Bc^J$
such that $\pen\subset\pi_J(\pen)$.
The group $G_J$ acts on $\Xc$, $\Xc^J$ on the left.
The map $\pi_J$ commutes with the left $G_J$-action.

Let $\Bcu$ be the variety of all Borel subalgebras in $\genu$,
$\dot\genu=\{(x,\pen)\in\genu\times\Bcu\,;\, x\in\pen\}$,
$p\,:\,\dot\genu\to\genu$ be the 1-st projection,
$\Ncu\subset\genu$ be the nilpotent cone,
$\dot\Ncu=\dot\genu\cap p^{-1}(\Ncu)$,
$\Bcu_x=p^{-1}(x)$ be the Springer fiber,
and $\Zcu=\{(x,\pen;x',\pen')\in\dot\Ncu\times\dot\Ncu\,;\,
x=x'\}$ be the Steinberg variety.

Taking $G_J$ in place of $\Gu$
we get the varieties $\Bc_J$, $\dot\Nc_J$, $\Zc_J$.
If $y\in W_J$ we put 
$$\Oc_{J,y}=G_J\cdot(\ben_{y^{-1}},\ben),\quad
\Zc_{J,y}=\{(x,\pen,\pen')\in\Zc_J\,;\,(\pen,\pen')\in\Oc_{J,y}\}.$$
The set $\Zc_{J,\leq y}=\bigcup_{y'\leq y}\Zc_{J,\leq y'}$
is closed in $\Zc_J$. 
If $y\in{}^J W$ we put ${}^J\Bc_y=\ad_{B_J}(\ben_y)\subset\Bc$.
We write ${}^J\Bc$ for ${}^J\Bc_1$.
The map ${}^J\Bc\to\Bc_J$, $\pen\mapsto\pen/\uen_J$,
is an isomorphism of $G_J$-varieties.
Hereafter we identify ${}^J\Bc$ and $\Bc_J$.

From now on we write $\flat$ for $\Iu$, i.e.
$\uen_\flat=\uen_\Iu$, $\ben_\flat=\ben_\Iu$, 
$B_\flat=B_\Iu$, $W^\flat=W^\Iu$, etc.
In particular the map ${}^\flat\Bc\to\Bcu$, $\pen\mapsto\pen/\uen_\flat$
is an isomorphism of $\Gu$-varieties.

Recall that $\ben_\flat/\uen_\flat\simeq\genu$.
For any Iwahori Lie subalgebra $\pen\subset\genu\otimes K$
the subspace $\pen^!=(\pen\cap\ben_\flat)/(\pen\cap\uen_\flat)\subset\genu$
is a Lie subalgebra.
Let $!\,:\,\Bc\to\Bcu$ be the map 
taking $\pen$ to $\pen^!.$

\subhead 2.10\endsubhead
Set $\Pc^J=\Gc/U_J$.
It is a scheme because $U_J$ acts locally freely on $\Gc$.
The group $G_J$ acts on $\Pc^J$ on the right 
because $G_J$ normalizes $U_J$.
This action is locally free.
Hence, the corresponding quotient map 
$\Pc^J\to\Xc^J$ is a $G_J$-torsor.
In particular $\Pc^J$ is a separated scheme.

Assume that $k\gg 0$, $w\in W^Jw_J$, and $\uen_\ell=\eps^\ell\cdot\uen_J$.
We set ${}^w\Pc^{J,k}=U_k^-\setminus{}^w\Gc/U_J$.
It is a scheme.
The group $G_J$ acts on ${}^w\Pc^{J,k}$ on the left and on the right,
because $G_J$ normalizes $U_J$ and $U^-_k$.
Both actions commute.
The quotient ${}^w\Pc^{J,k}\to{}^w\Xc^{J,k}$ is a $G_J$-torsor. 

\subhead 2.11\endsubhead
The group $G\times\CC^\times$ acts on $\gen$ by 
$\ad_{(g,\zeta)}x=\zeta^{-1}\ad_gx$, and on
$G$ by $\ad_{(g,\zeta)}h=\ad_gh$.
To avoid confusion between the different one parameter subgroups
we write $\CC^\times_q$ for the subgroup
$\{1\}\times\CC^\times\subset G\times\CC^\times$.
Set $A=H\times\CC_q^\times$.
The group $A$ acts on $\Tc$ by
$(s,\zeta)\cdot [g:x]=[s\cdot g:\zeta^{-1}\ad_sx]$.

\vskip3mm

\noindent{\bf 2.12. Convention.}
From now on we use the following conventions :
$\su\in\Hu$, $s'\in H'$ with component in $\Hu$ equal to $\su$,
$s=(s',\tau)\in H$ with $\tau\in\CC^\times_\delta$,
and $a=(s,\zeta)\in A$ with $\zeta\in\CC^\times_q$.
The corresponding fixed point sets are denoted
as usual (i.e. $G^s$, $\Bc^s$, $\gen^s$, $\Tc^a$, etc).
We also assume that $\tau$ is not a root of unity,
although a few statements still hold without this assumption.
Finally, an element in $G$ is called 
{\sl semisimple} if it belongs to $\ad_GH$.

\subhead 2.13\endsubhead
The action of $s$ on $\Bc\subset\Xc$ takes $\pen$ to $\ad_s\pen$.
The scheme ${}^w\Xc^k$ is endowed with the unique 
$H$-action such that the map
$p^k\,:\,{}^w\Xc\to{}^w\Xc^k$ is $H$-equivariant.
Similarly ${}^w\Tc^{k\ell}$ is endowed with the unique $A$-action
such that the map ${}^w\Tc\to{}^w\Tc^{k\ell}$ is $A$-equivariant.
We write ${}^w\Xc^s$ for $({}^w\Xc)^s$,
and ${}^w\Tc^a$ for $({}^w\Tc)^a$.
Finally, let $\{\Xc_i^{s}\,;\,i\in\Xi_s\}$ 
be the set of connected components of $\Xc^s$.

\proclaim{Lemma}

(i)
$G^s$ is connected, reductive, with Cartan subgroup $H$ and Weyl group $W_s$.

(ii)
$\Xc^s=\Bc^s$;
it is a smooth separated
$G^s$-scheme locally of finite type.
Each component $\Xc_i^s\subset\Xc^s$
is isomorphic to the flag variety of $G^s$
as a $G^s$-variety.

(iii)
$p^k$ restricts to a closed immersion
${}^w\Xc^{s}\to{}^w\Xc^k$ onto $({}^w\Xc^k)^{s}$ if $k\gg 0$.

(iv)
$\Tc^a=\dot\Nc^a;$ 
it is a smooth separated $G^s\times\CC^\times_q$-scheme
locally of finite type.

(v)
$\Tc^a$ is a vector bundle of finite rank over $\Xc^s$.

(vi)
$p^k_\ell$ restricts to a closed immersion 
${}^w\Tc^a\to{}^w\Tc^{k\ell}$ onto $({}^w\Tc^{k\ell})^a$
if $k,\ell\gg 0$.
\endproclaim

\noindent{\sl Proof:}
For any $w\in W$ we fix a representative $n_w$ of $w$ in $N_G(H)$.
The Bruhat decomposition for $G'$ implies that $G^{s}$ 
is generated by $H$, the $U_\a$'s with $\a$ a real root such that
$\gen_\a\subset\gen^s$, 
and the $n_w$'s such that $w(s)=s$.
Since $\tau$ is not a root of unity,
the number of $\a$ and $n_w$ as above is finite.
Hence $G^s\subset G$ is a linear closed subgroup.
It is reductive because $U_\a\subset G^s$ if and only if
$U_{-\a}\subset G^s$ 
(if there was a nontrivial unipotent radical $R$,
then $R$ would be normalized by $H$ and therefore contain some
$U_\a$; in turn, $R$ would be normalized by $U_{-\a}$,
yielding non-unipotent elements in $R$).
In particular $G^s\subset G_{KM}\rtimes\CC_\delta^\times$.
Hence there are $g\in G_{KM}$, $J\subsetneq I$ such that
$\ad_gG^s\subseteq G_J$, by \cite{KcP, Proposition 3.1}.
Then $\ad_gG^s=G_J^{\ad_gs}$,
and $G^s$ is connected by \cite{SS, Theorem 3.9},
because $\ad_gs$ is a semisimple element of $G_J$ 
and $G_J$ is a connected reductive group whose
derived subgroup is simply connected, see 2.9.
Claim $(i)$ is proved. 

Consider the open covering
$\bigcup_{w'\leq w}(\Uc^{w'})^s$ of ${}^w\Xc^s$.
We have $(\Uc^{w'})^s=w'V_{w'}\cdot 1_\Xc$ for a 
closed subscheme $V_{w'}\subset U^-$ of finite type.
In particular $V_{w'}\subset G$. 
Thus $(\Uc^{w'})^s\subset\Bc$ and $\Xc^s=\Bc^s$.
Moreover ${}^w\Xc^s$ is smooth and of finite type
because $V_{w'}$ is of finite type.

If $k\gg 0$ the $U^-_k$-action on $\Uc^w$ is free,
i.e. $\Uc^w\simeq U^-_k\times p^k(\Uc^w)$.
The isomorphism commutes with the $s$-action.
Hence $p^k(\Uc^w)^s=p^k((\Uc^w)^s)$.
Thus $({}^w\Xc^k)^s=p^k({}^w\Xc^s)$.
The map $p^k$ restricts to an immersion ${}^w\Xc^s\to{}^w\Xc^k$,
because ${}^w\Xc^s$ is of finite type.
Let us prove that $p^k({}^w\Xc^s)$ is closed in ${}^w\Xc^k$.
Since $G^s$ and ${}^w\Xc^s$ are of finite type
we may fix an element $w'\geq w$
such that $G^s\cdot{}^w\Xc^s\subset{}^{w'}\Xc$ 
and $p^k$ yields an immersion of $G^s\cdot{}^w\Xc^s$ into ${}^{w'}\Xc^k.$ 
The isomorphism $\Uc^{w'}\simeq U^-_k\times p^k(\Uc^{w'})$
takes $G^s\cdot{}^w\Xc^s$ onto $U_k^-\times p^k(G^s\cdot{}^w\Xc^s)$.
Hence $U_k^-\times p^k(G^s\cdot{}^w\Xc^s)$ is closed in $\Uc^{w'}$,
because $\Uc^{w'}$ is separated and $G^s\cdot{}^w\Xc^s$ is proper. 
Hence $p^k(G^s\cdot{}^w\Xc^s)$ is closed in ${}^{w'}\Xc^k$,
because $p^k$ is a quotient morphism.
Therefore $p^k({}^w\Xc^s)={}^w\Xc^k\cap p^k(G^s\cdot{}^w\Xc^s)$ 
is closed in ${}^w\Xc^k$.
Thus $(ii)$, $(iii)$ are proved,
and $(iv),(vi)$ follow, using $(v)$.

Finally, to prove $(v)$ it is enough to observe that
the projection $\Tc^a\to\Xc^s$ is the restriction
of the first projection $\Gc\times_B\uen^a\to\Gc/B$ to $\Tc^a$,
and that $\uen^a$ is finite dimensional 
because $\tau$ is not a root of unity, see Lemma 2.14 below.
\qed

\vskip3mm

The projection $\dot\Nc^a\to\Xc^s$ is surjective.
For a future use, let $\dot\Nc_i^a\subset\dot\Nc^a$ be the unique 
component above $\Xc_i^s$ for each $i\in\Xi_s$.

\vskip3mm

\noindent{\bf Remark.}
It is essential to consider the group $G$ (or $G'$)
rather than $\Gu(K)$. 
Otherwise Lemma 2.13.$(i)$ does not hold : 
it was noticed by D. Vogan that if $\Gu$ is of type $D_4$,
$\tau^{1/2}$ is a square root of $\tau$, and
$\su=(\a_1^\vee\otimes(-1))(\a_2^\vee\otimes\tau^{1/2})
(\a_3^\vee\otimes(-1))(\a_4^\vee\otimes(-\tau^{1/2}))
\in\Yu^\vee\otimes_\ZZ\CC^\times$,
then the reductive group $\Gu(K)^s$ is not connected, 
see \cite{BEG1, Section 5}.
However $G^s$ and $(G')^s$ are connected by Lemma 2.13.$(i)$.

\subhead 2.14\endsubhead
Recall that $\tau$ is not a root of unity by Convention 2.12. 

\proclaim{Lemma}
$(\genu\otimes K)^a$ is finite dimensional over $\CC$, 
$\Nc^a\subseteq(\genu\otimes K)^a$ is closed.
If $\zeta$ is not a root of unity there are 
$d$, $\g_d\in\Gu(K_d)$, $S\subset\Hu\times\CC^\times_q$ finite,
$r\in\QQ$,
such that 
$\genu^S\otimes\eps^r=\ad_{\g_d}((\genu\otimes K)^a)$
and $\Gu^S\subseteq\ad_{\g_d}(\Gu(K)^s)$.

(i) 
If $\zeta^{\ZZ\setminus\{0\}}\cap\tau^\ZZ=\emptyset$ then $r=0$, 
$(\genu\otimes K)^a=\Nc^a$ consists of nilpotent elements,
and there are a finite number of $G^s$-orbits in $\Nc^a$.

(ii)
If $\zeta^m=\tau^k$ with $m\in\ZZ_{>0}$, $k\in\ZZ_{<0}$, $(m,k)=1$,
then $r=k/m$, $\Nc^a$ consists of nilpotent elements, 
and there are a finite number of $G^s$-orbits in $\Nc^a$.

(iii)
If $\zeta^m=\tau^k$ with $m,k\in\ZZ_{>0}$, $(m,k)=1$, then $r=k/m$ 
and $\Nc^a=(\genu\otimes K)^a$.
\endproclaim

\noindent{\sl Proof:}
In the proof we use the following notations :
$\gen_i=\{x\in\genu\otimes K\,;\,\ad_sx=\zeta^ix\}$, 
$\ad^{(d)}$ is the adjoint action of $G_d$ on 
$\genu\otimes K_d$ and $\Gu(K_d)$, see 2.2,
and $\mu_d=\exp(2i\pi/d)$.
Finally, given $\g\in\Yu^\vee$ 
we may view $\g$ as a group homomorphism
$K_d^\times\to\Hu(K_d)$ and write $\g(z)$ for the image
of $z\in K_d^\times$.

The space $\gen_i$ is finite dimensional for each $i$,
because $\tau$ is not a root of unity.
Hence, $(\genu\otimes K)^a$ is finite dimensional and
$\Nc^a\subseteq(\genu\otimes K)^a$ is a closed subset.
Observe that if $(i)$, $(ii)$ or $(iii)$ holds 
then $\zeta$ is not a root of unity.

Assume that $\zeta^m\notin\tau^\ZZ$ for all $m\neq 0$.
Let $\Phi_{s,i}\subset\Phiu$ be the set of roots 
$\a$ such that $(\genu_\a\otimes K)\cap\gen_i\neq 0.$
Put $\Phi_s=\bigcup_{i\in\ZZ}\Phi_{s,i}$.
The set $\Phi_s$ is closed and symmetric,
i.e. $\a,\b\in\Phi_s,\,\a+\b\in\Phiu\Rightarrow\a+\b\in\Phi_s$,
and $\Phi_s=-\Phi_s$.
Thus it is a root system, see \cite{B, Ch. VI, \S 1, Proposition 23}. 
Fix a basis $\b_1,...\b_r$ of $\Phi_s$.
Note that $\Phi_{s,i}\cap\Phi_{s,j}=\emptyset$ if $i\neq j$
because $\zeta^{i-j}\notin\tau^\ZZ$ if $i\neq j$.
Hence there are unique $n_1,...n_r\in\ZZ$ and $i_1,...i_r\in\ZZ$
such that $\genu_{\b_t}\otimes\eps^{n_t}\subset\gen_{i_t}$ 
for all $t=1,2,...,r$.
Fix $d\in\ZZ_{>0}$, $\g\in\Yu^\vee$ 
such that $(\b_t:\g)=-dn_t$ for all $t$.
Fix a $d$-th root of $\tau$.
Set 
$$S_i=\{(\g(\mu_d),1),(\su\g(\tau^{1/d})^{-1},\zeta^i)\}
\subset\Hu\times\CC^\times_q.$$
Put $\g_d=\g(\eps^{1/d})\in\Gu(K_d)$.
Then $\ad^{(d)}_{\g_d}\gen_i=\genu^{S_i}$, because
$$\gen_i=\{x\in\genu\otimes K_d\,;\,
\ad^{(d)}_{(1,\mu_d)}x=\zeta^{-i}\ad^{(d)}_{(\su,\tau^{1/d})}x=x\}.$$
Let $L\subset\Gu(K_d)$ be the subgroup consisting of the
elements fixed by the automorphisms
$\ad^{(d)}_{(1,\mu_d)}$ and $\ad^{(d)}_{(\su,\tau^{1/d})}.$
The Bruhat decomposition for $\Gu(K_d)$
implies that $L$ is generated by the group $\Hu$, 
the root subgroups in the `loop group' $\Gu(K_d)$
associated to the real affine root whose root subspace in
$\genu\otimes K_d$ is contained in $\gen_0$, 
and the elements in $N_{G_d}(\Hu\times\CC^\times_\delta)$ fixing 
$(1,\mu_d)$ and $(\su,\tau^{1/d}).$
Similarly $\Gu(K)^s$ is generated by $\Hu$, the $U_\a$'s
with $\a\in\Phi$ a real root such that $\gen_\a\subset(\genu\otimes K)^s$,
and the $n_w$'s such that $w(s)=s$.
Hence $\Gu(K)^s=L.$ 
Thus $\ad^{(d)}_{\g_d}(\Gu(K)^s)=\Gu(K_d)^{S_0}$.
In particular $\ad^{(d)}_{\g_d}(\Gu(K)^s)$ contains $\Gu^{S_0}$.
Further, we have $\genu^{S_i}=\Ncu^{S_i}$ 
and there is a finite number of 
$\Gu^{S_0}$-orbits in $\Ncu^{S_i}$ for each $i\neq 0$
by \cite{KL1, 5.4.$(c)$},
because $\zeta$ is not a root of unity.
Finally observe that $\Gu^{S_0}=\Gu^{S_1}$ because $\CC^\times_q$ acts trivially
on $\Gu$ by definition.
Claim $(i)$ is proved.

Assume that 
$\zeta^m=\tau^k$ with $m\in\ZZ_{>0}$, $k\in\ZZ$, and $(m,k)=1$.
Fix a $m$-th root of $\tau$ such that $\zeta=\tau^{k/m}$.
We have 
$$\gen_i=
\{x\cdot\eps^{ik/m}\,;\,x\in\genu\otimes K_m,\,\ad^{(m)}_{(\su,\tau^{1/m})}x=x\}
\cap(\genu\otimes K).$$
The set of roots $\a\in\Phiu$ such that 
$$\{x\in\genu_\a\otimes K_m
\,;\,\ad^{(m)}_{(\su,\tau^{1/m})}x=x\}\neq \{0\}$$
is a root system.
Fix a basis $\b_1,\b_2,...\b_r$ of this root system.
There are unique integers $n_1,...n_r$ such that
$\genu_{\b_t}\otimes\eps^{n_t/m}$ is fixed by
$\ad^{(m)}_{(\su,\tau^{1/m})}$ for each $t=1,2,...r$,
because $\tau$ is not a root of unity.
Fix $e\in\ZZ_{>0}$, $\g\in\Yu^\vee$ 
such that $(\b_t:\g)=-en_t$ for all $t$.
Then
$$\ad^{(d)}_{\g_d}(\gen_{\b_t}\otimes \eps^{n_t/m})\subseteq\genu$$
for each $t$, where $d=em$ and $\g_d=\g(\eps^{1/d})$.
Hence
$$\ad^{(d)}_{\g_d}\{x\in\genu\otimes K_m\,;
\,\ad^{(m)}_{(\su,\tau^{1/m})}x=x\}=\genu^{S'},$$
where $S'=\{\g(\mu_e),\su\g(\tau^{1/d})^{-1}\}\subset\Hu.$
Thus,
$$\matrix
\ad^{(d)}_{\g_d}(\gen_i)&=
(\genu^{S'}\otimes\eps^{ik/m})\cap
\ad^{(d)}_{\g_d}\{x\in\genu\otimes K_d\,;\,\ad^{(d)}_{(1,\mu_d)}x=x\},\hfill\cr
&=(\genu^{S'}\otimes\eps^{ik/m})\cap
\{x\in\genu\otimes K_d\,;\,\ad^{(d)}_{(\g(\mu_d)^{-1},\mu_d)}x=x\}.\hfill\cr
\endmatrix$$
Set $S=S'\cup\{(\g(\mu_d),(\mu_{dm})^k)\}\subset\Hu\times\CC^\times_q$.
Then $\ad^{(d)}_{\g_d}\gen_1=\genu^S\otimes\eps^{k/m}$ and
$\ad^{(d)}_{\g_d}\gen_0=\genu^{S}.$
Using the Bruhat decomposition for $\Gu(K_d)$ and $\Gu(K)$
we get also $\Gu^S\subset\ad^{(d)}_{\g_d}(\Gu(K)^s)$.
The last statement in $(ii)$ follows from \cite{KL1, 5.4.$(c)$}. 
Note that the condition in $(ii)$ says that $r<0$.
Thus $\Nc^a$ consists of nilpotent elements, because
any element in $\Nc\cap(\genu\otimes\eps^r)$ is nilpotent if $r<0$.
\qed

\vskip3mm

\proclaim{Definition}
Assume that $\tau,\zeta$ are not roots of unity.
If $(i)$ or $(ii)$ holds
we say that the pair $(\tau,\zeta)$ is regular.
If $(i)$ holds we say that $(\tau,\zeta)$ is generic. 
\endproclaim

\head 3. Reminder on double affine Hecke algebras\endhead
\subhead 3.1\endsubhead
Let $\Hb$ be the unital associative 
$\CC[q,q^{-1}, t, t^{-1}]$-algebra generated by 
$\{t_i, x_\l\,;\, i\in I, \l\in X\}$
modulo the following defining relations
$$\matrix
x_\delta=t\hfill&\cr
x_\l x_\mu=x_{\l+\mu}\hfill&\cr
(t_i-q)(t_i+1)=0\hfill&\cr
t_it_jt_i\cdots=t_jt_it_j\cdots\hfill&\text{if}\ i\neq j\ 
\text{($m_{ij}$\ factors\ in\ both\ products)}\hfill\cr
t_ix_\l-x_\l t_i=0\hfill&\text{if}\ (\l:\a_i^\vee)=0\hfill\cr
t_ix_\l-x_{s_i(\l)}t_i=(q-1)x_\l\hfill&
\text{if}\ (\l:\a_i^\vee)=1\hfill\cr
\endmatrix$$
for all $i,j\in I$, $\l,\mu\in X$.
The complexified Grothendieck ring $\Rb_H$ is naturally identified with
the $\CC$-linear span of $\{x_\l\,;\,\l\in X\}$,
and $\Rb_A$ is identified with $\Rb_H[q,q^{-1}]$.
In particular it is viewed as a subring of $\Hb$.

For each reduced decomposition $y=s_{i_1}s_{i_2}\cdots s_{i_\ell}\in W$
we set $t_y=t_{i_1}t_{i_2}\cdots t_{i_\ell}.$
Then $(t_y\,;\,y\in W)\subset\Hb$ is a $\Rb_A$-basis.

Let $\IM\,:\,\Hb\to\Hb$ be the 
$\CC[q,q^{-1}]$-linear involution such that
$\IM(t_i)=-qt_i^{-1}$ and $\IM(x_\l)=x_\l^{-1}$.

For any $\zeta,\tau\in\CC^\times$ we set 
$\Rb_{\tau,\zeta}=\Rb_A/(q-\zeta,t-\tau)$,
and $\Hb_{\tau,\zeta}=\Hb/(q-\zeta,t-\tau)$.

On $\Hb_{\tau,\zeta}$ we introduce a filtration $\Hb_{\leq w}$ , $w\in W$,
setting $\Hb_{\leq w}$ to be the span of the basis elements 
$\{t_yx_\l\,;\,\l\in X,\,y\leq w\}$.

\vskip3mm

\noindent{\bf Remarks.}
$(i)$
One important feature of the ring $\Hb$ is that it is not
of finite type over its center.
This is the reason why in the proof of Theorem 7.6 below we use
the support in $\Spec\Rb_A$ of some $\Hb$-module
rather than the `central character'.

$(ii)$
Note that $\Hb$ differs slightly from the
double affine Hecke algebra used by Cherednik. 
The latter is the $\CC[q,q^{-1}]$-algebra $\Hb^\ch$ 
generated by $\{t_\o,t_i,x_\l\,;\,
\o\in\O,i\in I,\l\in\Xu\oplus\ZZ\delta\}$
with relations analogous to the previous one,
see also 8.3.

\subhead 3.2\endsubhead
For any non empty subset $J\subset I$ let $\Hb_J$ be the unital associative
$\CC[q,q^{-1}]$-algebra generated by 
$\{t_i, x_\l\,;\,i\in J,\l\in X'\}$ 
modulo the same relations as above, with $X'$ as in 2.2. 
If $J\neq I$ the
$\CC[q,q^{-1}]$-algebra $\Hb_J$ is the affine Hecke algebra associated to
the group $G_J$.
It is well-known that 
$\Hb_J$ is isomorphic to the 
$\CC[q,q^{-1}]$-subalgebra of $\Hb$ generated by 
$\{t_i,x_\l\,;\, i\in J,\l\in X'\}$.
The proof of Theorem 4.9 is based on the following easy lemma.

\proclaim{Lemma}
Assume that $|I|>2$.
Given any ring $\Ab$, 
the restriction map yields a bijection from the set of ring homomorphisms
$\Psi\,:\,\Hb\to\Ab$ onto the set of famillies $(\Psi_J)$
of ring homomorphisms $\Psi_J\,:\,\Hb_J\to\Ab$ such that
$J\neq I$ and $(\Psi_{J_2})|_{J_1}=\Psi_{J_1}$ whenever
$J_1\subset J_2$.
\endproclaim

We will write $\Hbu$ for $\Hb_\Iu$,
and $\Hbu_\zeta$ for $\Hbu/(q-\zeta)$.
The induction functor 
$\Lb\to\Hb_{\zeta,\tau}\otimes_{\Hbu_\zeta}\Lb$
from $\Hbu_\zeta$-modules to $\Hb_{\zeta,\tau}$-modules
is exact. The restriction functor is also exact.

\head 4. The convolution algebra\endhead
The affine analogue, $\Zc$, of the Steinberg variety is introduced in 4.1-3.
Fix an element $a\in A$ as in 2.12.
The scheme $\Zc^a$ is the disjoint union of the connected components
$\Zc^a_{ij}=\Zc^a\cap(\dot\Nc^a_i\times\dot\Nc^a_j)$
with $i,j\in\Xi_s$, see 2.13.
Each component is a scheme of finite type.
In 4.4 we introduce a subspace
$\Kb_a\subseteq\prod_j\bigoplus_i\Kb^A(\Zc^a_{ij})_a,$
where $\Kb^A(\Zc^a_{ij})$ is the complexified Grothendieck group,
see 1.4.
The convolution product on $\Kb_a$ is defined in 4.4.
A ring homomorphism $\Psi_a\,:\,\Hb\to\Kb_a$ is given in 4.9.
It is surjective after a suitable completion of $\Hb$, 
see the proof of Theorem 4.9.$(iv)$.
The main ingredient in the proof of Theorem 4.9 is the ring homomorphism 
$\Theta_J\,:\,\Kb^{G_J\times\CC^\times_q}(\Zc_J)\to\Kb_a$.
It is the composition of a chain of maps
$$\Kb^{G_J\times\CC^\times_q}(\Zc_J)\to
\Kb^{A}(\Zc^0_{\leq w_J})\to
\Kb^{A}(\Zc^a_{\leq w_J})_a\to\Kb_a,$$
where
$\Zc^0_{\leq w_J}=\Pc^J\times_{G_J}\Zc_J.$
The first map is induction. See (A.4.2). 
The convolution product on $\Kb^A(\Zc^0_{\leq w_J})$ is defined in 4.5.
The second map is the concentration map and the Thom isomorphism, 
because $\Zc_{\leq w_J}$ is a vector bundle over $\Zc^0_{\leq w_J}$,
see 4.5, 4.7 and the appendix.
The compatibility between the maps $\Theta_J$, see Step 3 in 1.1, 
is proved by computing explicitely the images
$D_{i,a}$, $D_{\l,a}$ of certain elements
in $\Kb^{G_J\times\CC^\times_q}(\Zc_J)$.
This is done in 4.6-7.
Further results on the ring $\Kb_{a,\leq 1}$ are given in 4.8
for a future use.

Once again, Kashiwara's manifold $\Xc^J$ seems to be more appropriate 
than the ind-scheme $\Bc^J$, in particular to define $\Kb^A(\Xc^J)$ and the 
concentration map $\Kb^A(\Xc^J)_a\to\Kb^A((\Xc^J)^a)_a,$ 
see Lemma 4.5.$(iv)$.

\subhead 4.1\endsubhead
Given $y\in W$ we set
$$\Oc_y=\{(\pi(g\cdot\dot y^{-1}),\pi(g))\,;\,g\in\Gc\},
\and \Oc_{\leq y}=\bigcup_{y'\leq y}\Oc_{y'},$$
where $\pi$ is as in 2.3.
For each $w\in W$ we set also 
${}^w\Oc_y=\Oc_y\cap(\Xc\times{}^w\Xc)$,
and
${}^w\Oc_y^k=(p^k\times p^k)({}^w\Oc_y)$.
Let $q_i, q_{ij}$ be as in 1.4.

\proclaim{Lemma}
Assume that $w'\succeq(w,y^{-1})$, $y''\succeq(y',y)$, and 
$k\gg 0$.

(i)
$\Oc_{\leq y}$ is closed in $\Xc^2$,
$\Oc_y$ is open in $\Oc_{\leq y}$.

(ii)
${}^w\Oc_{\leq y}\subset{}^{w'}\Xc\times{}^w\Xc$,
and ${}^w\Oc_y^k$ is locally closed in ${}^{w'}\Xc^k\times{}^w\Xc^k$.

(iii)
$q_{13}(q_{12}^{-1}(\Oc_{\leq y'})\cap 
q_{23}^{-1}(\Oc_{\leq y}))\subset\Oc_{\leq y''}.$

(iv) 
$q_2|_{\Oc_y}$ is a vector bundle of rank $\ell(y)$,
$q_2|_{\Oc_{\leq y}}$ is a locally trivial fibration 
with fibers $\simeq\Xc_{\leq y^{-1}}$.
\endproclaim

\noindent{\sl Proof:}
The set $\Oc_y$ is preserved by the diagonal action of
$\dot W$, $U^-$ on $\Xc^2$.
Hence,
$$\Oc_y\cap(\Xc\times\Uc^w)=
wU^-\cdot\{(\pi(1_\Gc\cdot u\dot y^{-1}), 1_\Xc)\,;\,u\in U\}.$$
Thus, using Lemma 2.3.$(i)$ we get
$$\Oc_y\cap(\Xc\times\Uc^w)=
wU^-\cdot(\Xc_{y^{-1}}\times\{1_\Xc\})
\simeq\Xc_{y^{-1}}\times U^-.\leqno(4.1.1)$$
Claim $(i)$ follows from (4.1.1)
because the sets $\Uc^w$ with $w\in W$ form an open covering of $\Xc$.

By (4.1.1) we have 
$$(\Oc_{\leq y'}\times\Uc^w)\cap(\Xc\times\Oc_{\leq y})\subset
(\Oc_{\leq y'}\times\Xc)\cap 
(\Xc\times wU^-\cdot(\Xc_{\leq y^{-1}}\times\{1_\Xc\})).$$
By Lemma 2.3.$(iii)$ we have 
$$(\Oc_{\leq y'}\times\Xc)\cap
(\Xc\times\Xc_{\leq y^{-1}}\times\{1_\Xc\})\subset
q_1^{-1}(\Xc_{\leq {y''}^{-1}})\times\{1_\Xc\}.$$
Hence
$q_{13}((\Oc_{\leq y'}\times\Uc^w)\cap(\Xc\times\Oc_{\leq y}))
\subset\Oc_{\leq y''}$
because $\Oc_{\leq y'}$ and $\Oc_{\leq y''}$ 
are preserved by the diagonal 
action of $U^-$ and $\dot W$.
Claim $(iii)$ follows.

By Lemma 2.3.$(ii)$ and (4.1.1) we have 
$\Oc_{\leq y}\cap(\Xc\times\Uc^w)\subset{}^{w'}\Xc\times\Uc^w$.
The 1-st part of $(ii)$ is proved. 

The automorphism $\varphi$ of $\Xc\times\Uc^w$ such that
$\varphi(\pi(g),\pi(wu))=(\pi(u^{-1}w^{-1}g),\pi(wu))$
for all $u\in U^-$, $g\in\Gc$,
takes $\Oc_y\cap(\Xc\times\Uc^w)$
onto $\Xc_{y^{-1}}\times\Uc^w.$
Claim $(iv)$ follows.

By $(iv)$ 
the restriction of $q_2$
to ${}^w\Oc_{y}$ is a morphism of finite type.
Hence, if $k\gg 0$ the map $p^{k}\times\Id$ 
restricts to an immersion of ${}^w\Xc$-schemes
${}^w\Oc_{y}\to{}^{w'}\Xc^{k}\times{}^w\Xc$.
The $U_k^-$-action on ${}^{w'}\Xc^{k}\times{}^w\Xc$
on the right factor preserves
$(p^k\times\Id)({}^w\Oc_y)$
because ${}^w\Oc_y$ is preserved by the diagonal $U^-$-action.
The quotient is a locally closed subset of
${}^{w'}\Xc^k\times{}^w\Xc^k$ isomorphic to
$\Xc_{y^{-1}}\times(U_k^-\setminus wU^-w^{-1})$.
The 2-nd part of $(ii)$ is proved.
\qed

\vskip3mm

\noindent{\bf Remark.}
Note that $\Oc_y$, $\Oc_{\leq y}$ and ${}^w\Oc_y^k$
are endowed with the reduced scheme structures induced from 
$\Xc^2$ and ${}^{w'}\Xc^k\times{}^w\Xc^k$ respectively.

\subhead 4.2\endsubhead
To simplify the notations, for any $w\in W$ it is convenient to write
$w(x)$, $\dot w(x)$ for the adjoint action
$\ad_wx$, $\ad_{\dot w}x$.
See also 2.3.
Set
$$\Zc_y=\{([g\cdot\dot y^{-1}:\dot y(x)],[g:x])\,;\,
g\in\Gc,\,x\in\uen\cap y^{-1}(\uen)\}.$$
Set also $\Zc_{\leq y}=\bigcup_{y'\leq y}\Zc_{y'}$, $\Zc=\bigcup_y\Zc_y$,
${}^w\Zc_y=\Zc_y\cap(\Tc\times{}^w\Tc)$,
${}^w\Zc_{y}^{k\ell}=(p^k_\ell\times p^k_\ell)({}^w\Zc_{y})$, etc.

\proclaim{Lemma}
Let $y,y',y'',w,w'$ be as in Lemma 4.1.

(i) 
$\Zc_{\leq y}$ is closed in $\Tc^2$,
$\Zc_y$ is open in $\Zc_{\leq y}$, 
the second projection $\Zc_{\leq y}\to\Tc$ is proper,
the natural projection $\Zc_y\to\Oc_y$ is a vector bundle.

(ii)
${}^w\Zc_y\subset{}^{w'}\Tc\times{}^w\Tc$,
${}^w\Zc^{k\ell}_y$ is locally closed in 
${}^{w'}\Tc^{k\ell}\times{}^w\Tc^{k\ell}$ if $k,\ell\gg 0$,
the projection ${}^w\Zc^{k\ell}_y\to{}^w\Oc^k_y$ 
is a vector bundle with fiber 
$(\uen\cap y^{-1}(\uen))/(\uen_\ell\cap y^{-1}(\uen_\ell))$.

(iii)
$q_{13}$ restricts to a proper map
$q_{12}^{-1}({}^{w'}\Zc_{\leq y'})\cap
q_{23}^{-1}({}^{w}\Zc_{\leq y})\to{}^{w}\Zc_{\leq y''}.$
\endproclaim

\noindent{\sl Proof:}
The set $\Zc_y$ is preserved by the 
diagonal action of $\dot W$, $U^-$ on $\Tc^2$.
Hence Lemma 2.3.$(i)$ gives
$$\matrix
\Zc_y\cap(\Tc\times\Vc^w)
&=wU^-\cdot\{([1_\Gc\cdot u\dot y^{-1}:\dot y(\ad_u^{-1}x)],
[1_\Gc:x])\,;\,u,x\}\hfill\cr
&=
wU^-\cdot\{([\dot y^{-1}u_2\cdot 1_\Gc:\ad_{u_2}^{-1}\dot y(x)],
[1_\Gc:x])\,;\,u_2,x\},\hfill
\endmatrix$$
where $u\in U,$ $x\in\uen\cap\ad_uy^{-1}(\uen)$,
and $u_2\in U(\Phi^-\cap y\Phi^+)$, $x\in\uen\cap y^{-1}(\ad_{u_2}\uen)$
respectively.
Hence the natural projection $\Zc_y\to\Oc_y$ is a vector bundle 
whith fibers isomorphic to $\uen\cap y^{-1}(\uen).$
The automorphism $\varphi$ of $\Tc\times\Vc^w$ such that
$$(t,[\dot w u_1\cdot 1_\Gc:x])\mapsto
(u_1^{-1}\dot w^{-1}\cdot t,[\dot w u_1\cdot 1_\Gc:x]),
\quad\forall u_1\in U^-,$$
takes $\Zc_y\cap(\Tc\times\Vc^w)$ to
$$\Zc^\varphi_y=\{([\dot y^{-1}u_2\cdot 1_\Gc:\ad_{u_2}^{-1}\dot y(x)],
[\dot wu_1\cdot 1_\Gc:x])\,;\,u_1,u_2,x\}$$
with $u_1,u_2,x$ as above.
Hence $\Zc_y^\varphi$ is closed in $\Vc^{y^{-1}}\times\Vc^w$,
because $\uen\cap y^{-1}(\ad_{u_2}\uen)$ is closed in $\uen$.
Thus $\Zc_{\le y}$ is closed in $\Tc^2$
and $\Zc_y$ is open in $\Zc_{\le y}$,
because the sets $\Vc^w$ with $w\in W$ form an open covering of $\Tc$.
The map $\rho\times\Id$ restricts to a closed immersion 
$\Zc^\varphi_{\leq y}\to\Xc_{\leq y^{-1}}\times\Vc^w$,
where $\Zc^\varphi_{\leq y}=\bigcup_{y'\leq y}\Zc^\varphi_{y'}$.
Hence the second projection 
$\Zc^\varphi_{\leq y}\to\Vc^w$ is proper.
Thus the second projection 
$\Zc_{\leq y}\to\Tc$ is also proper.

The first part of $(ii)$ follows from Lemma 4.1.$(ii)$.
The projection $\Zc_y\to\Oc_y$ is a vector bundle, 
and the fiber at $(y^{-1}u_2\cdot 1_\Xc,1_\Xc)$ is 
$\uen\cap y^{-1}(\ad_{u_2}\uen)$.
Since $\ell\gg 0$ we can assume that 
$\uen_\ell\subset\uen\cap y^{-1}(\ad_{u_2}\uen)$.
Set $\Zc^\ell_y=(p_\ell\times p_\ell)(\Zc_y).$ 
Then the projection $\Zc^\ell_y\to\Oc_y$ is a vector bundle again.
It is $U^-$-equivariant and $\dot W$-equivariant, 
and the fiber over $(y^{-1}u_2\cdot 1_\Xc,1_\Xc)$ is 
$(\uen\cap y^{-1}(\ad_{u_2}\uen))/(\uen_\ell\cap y^{-1}(\ad_{u_2}\uen_\ell))$.
We have $\uen_\ell\cap y^{-1}(\ad_{u_2}\uen_\ell)=
\ad_{\dot y^{-1}(u_2)}(\uen_\ell\cap y^{-1}(\uen_\ell))$
because $[\uen,\uen_\ell]\subset\uen_\ell$,
and $\uen\cap y^{-1}(\ad_{u_2}\uen)=
\ad_{\dot y^{-1}(u_2)}(\uen\cap y^{-1}(\uen)).$
Taking the quotient by the left $(U^-_k)^2$-action we get $(ii)$.

As $q_{23}^{-1}({}^w\Zc_{\leq y})$
is contained in $\Tc\times{}^{w'}\Tc\times{}^w\Tc$ by $(ii)$,
Lemma 4.1.$(iii)$ implies that 
$q_{13}\bigl(q_{12}^{-1}(\Zc_{\leq y'})\cap
q_{23}^{-1}({}^w\Zc_{\leq y})\bigr)$
is in ${}^w\Zc_{\le y''}$.
Hence $(iii)$ is a consequence of $(i)$.
\qed

\vskip3mm

\noindent{\bf Remarks.}
$(i)$
$\Zc_{\leq y}$, $\Zc^{k\ell}_{\leq y}$, etc,
are endowed with the reduced scheme structure induced by $\Tc^2$
and $({}^{w}\Tc^{k\ell})^2$ respectively. 

$(ii)$
We do not know if 
${}^w\Zc^{k\ell}_{\leq y}$
is closed in ${}^{w'}\Tc^{k\ell}\times{}^w\Tc^{k\ell}$ 
because ${}^w\Xc^k$ may be not separated.

\subhead 4.3\endsubhead
Set
$\Zc^a_{ij}=\Zc^a\cap(\dot\Nc^a_i\times\dot\Nc^a_j)$
for all $i,j\in\Xi_s$, see 2.13.
In the rest of Section 4 we set also
$\uen_\ell=\eps^\ell\cdot\uen_J$ with $J\subsetneq I$.
Moreover we assume assume that $w\in W^Jw_J$.

\proclaim{Proposition}
Assume that $k,k_2,k_1\gg 0$, 
$k_2\geq k_1$, $\ell_2\geq \ell_1$, and $y\in W_J$.

(i)
${}^w\Zc_{\leq y}$, ${}^w\Zc^\ell_{\leq y}$, and ${}^w\Zc^{k\ell}_{\leq y}$ 
are $G_J$-subschemes of
$({}^w\Tc)^2$, $({}^w\Tc^\ell)^2$,
and $({}^w\Tc^{k\ell})^2$.

(ii)
$q_{13}$ restricts to a proper map
$q_{12}^{-1}({}^w\Zc_{\leq w_J})\cap 
q_{23}^{-1}({}^w\Zc_{\leq w_J})\to{}^w\Zc_{\leq w_J}$.

(iii) 
The projections
${}^w\Zc^{\ell_2}_{\leq y}\to{}^w\Zc^{\ell_1}_{\leq y}$,
${}^w\Zc^{k\ell_2}_{\leq y}\to{}^w\Zc^{k\ell_1}_{\leq y}$ 
are vector bundles with fiber $\uen_{\ell_2}/\uen_{\ell_1}$.

(iv)
$\Zc^a=\{(x,\pen;x',\pen')\in\dot\Nc^a\times\dot\Nc^a\,;\,x=x'\}.$

(v)
The map $\rho\circ q_2$
yields a locally trivial fibration 
$\Zc^a_{ij}\cap\Zc_{\le y}\to\Xc^s_j$.
\endproclaim

\noindent{\sl Proof:}
Immediate from Lemma 2.13.$(iv)$ and Lemma 4.2.
Note that our assumptions imply that $w\succeq(w,y^{-1})$ and
$w_J\succeq(w_J,w_J)$.
\qed

\vskip3mm

\noindent{\bf Remark.}
Note also the following, to be compared with (A.6), (A.7).

\vskip2mm

$(i)$
$\Id\times p_{\ell_2\ell_1}$ restricts to an isomorphism 
${}^{w}\Zc_{y}^{k\ell_2}\to
(p_{\ell_2\ell_1}\times\Id)^{-1}({}^{w}\Zc_{y}^{k\ell_1})$.

\vskip2mm

$(ii)$
$p^{k_2k_1}\times\Id$ restricts to an isomorphism 
${}^{w}\Zc_{y}^{k_2\ell}\to(\Id\times 
p^{k_2k_1})^{-1}({}^{w}\Zc_{y}^{k_1\ell})$.

\vskip2mm

\noindent
We will not use it.
Proof is left to the reader. 

\subhead 4.4\endsubhead
Fix $y,y',y'',w,w'$ as in Lemma 4.1.
Write ${}^w\Zc_{\leq y}^a$ for the fixed points
set $({}^w\Zc_{\leq y})^a$, and similarly for the other sets we consider. 
Let $D\subset A$ be a closed subgroup.
Then ${}^w\Zc_{\leq y}^a$ is closed in ${}^{w'}\Tc^a\times{}^w\Tc^a$,
and $q_{13}$ restricts to a proper map
$q_{12}^{-1}({}^{w'}\Zc_{\leq y'}^a)\cap
q_{23}^{-1}({}^{w}\Zc_{\leq y}^a)\to {}^{w}\Zc_{\leq y''}^a$,
yielding a $\star$-product 
$$\Kb^{D}({}^{w'}\Zc^a_{\leq y'})
\times\Kb^{D}({}^{w}\Zc^a_{\leq y})\to
\Kb^{D}({}^{w}\Zc^a_{\leq y''})\leqno(4.4.1)$$
relative to $(\Tc^a)^3$, see A.2.
If $w\leq w'$ the open immersion 
${}^{w}\Tc^a\subset{}^{w'}\Tc^a$ 
gives a map 
$$h^{w'w}\,:\,\Kb^{D}({}^{w'}\Zc^a_{\leq y})\to\Kb^{D}({}^{w}\Zc^a_{\leq y}).$$
If $y\leq y'$ the closed immersion ${}^w\Zc_{\leq y}^a\to{}^w\Zc^a_{\leq y'}$
gives a map
$$h_{yy'}\,:\,\Kb^D({}^w\Zc_{\leq y}^a)\to\Kb^D({}^{w}\Zc^a_{\leq y'}).$$
The maps $h^{w'w}$ and $h_{yy'}$ commute.
Set $$\Kb^{D}(\Zc_{\leq y}^a)=\pro_w(\Kb^{D}({}^w\Zc^a_{\leq y}),h^{w'w}).$$
Then we have
$h_{yy'}\,:\,\Kb^D(\Zc_{\leq y}^a)\to\Kb^D(\Zc^a_{\leq y'}),$
and we set $$\Kb^{D}(\Zc^a)=\ind_y(\Kb^{D}(\Zc^a_{\leq y}),h_{yy'}).$$
To simplify we write $\Kb_a=\Kb^A(\Zc^a)_a$
and $\Kb_{a,\leq y}=\Kb^A(\Zc^a_{\leq y})_a$.
Recall that the subscript $a$ means specialization at the maximal ideal $\Jb_a$,
see 1.4.

\proclaim{Proposition}
(i)
(4.4.1) induces a $\star$-product
$\Kb_{a,\leq y'}\times\Kb_{a,\leq y}\to\Kb_{a,\leq y''}$,
yielding a ring $(\Kb_a,\star)$.
The natural map
$\Kb_{a,\leq w_J}\to\Kb_a$ is a ring homomorphism.

(ii)
If $s_1\in\ad_G(s_2)$ and $a_i=(s_i,\zeta)$
the rings $\Kb_{a_1}$, $\Kb_{a_2}$ are isomorphic.
\endproclaim

\noindent{\sl Proof:}
The maps $h^{w'w}$ and $h_{yy'}$ are $\star$-homomorphisms, 
see Example (A.3.1). Thus the bilinear map
$\star\,:\,\Kb_{a,\leq y'}\times\Kb_{a,\leq y}\to\Kb_{a,\leq y''}$
such that $(a'_w)\star(a_w)=(a'_w\star a_w)$ is well-defined
and is compatible with the maps $h_{yy'}$.
This yields Claim $(i)$.
It is known that $s_1\in\ad_G(s_2)$ if and only if
$s_1$, $s_2$ are $W$-conjugated.
Hence $\Kb_{a_1}\simeq\Kb_{a_2}$,
because $\dot W$ acts on $\Zc$ and normalizes $A$.
\qed

\vskip3mm

\noindent{\bf Remark.}
Observe that, for each $w$ there are finite subsets $S,T\subset\Xi_s$
such that 
$${}^w\Zc^a_{\leq y}\subseteq\bigcup_{i\in S}\bigcup_{j\in T}
\Zc_{ij}^a\cap{}^w\Zc_{\leq y}.$$
Inversely, given $S, T$ there is a $w$ such that the reverse inclusion holds.
Thus
$$\Kb^D(\Zc^a_{\leq y})=\prod_{i,j}\Kb^D(\Zc^a_{ij}\cap\Zc_{\leq y}).$$
Given $j$, the group
$\Kb^D(\Zc^a_{ij}\cap\Zc_{\leq y})$
vanishes except for a finite number of $i$.
Thus
$$\prod_{i,j}\Kb^D(\Zc^a_{ij}\cap\Zc_{\leq y})=
\prod_j\bigoplus_i\Kb^D(\Zc^a_{ij}\cap\Zc_{\leq y}).$$
Moreover, for each $i,j\in\Xi_s$ there is an element $y$ such that 
$\Zc^a_{ij}\subseteq\Zc_{\le y}^a.$
Therefore we have a natural injective map
$$\Kb_a\to
\prod_j\ind_y\bigoplus_i\Kb^A(\Zc^a_{ij}\cap\Zc_{\le y})_a=
\prod_j\bigoplus_i\Kb^A(\Zc^a_{ij})_a.$$

\subhead 4.5\endsubhead
Let ${}^w\Tc$, $\Tc^\ell$, ${}^w\Pc^{J,k}$ be as in 2.5, 2.6, and 2.10.
Recall that $\uen_0=\uen_J$, see 2.9, 4.3.
There is an isomorphism
$${}^w\Pc^{J,k}\times_{G_J}\dot\Nc_J\to{}^w\Tc^{k0},\quad
[g\cdot U_J:g_J\cdot(x+\uen_J,\ben)]\mapsto[g\cdot g_J:x+\uen_J]$$
where $g\in U_k^-\setminus{}^w\Gc$, $g_J\in G_J$, and $x\in\uen$.
The composed map
$${}^w\Pc^{J,k}\times_{G_J}\Zc_J\to
{}^w\Pc^{J,k}\times_{G_J}(\dot\Nc_J)^2\to
({}^w\Tc^{k0})^2\leqno(4.5.1)$$
is an immersion.

\proclaim{Lemma}
Let $y\in W_J$.

(i)
(4.5.1) restricts to an isomorphism
${}^w\Pc^{J,k}\times_{G_J}\Zc_{J,\leq y}\to{}^w\Zc^{k0}_{\leq y}$.

(ii)
${}^w\Zc_{\leq y}^{k\ell}$ is closed in ${}^w\Zc^{k\ell}_{\leq w_J}$.
\endproclaim

\noindent{\sl Proof:}
Observe that ${}^w\Zc^{k0}_{\le y}$ is a $G_J$-subscheme of
$({}^w\Tc^{k0})^2$ by Proposition 4.3.$(i)$.
The immersion (4.5.1)
takes $[g\cdot U_J:g_J\cdot(x+\uen_J,\ben_{y^{-1}},\ben)]$
to $([g\cdot g_J\dot y^{-1}:\dot y(x)+\uen_J], [g\cdot g_J:x+\uen_J])$
for any $x\in\ben_{y^{-1}}\cap\ben$.
Hence ${}^w\Pc^{J,k}\times_{G_J}\Zc_{J,y}$ surjects to ${}^w\Zc^{k0}_y$,
yielding $(i)$.
Claim $(ii)$ is immediate because the natural projection
${}^w\Zc_{\leq y}^{k\ell}\to{}^w\Zc_{\leq y}^{k0}$
is a vector bundle, see 4.3, 
and ${}^w\Zc_{\leq y}^{k0}$ is closed in ${}^w\Zc^{k0}_{\leq w_J}$ by $(i)$. 
\qed

\vskip3mm

The $\star$-product on 
$\Kb^{A}({}^w\Zc_{\leq w_J}^{k0})$ 
relative to 
${}^w\Pc^{J,k}\times_{G_J}(\dot\Nc_J)^3$
makes sense, see (A.4.1), because ${}^w\Xc^{J,k}$ is smooth.
Let $p^{k_2k_1}$ denote also the obvious projection 
${}^w\Pc^{J,k_2}\to{}^w\Pc^{J,k_1}.$
It is a smooth map.
Set
$$\Kb^A({}^w\Zc^0_{\leq w_J})
=\ind_k\bigl(\Kb^A({}^w\Zc^{k0}_{\leq w_J}),
(p^{k_2k_1}_{\Zc_J})^*\bigr),$$ 
where $p^{k_2k_1}_{\Zc_J}=(p^{k_2k_1})_{\Zc_J}$, see (A.4.4).

On the other hand, we endow $\Kb^A(({}^w\Zc^0_{\leq w_J})^a)$ with the
$\star$-product relative to $(({}^w\Tc^0)^a)^3$.
This is possible, since if $k$ is large then 
$({}^w\Tc^{k0})^a\simeq ({}^w\Tc^0)^a$
by Lemma 2.13.$(vi)$.
Thus $({}^w\Zc^0_{\leq w_J})^a$ 
is a scheme of finite type
and $\Kb^A(({}^w\Zc^0_{\leq w_J})^a)$ 
is well-defined, without any use of direct or inverse limit, see 1.4.
Moreover $({}^w\Zc^0_{\leq w_J})^a$ 
is closed in $(({}^w\Tc^0)^a)^2$ 
and $({}^w\Tc^0)^a$ is smooth, see Proposition 4.3,
hence the $\star$-product is well-defined.
The concentration map
$$\Kb^A({}^w\Zc^{k0}_{\leq w_J})_a\to\Kb^A(({}^w\Zc^0_{\leq w_J})^a)_a
\leqno(4.5.2)$$
relative to ${}^w\Pc^{J,k}\times_{G_J}(\dot\Nc_J)^2$ 
commutes with $\star$, see (A.4.3).

\proclaim{Lemma}
(iii)
$\star$ factorizes through a product on $\Kb^A({}^w\Zc^0_{\leq w_J}).$

%(iv)
%${}^\sharp p_{\ell 0}$
%yields a ring isomorphism
%${}^\sharp p_0\,:\,\Kb^A({}^w\Zc^0_{\leq w_J})
%\to\Kb^A({}^w\Zc_{\leq w_J}).$

(iv)
(4.5.2) factorizes through a ring isomorphism
$$\cb_a\,:\,\Kb^A({}^w\Zc^0_{\leq w_J})_a\to
\Kb^A(({}^w\Zc^0_{\leq w_J})^a)_a.$$
\endproclaim

\noindent{\sl Proof:}
Immediate by (A.4.4). 
Note that $\Kb^A(({}^w\Zc^0_{\leq w_J})^a)$ 
is also the direct limit
$\ind_k\Kb^A(({}^w\Zc^{k0}_{\leq w_J})^a),$ 
because the limit stabilizes.
\qed

\vskip3mm

\noindent{\bf Remark.}
We do not use a $\star$-product on
$\Kb^A({}^{w}\Zc^{k\ell}_{\leq w_J})$ 
relative to $({}^w\Tc^{k\ell})^3$,
because the scheme ${}^w\Xc^k$ may be not separated, 
see Remark 4.2.$(ii)$.

\subhead 4.6\endsubhead
Following \cite{KT2, 2.2}, for any $\l\in X$ 
we define $H$-equivariant invertible sheaves
${}^wL_\l^k$ on ${}^w\Xc^k$ whose local sections are
$$\Gamma(U;{}^wL^k_\l)=\{f\in\Gamma((p^k\pi)^{-1}(U);\Oc_\Gc)\,;\,
f(ugb)=f(g)b^{-\l},\,\forall u,g,b\},$$
where $(u,g,b)\in U_k^-\times(p^k\pi)^{-1}(U)\times B.$
Let ${}^wL_\l^{k\ell}$ denote also
the pull-back of ${}^wL_\l^k$ to ${}^w\Tc^{k\ell}$.
The sheaves ${}^wL_\l^{k\ell}$ 
are made $A$-equivariant with $\CC_q^\times$ acting trivially.
Set 
$$\Kb_H(\Xc)=\pro_w\Kb_H({}^w\Xc),\quad
\Kb_A(\Tc)=\pro_w\Kb_A({}^w\Tc),$$
where
$$\Kb_H({}^w\Xc)=\ind_k(\Kb_H({}^w\Xc^k),p^{k_2k_1*}),\quad
\Kb_A({}^w\Tc)=\ind_{k,\ell}
(\Kb_A({}^w\Tc^{k\ell}),p^{k_2k_1*},p^*_{\ell_2\ell_1}),$$
see A.5.
Let ${}^wL_\l\in\Kb_H({}^w\Xc),\Kb_A({}^w\Tc)$
be the element representing the systems of sheaves 
$({}^wL_\l^k)$ and $({}^wL_\l^{k\ell})$ respectively.
Let $L_\l\in\Kb_H(\Xc),\Kb_A(\Tc)$
be the element representing the system of sheaves $({}^wL_\l)$.

For a future use, we set also
$$\Kb_H(\Xc^s)=\pro_w\Kb_H({}^w\Xc^s),\quad
\Kb_A(\Tc^a)=\pro_w\Kb_A({}^w\Tc^a).$$

\vskip3mm

\noindent{\bf Remark.}
Let $\Lc_\l$ be the pull-back to $\Tc$ of the
$H$-equivariant invertible sheaf on $\Xc$ 
whose local sections on $U$ are
$\{f\in\Gamma(\pi^{-1}(U);\Oc_\Gc)\,;\,
f(gb)=f(g)b^{-\l},\ \forall(g,b)\in\pi^{-1}(U)\times B\}.$
Then the inverse image of ${}^wL^{k\ell}_\l$ by the projection
${}^w\Tc\to{}^w\Tc^{k\ell}$ equals $\Lc_\l|_{{}^w\Tc}$.
Thus ${}^wL_\l$ and $L_\l$ are well-defined.

\subhead 4.7\endsubhead
Fix $\l,\l',\l''\in X$ such that $\l=\l'+\l''$.
Let $D_\l$ be the restriction of
$L_{\l'}\boxtimes L_{\l''}$ to $\Zc_1$,
and ${}^wD^{k\ell}_\l$ be the restriction of
${}^wL^{k\ell}_{\l'}\boxtimes {}^wL^{k\ell}_{\l''}$
to ${}^w\Zc_1^{k\ell}$.
Since ${}^w\Zc_1^{k\ell}$ is closed in ${}^w\Zc^{k\ell}_{\leq w_J}$
by Lemma 4.5.$(ii)$, we may view ${}^wD^{k0}_\l$ 
as an element in $\Kb^A({}^w\Zc^{k0}_{\leq w_J})$.
We have
$$(p^{k_2k_1}_{\Zc_J})^*({}^wD_\l^{k_10})={}^wD_\l^{k_20}.$$
%$$\matrix
%(\Id\times p^{k_2k_1})^*({}^wD_\l^{k_1\ell})
%&=({}^wL_{\l'}^{k_1\ell}\boxtimes{}^wL_{\l''}^{k_2\ell})
%|_{(\Id\times p^{k_2k_1})^{-1}({}^w\Zc_1^{k_1\ell})}\hfill\cr
%&=({}^wL_{\l'}^{k_1\ell}\boxtimes{}^wL_{\l''}^{k_2\ell})
%|_{(p^{k_2k_1}\times\Id)({}^w\Zc_1^{k_2\ell})}\hfill\cr
%&=(p^{k_2k_1}\times\Id)_*((p^{k_2k_1}\times\Id)^*
%({}^wL_{\l'}^{k_1\ell}\boxtimes{}^wL_{\l''}^{k_2\ell})
%|_{{}^w\Zc_1^{k_2\ell}})
%\hfill\cr
%&=(p^{k_2k_1}\times\Id)_*({}^wD_\l^{k_2\ell}).
%\hfill\cr
%\endmatrix$$
Let ${}^wD^0_\l$ 
represent the system $({}^wD^{k0}_\l)$ 
in $\Kb^A({}^w\Zc^0_{\leq w_J})$.

Assume that $\l=-\a_i$ with $i\in J$,
and that $(\l':\a_i^\vee)=(\l'':\a_i^\vee)=-1$.
Let $D_i$ be the restriction of $q^{-1}L_{\l'}\boxtimes L_{\l''}$ 
to $\bar\Zc_{s_i}=\Zc_{s_i}\cup\Zc_1$.
The set $({}^w\bar\Zc_{s_i})^{k\ell}=(p^k_\ell\times p^k_\ell)
(\bar\Zc_{s_i}\cap(\Tc\times{}^w\Tc))$
is closed in ${}^w\Zc^{k\ell}_{\leq w_J}$ 
(as for Lemma 4.5.$(ii)$).
Let ${}^wD_i^{k\ell}$ be the restriction of 
$q^{-1}({}^wL^{k\ell}_{\l'}\boxtimes {}^wL^{k\ell}_{\l''})$ 
to $({}^w\bar\Zc_{s_i})^{k\ell}$.
We have
$$(p^{k_2k_1}_{\Zc_J})^*({}^wD_i^{k_10})={}^wD_i^{k_20}.$$
Let ${}^wD^0_i\in\Kb^A({}^w\Zc^0_{\leq w_J})$ represent the system
$({}^wD^{k0}_i)$ in $\Kb^A({}^w\Zc^0_{\leq w_J})$.

Let $p_a\,:\,\Tc^a\to(\Tc^0)^a$ be the restriction of the vector bundle 
$p_0\,:\,\Tc\to\Tc^0$ in 2.6 to the fixed point subsets.
The map
${}^\sharp p_a\,:\,\Kb^A(({}^w\Zc^0_{\leq w_J})^a)_a
\to\Kb^A({}^w\Zc^a_{\leq w_J})_a$
is invertible.
%Set 
%$\cb_a\,:\,\Kb^A({}^w\Zc_{\leq w_J})_a\to\Kb^A({}^w\Zc^a_{\leq w_J})_a$
%equal to
%${}^\sharp p_a\circ\cb_a\circ({}^\sharp p_0)^{-1}$.
Put $D_{\l,a}=({}^\sharp p_a\cb_a({}^wD^0_\l))$ and
$D_{i,a}=({}^\sharp p_a\cb_a({}^wD^0_i))$ in $\Kb_{a,\leq w_J}$,
where $\cb_a$ is as in Lemma 4.5.$(iv)$.
Let $D_{\l,a}$, $D_{i,a}$ denote also their images  
in $\Kb_a$ by the obvious map.

\proclaim{Lemma}
$D_{\l,a}$, $D_{i,a}$
do not depend on the choice of $\l'$, $\l''$, or $J$.
\endproclaim

\noindent{\sl Proof:}
The independence on $\l',\l''$ is proved as in \cite{L2, 7.19}.
We prove the independence on $J$ for $D_{\l,a}$ only.
The case of $D_{i,a}$ is similar.
Assume that $J_1\subset J_2\subsetneq I$,
and that $w\in W^{J_2}w_{J_2}$.
Fix $i=1$ or 2.
Given $\uen_\ell=\eps^\ell\cdot\uen_{J_i}$,
we write ${}^w\Tc^{k\ell}_i$ for ${}^w\Tc^{k\ell}$,
and ${}^w\Zc^{k\ell}_{i,\leq y}$ for ${}^w\Zc^{k\ell}_{\leq y}$.
Both schemes depend on $J_i$ because the ideal $\uen_\ell$ depends on it.
Write ${}^wL_{\l i}^{k\ell}$ for the element ${}^wL_\l^{k\ell}$
in $\Kb_A({}^w\Tc^{k\ell}_i)$,
and ${}^wD_{\l i}^{k\ell}$ for the element ${}^wD_\l^{k\ell}$
in $\Kb^A({}^w\Zc^{k\ell}_{i,\leq w_{J_1}})$.
The pro-objects 
$$({}^w\Tc^{k\ell}_1,p^{k_2k_1},p_{\ell_2\ell_1}),\quad
({}^w\Tc^{k\ell}_2,p^{k_2k_1},p_{\ell_2\ell_1})$$
are isomorphic, see 2.6 and A.5, and the systems 
$({}^wL_{\l 1}^{k\ell})$ and $({}^wL_{\l 2}^{k\ell})$ 
yield the same element in $\Kb_A({}^w\Tc)$.
Hence the ind-objects
$$\bigl(\Kb^A({}^{w}\Zc^{k0}_{1,\leq w_{J_1}}),
(p^{k_2k_1}_{\Zc_{J_1}})^*\bigr),\quad
\bigl(\Kb^A({}^{w}\Zc^{k0}_{2,\leq w_{J_1}}),
(p^{k_2k_1}_{\Zc_{J_2}})^*\bigr),$$
are isomorphic, and the systems
$({}^wD_{\l 1}^{k0})$ and $({}^wD_{\l 2}^{k0})$
yield the same element ${}^wD_\l^0$ in $\Kb^A({}^w\Zc^0_{\leq w_{J_1}})$.
See A.5 for generalities on pro-objects and ind-objects.
Remark that
${}^\sharp p_a\circ\cb_a\,:\,\Kb^A({}^w\Zc^0_{\leq w_{J_1}})_a\to
\Kb^A({}^w\Zc_{\leq w_{J_1}}^a)_a$ depends a priori on the choice of $J$
(=$J_1$ or $J_2$).
A direct computation yields
$$\matrix
{}^\sharp p_a\cb_a({}^wD^0_\l)
&={}^\sharp p_a\cb_a
((\CC\boxtimes {}^wL_\l^{k0})|_{{}^w\Zc_1^{k0}}),\hfill\cr
&={}^\sharp p_a
((\CC\boxtimes\Lc_\l)|_{({}^w\Zc^0_1)^a}),\hfill\cr
&=(\CC\boxtimes\Lc_\l)|_{{}^w\Zc_1^a}.\hfill
\endmatrix\leqno(4.7.1)$$
Hence the element $({}^\sharp p_a\cb_a({}^wD^0_\l))$ in $\Kb_a$
does not depend on this choice indeed.
\qed

\subhead 4.8\endsubhead
The $\CC$-linear map $\Rb_H\to\Kb(\Xc)$ such that $x_\l\mapsto L_\l$
is not surjective.
Let ${}^w\b_{s}\,:\,\Rb_H\to\Kb_H({}^w\Xc^{s})_{s}$
be the unique $\CC$-linear map taking $x_\l$ to $\Lc_\l|_{{}^w\Xc^{s}}$.
Let $I_{w,s}\subset\Rb_A$ be the ideal generated by $\Ker({}^w\b_{s})$
and $q-\zeta$, and let $\Ib_a\subset\Rb_A$ be the ideal
generated by $q-\zeta$ and $\{f-f(s)\,;\,f\in(\Rb_H)^{W_s}\},$
where $W_s$ is as in 2.2.
If $w'\leq w$ then $I_{w,s}\subseteq I_{w',s}$
because ${}^{w'}\Xc^s\subset{}^{w}\Xc^s$.
Set $\hat\Rb_a=\pro_w\Rb_A/I_{w,s}$.
Let $W(a)$ be the $W$-orbit of $a$.

\proclaim{Lemma}
(i)
${}^w\b_{s}$ is a surjective ring homomorphism.

(ii)
The $\Rb_A$-algebras
$\hat\Rb_a$, $\prod_{b\in W(a)}\Rb_A/\Ib_b$, and $\Kb_A(\Tc^a)_a$ 
are isomorphic. 

(iii)
The assignement $\Lc_\l|_{\Tc^a}\mapsto D_{\l,a}$, with $\l\in X$,
extends uniquely to a $\CC$-algebra isomorphism
$(\Kb_A(\Tc^a)_a,\otimes)\to(\Kb_{a,\leq 1},\star)$.
\endproclaim

\noindent{\sl Proof:}
The map ${}^w\b_{s}$ is a ring homomorphism. 
We must prove that it is surjective.
The $G$-action on $\Bc$ restricts to a $G^{s}$-action on
$\Xc_i^{s}$, and the $G^{s}$-schemes $\Xc_i^{s}$ are isomorphic to
the flag variety of the connected reductive group $G^{s}$, see Lemma 2.13.
Hence the ring homomorphism
$\Rb_{H}\to\Kb_H(\Xc_i^{s})_{s}$
taking $x_\l$ to $\Lc_\l|_{\Xc^{s}_i}$ is surjective.
It takes the element $t=x_\delta$ to $\tau$.
The RHS is a finite dimensional algebra. 
Hence it is the ring of functions of a finite length subscheme 
$\xi_i\subset H'\times\{\tau\}$.
It is sufficient to prove that the supports
of the schemes $\xi_i$ are disjoint.
Since $\Kb_H(\Xc_i^{s})_{s}$
is isomorphic to the De Rham cohomology ring of $\Xc_i^{s}$, 
it is $\ZZ_{\geq 0}$-graded with a one-dimensional zero component.
Thus $\xi_i$ has a punctual support.
Let $|\xi_i|\in H'\times\{\tau\}$ denote the support of $\xi_i$.
Each connected component
$\Xc_i^{s}$ contains an element $\ben_w$ for some $w\in W$.
We have $|\xi_i|=w^{-1}(s)$ whenever $\ben_w\in\Xc_i^{s}$.
The group $W$ is the quotient of $N_{G_{KM}}(H')$ by $H'$,
and the $N_{G_{KM}}(H')$-action on $H$ by conjugation
coincides with the $W$-action.
Hence if $w_1^{-1}(s)=w_2^{-1}(s)$,
then $\ben_{w_1}=\ad_g(\ben_{w_2})$ 
for some element $g\in G^s$,
and, thus, $\ben_{w_1}$, $\ben_{w_2}$ belong to the same connected
component of $\Xc^{s}$ by Lemma 2.13.$(i)$.
This proves $(i)$.

Using $(i)$ we get 
$\Kb_H(\Xc^s)_s=\pro_w\Rb_H/\Ker({}^w\b_s).$
Hence $\Kb_A(\Tc^a)_a=\hat\Rb_a$ because $\rho$ restricts
to a vector bundle ${}^w\Tc^a\to{}^w\Xc^s$, see Lemma 2.13.
By definition $\Kb_H(\Xc^s)_s=\prod_{i\in\Xi_s}\Kb_H(\Xc_i^s)_s$. 
Since $\Xc_i^s$ is isomorphic to the flag variety of $G^s$,
the $\Rb_H$-algebra $\Kb_H(\Xc_i^s)_s$ is isomorphic to 
$\Rb_H/(f-f(w(s))\,;\,f\in(\Rb_H)^{W_{w(s)}})$, see Lemma 2.13 again,
where $w\in W$ is such that
$\ben\in\Xc_i^{w(s)}$.
This proves $(ii)$.

Let us prove $(iii)$.
The variety ${}^w\Tc^a$ is smooth, 
hence $\Kb_A({}^w\Tc^a)\simeq\Kb^A({}^w\Tc^a)$.
The second projection ${}^w\Zc_1^a\to{}^w\Tc^a$ is an isomorphism,
hence $\Kb^A({}^w\Tc^a)\simeq\Kb^A({}^w\Zc_1^a)$.
The composed isomorphism
$\Kb_A({}^w\Tc^a)\to\Kb^A({}^w\Zc^a_1)$
takes the tensor product to $\star$.
It takes also $\Lc_\l|_{{}^w\Tc^a}$ to $(\CC\boxtimes\Lc_\l)|_{{}^w\Zc^a_1},$
i.e. to ${}^\sharp p_a\cb_a({}^wD^0_\l)$ by (4.7.1).
The unicity in $(iii)$ follows from $(ii)$.
\qed

\vskip3mm

\noindent{\bf 4.9.}
For any $y$ the space $\Kb_{a,\leq y}$ 
is endowed with the topology
induced by the product $\prod_w\Kb^A({}^w\Zc^a_{\leq y})_a$, 
where each term is given the discrete topology.
The ring $\Kb_a$ is endowed with the smallest topology such that
the natural maps $\Kb_{a,\leq y}\to\Kb_a$ are continuous.
Recall that a representation $\Mb$ of a topological ring $\Ab$
is {\sl smooth} if for any $m\in\Mb$ the
subset $\{x\in\Ab\,;\,x(m)=0\}$ is open.

\proclaim{Theorem}
(i)
There is a unique ring homomorphism
$\Psi_a\,:\,\Hb\to\Kb_a$ such that $t_i\mapsto -1-qD_{i,a}$,
and $x_\l\mapsto D_{\l,a}$.

(ii)
$\Kb_{a,\leq y}$ is a free right
$\Kb_{a,\leq 1}$-module of rank $\ell(y)$ with basis 
$(\Psi_a(t_{y'})\,; \,y'\leq y)$.

(iii)
For any closed subgroup $D\subset A$
the forgetful map $\Kb^A(\Zc^a_{\leq y})\to\Kb^{D}(\Zc^a_{\leq y})$
is an isomorphism.

(iv)
A smooth and simple $\Kb_a$-module pulls-back to a simple $\Hb$-module.
Non-isomorphic smooth $\Kb_a$-modules pull-back to non-isomorphic $\Hb$-modules.
\endproclaim

\noindent{\sl Proof of $(i)$:}
The group $G_J$ is connected 
with a simply connected derived group.
Thus there is a $\star$-homomorphism 
$\Hb_J\to\Kb^{G_J\times\CC_q^\times}(\Zc_J)$ 
such that $t_i\mapsto -1-qD_i$, $x_\l\mapsto D_\l$, and
$D_\l,D_i\in\Kb^{G_J\times\CC^\times_q}(\Zc_J)$ are as in 4.7,
see \cite{L2, Theorem 7.25}.
Therefore, due to Lemma 3.2 and Lemma 4.7, 
in order to prove $(i)$ it is sufficient to construct
a $\star$-homomorphism 
$\Theta_J\,:\,\Kb^{G_J\times\CC_q^\times}(\Zc_J)\to\Kb_a$
such that $\Theta_J(D_\l)=D_{\l,a}$
and $\Theta_J(D_i)=D_{i,a}$.
To do so we first construct, for each $w$, a map 
$${}^w\Theta_J\,:\,\Kb^{G_J\times\CC_q^\times}(\Zc_J)
\to\Kb^A({}^w\Zc^a_{\leq w_J})_a$$
such that ${}^w\Theta_J$ commutes with $\star$ and
$h^{w'w}\circ{}^{w'}\Theta_J={}^w\Theta_J$
for each $w$, $w'$.
Then we compose the resulting $\star$-homomorphism 
$\Kb^{G_J\times\CC_q^\times}(\Zc_J)\to\Kb_{a,\leq w_J}$
with the natural map $\Kb_{a,\leq w_J}\to\Kb_a$.

\vskip2mm

(a)
Let 
$\varphi_k\,:\,{}^w\Pc^{J,k}\times\Zc_J\to\Zc_J,$
$\psi_k\,:\,{}^w\Pc^{J,k}\times\Zc_J\to{}^w\Pc^{J,k}\times_{G_J}\Zc_J$
be the obvious projections. 
Set
$\Gamma_k=(\psi_k^*)^{-1}\circ(\varphi_k)^*\,:\,
\Kb^{G_J\times\CC_q^\times}(\Zc_J)\to
\Kb^A({}^w\Zc^{k0}_{\leq w_J}),$
see 4.5.
We have
$(p^{k_2k_1}_{\Zc_J})^*\circ\Gamma_{k_1}=\Gamma_{k_2}$.
Hence we get a map
$\Gamma\,:\,\Kb^{G_J\times\CC_q^\times}(\Zc_J)\to
\Kb^A({}^w\Zc^0_{\leq w_J}).$

\vskip2mm

(b)
Let
$\cb_a\,:\,\Kb^A({}^w\Zc^0_{\leq w_J})_a\to
\Kb^A(({}^w\Zc^0_{\leq w_J})^a)_a$ 
and
${}^\sharp p_a\,:\,\Kb^A(({}^w\Zc^0_{\leq w_J})^a)_a
\to\Kb^A({}^w\Zc^a_{\leq w_J})_a$
be as in 4.5, 4.7.

\vskip2mm

We put ${}^w\Theta_J={}^\sharp p_a\circ\cb_a\circ\Gamma$.
The compatibility with $h^{w'w}$ is immediate. 
The maps $\Gamma$, $\cb_a$, and ${}^\sharp p_a$
commute with $\star$ by (A.4.2), (A.4.3),
and (A.3.2) respectively.
We have $\Theta_J(D_\l)=D_{\l,a}$,
because $\Gamma_k(D_\l)={}^wD^{k0}_\l.$ 
Similarly $\Theta_J(D_i)=D_{i,a}$.
\qed

\vskip3mm

\noindent{\sl Proof of $(iii)$ :}
Fix $i,j\in\Xi_s$
such that $\Zc^a_{ij}\cap\Zc_y\neq\emptyset$.
%Set $\Kb^D(\Zc^a_y)=\pro_w\Kb_D({}^w\Zc^a_y)$.
Composing the projections
$\Zc_y\to\Oc_y$ and $\Oc_y\to\Xc,$ 
we get a vector bundle 
$\Zc^a_{ij}\cap\Zc_y\to\Xc^s_j$,
see Lemma 2.14 and Lemma 4.2.
Similarly, composing the projections
$\Zc_{\le y}\to\Oc_{\le y}$ and $\Oc_{\le y}\to\Xc,$ 
we get a locally trivial fibration 
$\Zc^a_{ij}\cap\Zc_{\le y}\to\Xc^s_j$, see Proposition 4.3.
Moreover the $\Rb_D$-module $\Kb^D(\Xc^s_j)$ is free
and $\Kb_1^D(\Xc_j^s)=0$, see Lemma 2.13.
Hence 
$$0\to\Kb^D(\Zc^a_{ij}\cap\Zc_{<y})
\to\Kb^D(\Zc^a_{ij}\cap\Zc_{\leq y})
\to\Kb^D(\Zc^a_{ij}\cap\Zc_{y})\to 0$$
is an exact sequence of $\Rb_D$-modules by the cellular fibration lemma,
and the forgetful map 
$\Rb_D\otimes_{\Rb_A}\Kb^A(\Zc^a_{ij}\cap\Zc_{\leq y})
\to\Kb^D(\Zc^a_{ij}\cap\Zc_{\leq y})$
is an isomorphism, see \cite{CG, 5.5}.
Therefore
$$0\to\Kb^D(\Zc^a_{<y})
\to\Kb^D(\Zc^a_{\leq y})
\to\Kb^D(\Zc^a_y)\to 0$$
is an exact sequence of $\Rb_D$-modules, and the forgetful map 
$\Rb_D\otimes_{\Rb_A}\Kb^A(\Zc^a_{\leq y})\to\Kb^D(\Zc^a_{\leq y})$
is an isomorphism, see Remark 4.4.
\qed

\vskip3mm

\noindent{\sl Proof of $(ii)$ :}
Set $D=A$.
Then 
$$0\to\Kb^A(\Zc^a_{<y})_a
\to\Kb_{a,\leq y}
\to\Kb^A(\Zc^a_{y})_a\to 0\leqno(4.9.1)$$
is an exact sequence of $\CC$-vector spaces.
The 2-nd map is a morphism of right $\Kb_{a,\leq 1}$-modules,
relatively to $\star$.
Let $\Kb_{a,\leq 1}$ act on $\Kb^A(\Zc^a_y)_a$ so that
$$E\star D_{\l,a}=E\otimes(\CC\boxtimes\Lc_\l),$$
see 4.8.
Then (4.9.1) is an exact sequence of
right $\Kb_{a,\leq 1}$-modules by (4.7.1).
Note that $\Kb^A(\Zc^a_y)_a$
is free of rank one over $\Kb_{a,\leq 1}$ by Lemma 4.8.
Therefore
$\Kb_{a,\leq y}$ is free of rank $\ell(y)$
over $\Kb_{a,\leq 1}$.

Fix $i_1\in I$.
Assume that $y\geq s_{i_1}y$. 
Let
$$\theta_{i_1}\,:\,\Kb_{a,\leq s_{i_1}y}
\to\Kb_{a,\leq y}
\to\Kb_{a,\leq y}
\to\Kb(\Zc^a_y)_a$$
be the composition of
the direct image by the natural inclusion,
the map $E\mapsto D_{i_1,a}\star E$,
and the restriction to the open subset 
$\Zc^a_y\subset\Zc^a_{\leq y}.$
Note that $D_{i_1,a}\star E$ still lies in $\Kb_{a,\leq y}$
because $s_{i_1}y\leq y$, see Proposition 4.4.$(i)$.

The map $\theta_{i_1}$ commutes with the right $\Kb_{a,\leq 1}$-action.
We have $\Psi_a(\Hb_{\leq y})\subset\Kb_{a,\leq y}$
by an easy induction on $\ell(y)$. 
Moreover $\Kb_{a,\leq y}$ is free of rank $\ell(y)$
over $\Kb_{a,\leq 1}$. 
Therefore it is sufficient to prove that $\theta_{i_1}$ is surjective.

There are invertible elements 
$R,S\in\Kb_A(\Tc^a)_a$ such that
the restriction of $D_{i_1,a}$ to the open subset 
$\Zc_{s_{i_1}}^a\subset\bar\Zc_{s_{i_1}}^a$
is $(R\boxtimes S)|_{\Zc_{s_{i_1}}^a}$.
The $\CC$-vector space $\Kb^A(\Zc^a_y)_a$ is spanned by 
$\{(R\boxtimes\Lc_\l)|_{\Zc^a_y}\,;\,\l\in X\}$
by Lemma 4.8.$(i)$, because the second projection 
$\Zc^a_{ij}\cap\Zc_y\to\Xc^s_j$ is a vector bundle for each $i,j$.
Hence it is sufficient to check that 
$(R\boxtimes\Lc_\l)|_{\Zc^a_y}$ lies in the image of $\theta_{i_1}$.

Let $\Uc\subset\Tc^3$ be open and such that
$$\Uc\cap(\bar\Zc_{s_i}\times\Tc)=
\Uc\cap(\Zc_{s_i}\times\Tc),$$
$$\Uc\cap(\Tc\times\Zc_{\leq s_iy})=
\Uc\cap(\Tc\times\Zc_{s_iy}),$$
$$\Uc\cap q_{12}^{-1}(\bar\Zc_{s_i})\cap 
q_{23}^{-1}(\Zc_{\leq s_iy})=
q_{12}^{-1}(\Zc_{s_i})\cap 
q_{23}^{-1}(\Zc_{s_iy}),$$
see \cite{L2, 8.3}.
Put $E=(S^{-1}\boxtimes\Lc_\l)|_{\Zc^a_{\leq s_{i_1}y}}$.
Then
$$\theta_{i_1}(E)=
q_{13*}((R\boxtimes S\boxtimes\CC)\otimes_{\Uc^a}
(\CC\boxtimes S^{-1}\boxtimes\Lc_\l)).$$
The schemes $\Uc\cap(\Zc_{s_{i_1}}\times\Tc), 
\Uc\cap(\Tc\times\Zc_{s_{i_1}y})$ 
intersect transversally along $\Zc_y$.
Hence $\theta_{i_1}(E)=(R\boxtimes\Lc_\l)|_{\Zc^a_y}$.
\qed

\vskip2mm

\noindent{\sl Proof of $(iv)$:}
Let ${}^w\Psi_a\,:\,\Hb\to\Kb^A({}^w\Zc^a)_a$ 
be the composition of $\Psi_a$ and the restriction.
The restriction of ${}^w\Psi_a$ to the subring $\Rb_H\subset\Hb$
is the composition of ${}^w\beta_s$, the pull-back by the vector bundle
${}^w\Tc^a\to{}^w\Xc^s$, and the isomorphism
$\Kb_A({}^w\Tc^a)_a\to\Kb^A({}^w\Zc^a_1)_a$ in Lemma 4.8.
Hence it yields a ring isomorphism 
$\Rb_A/I_{w,s}\to\Kb^A({}^w\Zc_1^a)_a.$
Taking the limit over $w$, we get an isomorphism of topological rings
$\hat\Psi_a\,:\,\hat\Rb_a\to\Kb_{a,\leq 1}$.
Set $\hat\Hb_a=\Hb\otimes_{\Rb_A}\hat\Rb_a$, and
$\hat\Hb_{a,\leq y}=\Hb_{\leq y}\otimes_{\Rb_A}\hat\Rb_a$.
By $(ii)$ there is a unique isomorphism of topological vector spaces 
$\hat\Hb_{a,\leq y}\to\Kb_{a,\leq y}$ such that
$x\otimes f\mapsto\Psi_a(x)\star\hat\Psi_a(f)$.
Let $\hat\Psi_a$ denotes also the resulting isomorphism
$\hat\Hb_a\to\Kb_a$.

Let $\Mb$ be a smooth $\Kb_a$-module.
Then $\Psi_a(\Rb_a)(m)=\Kb_{a,\leq 1}(m)$
for any element $m\in\Mb.$
Thus $\Psi_a(\Hb)(m)=\Kb_a(m)$.
In particular $\Psi_a(\Mb)$ is simple if $\Mb$ is simple.
This yields the first part of $(iv)$.
The second part is obvious because
an $\Hb$-isomorphism of smooth $\Kb_a$-modules 
is a linear isomorphism which
commutes with the action of the subring
$\Psi_a(\Hb)\subset\Kb_a$, hence with the action of $\Kb_a$ by density
($\Psi_a(\Hb)\subset\Kb_a$ is dense, 
and the $\Kb_a$ action on a smooth module is continuous). 
We are done.
\qed

\vskip3mm

\noindent{\bf Remark.}
Probably $\bigcap_w\Ker({}^w\b_s)=(t-\tau)$.
Then the natural ring homomorphism
$\Rb_{\tau,\zeta}\to\hat\Rb_a$ would be injective,
and the proof of Theorem 4.9.$(iv)$ 
would imply that $\Ker(\Psi_a)=(t-\tau,q-\zeta)\subset\Hb$.
We will not use this statement.

\head 5. Induced modules\endhead

Sections 5.1 to 5.3 contain definition and basic properties of a set $\Sc_\phi$.
The group $\Kb^A(\Sc^a_\phi)$ and the $\Kb_a$-action on
$\Kb^A(\Sc^a_\phi)$ is given in 5.4.
Sections 5.5, 5.7 contain basic properties of $\Kb^A(\Sc^a_\phi)$.
The induction theorem is proved in 5.6, 5.8.
The arguments are adapted from \cite{L3}.

\subhead 5.1\endsubhead
Let $\gen_{KM}\subset\gen$ be the Lie algebra of $G_{KM}$, see 2.2.
Given a $\sen\len_2$-triple $\phi=\{e,f,h\}\subset\gen_{KM}$,
let $\gen_\phi$ be the linear span of $\{e,f,h\}$. 
It is a reductive Lie subalgebra of $\gen_{KM}$ by \cite{KcP, Proposition 3.4}.
Hence it can be conjugated to $\gen_J$ for some $J\subsetneq I$ by \cite{KcP}.
In particular there is a unique connected algebraic subgroup 
$G_\phi\subset G_{KM}$
with Lie algebra $\gen_\phi$ in the sense of \cite{KcW, 2.11}, by
\cite{KcW, Lemma 2.13 and Proposition 2.14}.
Let $\phi\,:\,\SL_2\to G_\phi$ 
denote also the homomorphism of algebraic groups
whose tangent map at 1 takes 
$(\smallmatrix 0&1\cr0&0\endsmallmatrix)$ to $e$ and
$(\smallmatrix 0&0\cr1&0\endsmallmatrix)$ to $f$.
Set $\phi(z)=\phi(\smallmatrix z&0\cr 0&1/z\endsmallmatrix)$ 
for any $z\in\CC^\times$.
Let $Z_\phi(G)\subset G$ be the centralizer of $G_\phi$.

Set $A_e=\{a\in A\,;\,\ad_ae=e,\,\ad_af=f\}$,
and $C_\phi=\{(\phi(z),z^2)\,;\,z\in\CC^\times_q\}$.
Assume that $C_\phi\subset A_e$.
If $a\in A_e$ we set 
$s_\phi=a\cdot(\phi(\zeta^{-1/2}),\zeta^{-1})\in Z_\phi(G)$,
where $\zeta^{1/2}$ is a square root of $\zeta$. 

From now on we assume that $\phi\subset\genu$,
so that $\gen_\phi\subset\genu$.

\subhead 5.2\endsubhead
Recall that ${}^\flat\Bc_{\leq y^{-1}}\subset{}^w\Xc$
whenever $y\in W^\flat$, $w\in W$, and $w\geq w_\flat y^{-1}$.
See 2.9 for the meaning of the symbol $\flat$, 
and Lemma 2.3.$(ii)$ for the proof of the inclusion.
Fix $k\gg 0$.
Then $p^k$ yields a closed immersion
${}^\flat\Bc_{\leq y^{-1}}\to{}^w\Xc^k.$
We may omit $p^k$ and identify 
${}^\flat\Bc_{\leq y^{-1}}$ with its image in ${}^w\Xc^k.$
%Then $y'(\uen_\ell)\subset\ben$ for all $y'\leq yw_\flat$.
We set $\uen_\ell=\eps^\ell\cdot\uen_\flat$,
with $\ell\gg 0$.

\subhead 5.3\endsubhead
Let $\zen(f)\subset\gen$ be the centralizer of $f$,
$\senu_\phi=\zen(f)\cap\genu+e$,
and $\Scu_\phi\subset\dot\Ncu$ be the inverse image of
$\senu_\phi\cap\Ncu$ by the Springer map. 
The scheme $\Scu_\phi$ is smooth.

Put $\sen_\phi=\zen(f)\cap\ben_\flat+e$.
Note that $e\in\ben_\flat$ because $e\in\gen_\phi$ 
and $\gen_\phi\subset\genu$ by hypothesis.
We have 
$\sen_\phi\cap\Nc=\sen_\phi\cap(\Ncu+\uen_\flat)$,
because
$\ben_\flat\cap\Nc=\Ncu+\uen_\flat$.
In particular $\sen_\phi\cap\Nc$ is closed in $\ben_\flat$.
The $A$-action on $\gen$ restricts to a $A_e$-action on $\sen_\phi$.

Put $\Sc_{\phi,y}=\dot\Nc\cap(\sen_\phi\times{}^\flat\Bc_{y^{-1}}).$
Then $\Sc_{\phi,y}\subset{}^w\Tc.$
Moreover $\Sc_{\phi,y}$ is locally closed.
It is endowed with the induced reduced scheme structure.
The schemes $\Sc_{\phi,\leq y},$ $\Sc_{\phi,<y}$
are defined in the obvious way, and
$\Sc_{\phi}=\bigcup_y\Sc_{\phi,y}$, see 2.3.
Note that $\Sc_\phi$ is contained in the set $\dot\Nc$, 
which is itself naturaly embedded into the scheme $\Tc$, see 2.8.

Put $\Sc_{\phi,y}^{k\ell}=
\{(x+\uen_\ell,\pen)\in{}^w\Tc^{k\ell}\,;\,(x,\pen)\in\Sc_{\phi,y}\}.$
Since $\Sc_{\phi,y}^{k\ell}$,
$\Sc_{\phi,\leq y}^{k\ell}$, and $\Sc_{\phi,<y}^{k\ell}$
are locally closed in ${}^w\Tc^{k\ell}$,
see the lemma below, they 
are naturally subschemes of ${}^w\Tc^{k\ell}$.
They are separated, because they are contained in 
$p^k_\ell\rho^{-1}(\Xc_{\leq w_\flat y^{-1}})$, 
and $k$ is large enough for $p^k$ to yield a closed immersion
$\Xc_{\leq w_\flat y^{-1}}\to{}^w\Xc^k$ by 5.2. 

Let $\g^{k\ell}_y\,:\,\Sc_{\phi,y}^{k\ell}\to\Scu_\phi$ be the map 
$(x+\uen_\ell,\pen)\mapsto(x+\uen_\flat,\pen^!)$.
The group $A_e$ acts on $\Sc_{\phi,y}^{k\ell}$, $\Scu_\phi$
and $\g^{k\ell}$ commutes with the $A_e$-action.
Set also $\g_y=\g^{k\ell}_y\circ p_\ell^k\,:\,\Sc_{\phi,y}\to\Scu_\phi.$

\proclaim{Lemma}
(i)
$\Sc^{k\ell}_{\phi,\leq y}$ is closed in ${}^w\Tc^{k\ell}$,
$\Sc^{k\ell}_{\phi,y}$ is open in $\Sc^{k\ell}_{\phi,\leq y}$.

(ii)
$(\Sc^{k\ell}_{\phi,y})^{\CC^\times_\delta},
(\Sc^{k\ell}_{\phi,<y})^{\CC^\times_\delta}$
are both open and closed in
$(\Sc^{k\ell}_{\phi,\leq y})^{\CC^\times_\delta}$.

(iii) 
$\g^{k\ell}_y$ is a $A_e$-equivariant vector bundle.

(iv) 
$p^{k_2k_1}$ restricts to an isomorphism
$\Sc^{k_2\ell}_{\phi,\leq y}\to\Sc^{k_1\ell}_{\phi,\leq y}$,
and $p_{\ell_2\ell_1}$ to a vector bundle
$\Sc_{\phi,\leq y}^{k\ell_2}\to\Sc_{\phi,\leq y}^{k\ell_1}.$

(v)
$C_\phi$ acts on $\sen_\phi$ fixing $e$, and
$z\mapsto(\phi(z),z^2)$ has $<0$ weights on $\sen_\phi-e$.
\endproclaim

\noindent{\sl Proof:}
Observe that the fixed points subset 
$\Xc^{\CC^\times_\delta}\subset\Xc$ is a disjoint union of 
closed $\Gu$-orbits, see Lemma 2.13.$(ii)$, and that
${}^\flat\Bc_{y^{-1}},{}^\flat\Bc_{<y^{-1}}$ 
are preserved by $\Gu$. 
Thus
$({}^\flat\Bc_{y^{-1}})^{\CC^\times_\delta}$,
$({}^\flat\Bc_{<y^{-1}})^{\CC^\times_\delta}$
are open and closed in
$({}^\flat\Bc_{\leq y^{-1}})^{\CC^\times_\delta}$.
Hence
$(\Sc_{\phi,y})^{\CC^\times_\delta},
(\Sc_{\phi,<y})^{\CC^\times_\delta}$
are open and closed in
$(\Sc_{\phi,\leq y})^{\CC^\times_\delta}$.
This yields $(ii)$, 
because $p^k$ restricts to a closed immersion 
${}^\flat\Bc_{\leq y^{-1}}\to{}^w\Xc^k$ by 5.2.

As for $(i)$, note that 
$\Sc_{\phi,\leq y}^{k\ell}$ is contained in 
$p^k_\ell\rho^{-1}({}^ \flat\Bc_{\leq y^{-1}})$ and that  
$p^k$ restricts to a closed immersion 
$p_\ell\rho^ {-1}({}^\flat\Bc_{\leq y^{-1}})\to{}^w\Tc^{k\ell}$ by 5.2.
Thus, to prove that $\Sc_{\phi,\leq y}^{k\ell}$ is closed in ${}^w\Tc^{k\ell}$
it is enough to check that 
$p_\ell(\dot\Nc\cap(\sen_\phi\times{}^\flat\Bc_{\leq y^{-1}}))$
is closed in $\Tc^\ell$.
This is obvious, as well as the rest of $(i)$.
Claim $(iv)$ is also obvious.
Claim $(v)$ follows from the representation theory of $\slen_2$.

In order to prove $(iii)$
we first prove that $\Sc^{k\ell}_{\phi,y}$ is a smooth scheme.
The set $\Vc_y=\dot\gen\cap(\ben_\flat\times{}^\flat\Bc_{y^{-1}})$
is a closed subset of $\ben_\flat\times{}^\flat\Bc_{y^{-1}}$ 
containing $\Sc_{\phi,y}$.
Set $\Vc_y^{k\ell}=\{(x+\uen_\ell,\pen)\,;\,(x,\pen)\in\Vc_y\}.$
Note that $\Vc_y^{k\ell}$ is a smooth scheme because
the map $\epsilon\,:\,\Vc^{k\ell}_y\to\dot\genu$, 
$(x+\uen_\ell,\pen)\mapsto(x+\uen_\flat,\pen^!)$
is smooth.
Put
$\Vc_{\phi,y}^{k\ell}=\{(x+\uen_\ell,\pen)\in\Vc_{y}^{k\ell}\,;\,
x\in\sen_\phi\}$.
It is a closed subset of $\Vc_y^{k\ell}$ such that
$$\Sc_{\phi,y}^{k\ell}=\{(x+\uen_\ell,\pen)\in\Vc_{\phi,y}^{k\ell}\,;\,
x\in\Ncu+\uen_\flat\}.$$
Let $\cenu$ be the canonical Cartan Lie algebra of $\genu$.
Recall that $\cenu=\pen/[\pen,\pen]$ for each $\pen\in\Bcu$ and that
$\cenu$ is independent on the choice of $\pen$, up to a canonical
isomorphism, see \cite{CG, Lemma 3.1.26}.
The map $\delta\,:\,\dot\genu\to\cenu$, 
$(x,\pen)\mapsto x+[\pen,\pen]$ is smooth.
To prove that $\Sc_{\phi,y}^{k\ell}$ is smooth
it is sufficient to check that the restriction of
$\delta\circ\epsilon$ to $\Vc_{\phi,y}^{k\ell}$ is smooth because 
$\Sc_{\phi,y}^{k\ell}=
\Vc^{k\ell}_{\phi,y}\cap(\delta\circ\epsilon)^{-1}(0)$.
The element $g\in B_\flat$ acts on $\Vc_y^{k\ell}$ so that 
$g\cdot (x+\uen_\ell,\pen)=(\ad_gx+\uen_\ell,\ad_g\pen).$
The subgroup $U_\ell\subset B_\flat$ is normal 
and $\ad_{U_\ell}\ben_{y^{-1}}=\ben_{y^{-1}}.$
Hence $U_\ell$ acts trivially on ${}^\flat\Bc_{y^{-1}}$.
Thus $B_\flat/U_\ell$ acts on $\Vc_y^{k\ell}$.
The map $\delta\circ\epsilon$ is smooth
because it is the composition of two smooth maps.
It is constant along the $B_\flat/U_\ell$-orbits.
Hence it is sufficient to prove that the map
$$\Vc_{\phi,y}^{k\ell}\times B_\flat/U_\ell\to\Vc_y^{k\ell},\,
(x+\uen_\ell,\pen;g)\mapsto g\cdot(x+\uen_\ell,\pen)
\leqno(5.3.1)$$
is smooth.
The group $\SL_2$ acts on $\ben_\flat$ via $\phi$.
The subspace $\uen_\ell\subset\ben_\flat$ is a $\SL_2$-submodule.
The representation theory of $\SL_2$ implies that 
$$\ben_\flat/\uen_\ell=
[\ben_\flat/\uen_\ell,e]\oplus
(\zen(f)\cap\ben_\flat+\uen_\ell)/\uen_\ell,$$
i.e. $(\sen_\phi+\uen_\ell)/\uen_\ell$ is a transversal slice to the
$B_\flat/U_\ell$-orbit of $e$ in $\ben_\flat/\uen_\ell$ near $e$.
The $\CC^\times$-action in $(v)$ descends to 
$(\sen_\phi+\uen_\ell)/\uen_\ell$.
Hence, by $(v)$, the space $(\sen_\phi+\uen_\ell)/\uen_\ell$ 
is a transversal slice to the
$B_\flat/U_\ell$-orbit of $e$ in $\ben_\flat/\uen_\ell$.
In other words the map (5.3.1) is smooth.
Hence $\Sc^{k\ell}_{\phi,y}$ is smooth.

Note that $\epsilon$ is a $A$-equivariant vector bundle, 
with $\CC^\times_\delta$ acting trivially on the base.
The fiber of $\epsilon\circ p^k_\ell\,:\,\Vc_y\to\dot\genu$ 
containing the element $(x_1,\pen_1)$ is
$$\{g\cdot(x_1+x,\pen_1);
(g,x)\in U_\flat\times(\uen_\flat\cap\pen_1)\}.$$
Hence $\CC^\times_\delta$ acts with 
positive weights on the fibers of $\epsilon$.
Since $\Sc^{k\ell}_{\phi,y}\subset\Vc^{k\ell}_y$
is a smooth closed $\CC^\times_\delta$-subscheme,
the restriction of $\epsilon$ to $\Sc^{k\ell}_{\phi,y}$
is a vector bundle over $\epsilon(\Sc^{k\ell}_{\phi,y})$
by the following fact, which follows from \cite{BH, Theorem 9.1} : 
if $f:V\to B$ is a vector bundle over a smooth variety, 
with a fiber preserving linear $\CC^\times$-action with positive weights, 
and $M\subset V$
is a $\CC^\times$-stable smooth closed subvariety, 
then $M$ is a sub-bundle of $V$ restricted to $f(M)$.
Finally $\g^{k\ell}_y$
is the restriction of $\epsilon$ to $\Sc^{k\ell}_{\phi,y}$, and
$\epsilon(\Sc^{k\ell}_{\phi,y})=\Scu_\phi.$
\qed

\proclaim{Corollary}
Assume that $a\in A_e$.

(i)
$\Sc^a_{\phi,\leq y}$ is closed in ${}^w\Tc^a$,
$\Sc^a_{\phi,y}$ is open in $\Sc^a_{\phi,\leq y}$.

(ii)
$(\Sc^a_{\phi,y})^{\CC^\times_\delta},
(\Sc^a_{\phi,<y})^{\CC^\times_\delta}$
are both open and closed in 
$(\Sc^a_{\phi,\leq y})^{\CC^\times_\delta}$.

(iii) 
$\g_y$ restricts to an $A_e$-equivariant vector bundle
$\g_{y,a}\,:\,\Sc^a_{\phi,y}\to\Scu^a_\phi$.
\endproclaim

\subhead 5.4\endsubhead
From now on $D\subset A_e$ is a closed subgroup containing $a$.
For any $y\in W^\flat$ the scheme $\Sc^a_{\phi,\leq y}$ is of finite type,
hence $\Kb^{D}(\Sc^a_{\phi,\leq y})$ makes sense.
Given $y',y''\in W^\flat$, $w,w'\in W$ such that
$y''$, $w'$ are large enough and $w\geq w_\flat{y'}^{-1}$,
the 1-st projection is a proper map
${}^w\Zc_{\leq y}\cap
({}^{w'}\Tc\times\Sc_{\phi,\leq y'})
\to\Sc_{\phi,\leq y''}$
(as for Lemma 4.2.$(iii)$).
Thus there is a $\star$-product
$$\Kb^{D}({}^{w}\Zc^a_{\leq y})
\times\Kb^{D}(\Sc^a_{\phi,\leq y'})\to
\Kb^{D}(\Sc^a_{\phi,\leq y''})$$
relative to ${}^{w'}\Tc^a\times{}^w\Tc^a$.
If $y\leq y'$ the closed immersion 
$\Sc^a_{\phi,\leq y}\to\Sc^a_{\phi,\leq y'}$
gives a map
$h_{yy'}\,:\,\Kb^{D}(\Sc^a_{\phi,\leq y})\to\Kb^{D}(\Sc^a_{\phi,\leq y'}).$
We set 
$$\Kb^{D}(\Sc^a_\phi)=\ind_y(\Kb^{D}(\Sc^a_{\phi,\leq y}),h_{yy'}).$$
The ring $\Kb^{D}(\Zc^a)$ acts on $\Kb^{D}(\Sc^a_\phi)$ 
because the map $h_{yy'}$ commutes with $\star$,
see (A.3.1).
This representation is smooth.
See 4.9 for the definition of a smooth module.
Using $\Psi_a$ we get the following.

\proclaim{Proposition}
(i)
There is a unique $\Hb$-action on 
$\Kb^{D}(\Sc^a_\phi)_a$ such that
$-q^{-1}(t_i+1)$ takes $E$ to $D_{i,a}\star E$, 
and $x_\l$ takes $E$ to $D_{\l,a}\star E$.
 
(ii)
$\Hbu$ takes $\Kb^D(\Sc^a_{\phi,1})_a$ into itself.
\endproclaim

\proclaim{5.5. Lemma}
(i)
$\Kb_{1,\top}(\Scu^a_\phi)=0$,
and the comparison map
$\Kb_0(\Scu^a_\phi)\to\Kb_{\top}(\Scu^a_\phi)$ is invertible. 
The same holds for $\Sc_\phi^a$.

(ii) 
If $C_\phi\subset D$ and $D$ is generated by $D^\circ\cup\{a\}$, 
then $\Kb^{D}(\Scu^a_\phi)$ is a free $\Rb_{D}$-module, 
the forgetful map
$\Kb^{D}_{\top}(\Scu^a_\phi)_a\to\Kb^{\la a\ra}_{\top}(\Scu^a_\phi)_a$ 
and the comparison map
$\Kb^{D}(\Scu^a_\phi)\to\Kb^{D}_{\top}(\Scu^a_\phi)$
are invertible, 
and $\dim\Kb^{D}(\Scu^a_\phi)_a=\dim\Kb_\top(\Scu^a_\phi)$.
The same holds for $\Sc_\phi^a$.

(iii)
The natural sequence 
$\Kb^{D}(\Sc^a_{\phi,<y}){\buildrel\nu\over\to}
\Kb^{D}(\Sc^a_{\phi,\leq y})
\to\Kb^{D}(\Sc^a_{\phi,y})\to 0$
is exact.
If $C_\phi,\CC^\times_\delta\subset D$ 
and $D$ is generated by $D^\circ\cup\{a\}$,
then $\nu$ is injective.
\endproclaim

\noindent{\sl Proof:}
We say that a variety $\Xc$ has property $(S_\QQ)$
if $\Kb_{1,\top}(\Xc)=0$,
and the comparison map
$\Kb_0(\Xc)\to\Kb_{\top}(\Xc)$ is an isomorphism.
Recall that a $\CC^\times$-action on $\Xc$ is a contraction to
the subset $\Yc\subset\Xc$ if
$\lim_{z\to 0}z\cdot x\in\Yc$ for all $x\in\Xc$, or
$\lim_{z\to\infty}z\cdot x\in\Yc$ for all $x\in\Xc$.
Using the Bialynicki-Birula decomposition we get the following, 
see \cite{DLP, Section 1}. 

\vskip2mm

\noindent{\sl Claim 1.
Let $\Xc$ be a smooth quasi-projective $\CC^\times$-variety.
Assume that the $\CC^\times$-action is a contraction to a compact subset
and that $\Xc^{\CC^\times}$ has property $(S_\QQ)$,
then $\Xc$ has property $(S_\QQ)$.
}

\vskip2mm

\noindent
Recall that $a\in A_e$ and $e\in\senu_\phi$.
Set $\Bcu_e^s=\Bcu_e\cap\Bcu^s$.
Then $\Scu^a_\phi\cap\Bcu_e=\Bcu_e^s$, 
where $\Bcu_e$ is viewed as a subset in the fiber at $e$
of the Springer map $\Scu_\phi\to\senu_\phi$.
The $C_\phi$-action on $\Scu^a_\phi$ preserves $\Bcu_e^s$, and
$(\Bcu_e^s)^{C_\phi}=(\Bcu_e^{s_\phi})^{C_\phi}$
with $s_\phi$ as in 5.1. 
Since $s_\phi$ and $e$ commute, 
each connected component
of $\Bcu^{s_\phi}_e$ is isomorphic to the variety of Borel subgroups of
$\Gu^{s_\phi}$ containing $e$.
Hence $(\Bcu^{s_\phi}_e)^{C_\phi}$ has property $(S_\QQ)$ by \cite{DLP}.
Thus $\Scu^a_\phi$ has property $(S_\QQ)$, because the $C_\phi$-action
on $\Scu^a_\phi$ is a contraction to $\Bcu_e^s$ 
and $(\Scu^a_\phi)^{C_\phi}=(\Bcu^{s_\phi}_e)^{C_\phi}$.
The proof for $\Sc^a_\phi$ is identical because $G^{s_\phi}$
is a reductive group.
$(i)$ is proved.

Given a linear group $H$, we say that a variety $\Xc$ has property $(T_H)$
if $\Kb^{H'}_{1,\top}(\Xc)=0$,
$\Kb^{H'}_{\top}(\Xc)$ is a free $\Rb_{H'}$-module,
the comparison map
$\Kb^{H'}(\Xc)\to\Kb^{H'}_{\top}(\Xc)$ is an isomorphism,
and the forgetful map
$\Kb^H(\Xc)\otimes_{\Rb_H}\Rb_{H'}\to\Kb^{H'}(\Xc)$ is an isomorphism,
for any closed subgroup $H'\subset H$.
Observe that if $\Xc$ has the property $(S_\QQ)$ and $H$ acts trivially
on $\Xc$, then $\Xc$ has also the property $(T_H)$,
because $\Kb^{H'}_0(\Xc)=\Rb_{H'}\otimes\Kb_0(\Xc)$ and 
$\Kb^{H'}_{i,\top}(\Xc)=\Rb_{H'}\otimes\Kb_{i,\top}(\Xc)$ for any $H'$.
Using the Bialynicki-Birula decomposition we get the following
equivariant annalogue of Claim 1, see \cite{N, 7.1} for instance.

\vskip2mm

\noindent{\sl Claim 2.
Given a diagonalizable group $H$,
a one-parameter subgroup $\CC^\times\subset H$, and
a smooth quasi-projective $H$-variety $\Xc$ such that 
the $\CC^\times$-action on $\Xc$ is a contraction to a compact subset,
if $\Xc^{\CC^\times}$ has property $(T_H)$
then $\Xc$ has property $(T_H)$.
}

\vskip2mm

\noindent
Fix a one-parameter subgroup $\psi\,:\,\CC^\times\to D^\circ$ 
in general position so that 
$(\Scu^a_\phi)^{\CC^\times_\psi}=(\Scu^a_\phi)^{D^\circ}$,
and such that the $\CC^\times_\psi$-action on $\Scu^a_\phi$ 
is a contraction to $\Bcu^s_e$
(this is possible by Lemma 5.3.$(v)$ because $C_\phi\subset D^\circ$).
Then $D$ acts trivially on $(\Scu^a_\phi)^{\CC^\times_\psi}$, 
because $D$ is generated by $D^\circ\cup\{a\}$.
On the other hand $(\Scu^a_\phi)^{\CC^\times_\psi}=
(\Bcu_e^{s_\phi})^{\CC^\times_\psi}$.
Hence $(\Scu^a_\phi)^{\CC^\times_\psi}$ 
has the property $(S_\QQ)$ by \cite{DLP}. 
Thus it has the property $(T_{D})$ by the remark above. 
Hence $\Scu^a_\phi$ has the property $(T_{D})$. 
The proof for $\Sc^a_\phi$ is identical.

The first part of $(iii)$ follows from the 
localization long exact sequence and Corollary 5.3.$(i)$.
The second one follows by induction on $y$ from 
$(ii)$, Corollary 5.3.$(ii),(iii)$,
and the following general fact.

\vskip2mm

\noindent{\sl Claim 3.
Let $H$ be a diagonalizable group,
$\Xc$ a $H$-variety, $\Yc\subset\Xc$ a closed subset
preserved by the $H$-action, 
$\CC^\times\subset H$ a subgroup such that
$\Yc^{\CC^\times},\Xc^{\CC^\times}\setminus\Yc^{\CC^\times}
\subset\Xc^{\CC^\times}$
are open and closed.
If the $\Rb_H$-module $\Kb^H(\Yc)$ is free then the natural map
$\Kb^{H}(\Yc)\to\Kb^{H}(\Xc)$ is injective. 
}

\vskip2mm

\noindent
Let $\nu$ be the map $\Kb^{H}(\Yc)\to\Kb^{H}(\Xc).$
In any connected component $H_i\subset H$ we can fix an element $h_i$ such that
$\Xc^{h_i}\subset\Xc^{\CC^\times}$.
Using the concentration theorem and the natural 
commutative square
$$\matrix
\Kb_1^H(\Xc^{h_i}\setminus\Yc^{h_i})&{\buildrel 0\over\lra}&
\Kb^H(\Yc^{h_i})\cr
\bda&&\bda\cr
\Kb_1^H(\Xc\setminus\Yc)&\lra&\Kb^H(\Yc),\cr
\endmatrix$$
we prove that for any $x\in\Ker(\nu)$ and for any $i$,
there is an element $u_i\in\Rb_H\setminus\Jb_{h_i}$ such that $u_i\cdot x=0$. 
Since the $\Rb_H$-module $\Kb^H(\Yc)$ is free,
we get $x=0$. Claim 3 is proved.
\qed

\vskip3mm

From now on we assume that
$C_\phi,\CC^\times_\delta\subset D$ 
and that $D$ is generated by $D^\circ\cup\{a\}$.

\subhead 5.6\endsubhead
The group $\Kb^D(\Scu_\phi)$ is endowed with
the $\Hbu$-action in \cite{L3, Proposition 4.4}. 
The ring $\Hbu$ acts also on $\Kb^D(\Sc^a_{\phi,1})_a$
by Proposition 5.4.$(ii)$.
Note that $\Scu_\phi\simeq\Sc^{k0}_{\phi,1}$ because $k$ is large.
In particular $\Scu_\phi$ is identified with a closed subset in 
${}^{w_\flat}\Tc^{k0}$.
Let 
$$\tilde\cb_a\,:\,
\Kb^D(\Scu_\phi)_a\to\Kb^D(\Scu^a_{\phi})_a$$
be the concentration map relative to
${}^{w_\flat}\Tc^{k0}$,
i.e. the map (A.3.3) with $\Zc\subset\Tc_1\times\Tc_2$ 
equal to $\Scu_{\phi}$,
$\Tc_1={}^{w_\flat}\Tc^{k0}$, and $\Tc_2=\{\bullet\}$.
Set $\tilde\Theta_\flat=(\g_{1,a})^*\circ\tilde\cb_a$,
where $\g_{1,a}$ is the value of $\g_{y,a}$ at $y=1$, see Corollary 5.3.

\proclaim{Lemma}
$\tilde\Theta_\flat$ is an isomorphism of $\Hbu$-modules
$\Kb^D(\Scu_\phi)_a\to\Kb^D(\Sc^a_{\phi,1})_a$.
\endproclaim

\noindent{\sl Proof:}
$\tilde\Theta_\flat$ is invertible because $\g_{1,a}$ is an vector bundle,
see Corollary 5.3.
Let $\Psiu\,:\,\Hbu\to\Kb^{\Gu\times\CC_q^\times}(\Zcu)$ 
be the ring homomorphism in \cite{L2, Theorem 7.25}.
The $\Hbu$-action on $\Kb^D(\Scu_\phi)$ is the composition
of $\Psiu$ and the $\star$-product
$$\Kb^D(\Zcu)\times\Kb^D(\Scu_\phi)\to\Kb^D(\Scu_\phi)$$
relative to $\dot\Ncu^2$.
The $\Hbu$-action on $\Kb^D(\Sc^a_{\phi,1})$ is the composition
of $\Psi_a$ and the $\star$-product
$$\Kb^{D}({}^{w_\flat}\Zc^a_{\leq w_\flat})
\times\Kb^{D}(\Sc^a_{\phi,1})\to
\Kb^{D}(\Sc^a_{\phi,1})$$
relative to $(\Tc^a)^2$.
Let $\Theta_\flat$, $\cb_a$, $\Gamma$ and $p_a$ be as in 4.9.
We must prove that the diagram
$$\matrix
\Kb^{\Gu\times\CC_q^\times}(\Zcu)&\times&\Kb^{D}(\Scu_\phi)_a&
\lra&\Kb^{D}(\Scu_\phi)_a\cr
{\ss\Theta_\flat}\bda&&{\ss\tilde\Theta_\flat}\bda
&&{\ss\tilde\Theta_\flat}\bda\cr
\Kb^{D}({}^{w_\flat}\Zc_{\leq w_\flat}^a)_a&\times
&\Kb^{D}(\Sc_{\phi,1}^a)_a
&\lra&\Kb^{D}(\Sc_{\phi,1}^a)_a
\endmatrix$$
is commutative.
Recall that
$${}^{w_\flat}\Tc^{k0}\simeq{}^{w_\flat}\Pc^{\flat,k}\times_\Gu\dot\Ncu,\quad
{}^{w_\flat}\Zc^{k0}_{\leq w_\flat}
\simeq{}^{w_\flat}\Pc^{\flat,k}\times_\Gu\Zcu.$$
We can complete the diagram as follows
$$\matrix
\Kb^{\Gu\times\CC_q^\times}(\Zcu)&\times&\Kb^{D}(\Scu_\phi)_a&
\lra&\Kb^{D}(\Scu_\phi)_a\cr
{\ss\Gamma}\bda&&\Arrowvert&&\Arrowvert\cr
\Kb^{D}({}^{w_\flat}\Zc^{k0}_{\leq w_\flat})_a&\times&
\Kb^{D}(\Scu_\phi)_a&
\lra&\Kb^{D}(\Scu_\phi)_a\cr
{\ss\cb_a}\bda&&{\ss\tilde\cb_a}\bda&&{\ss\tilde\cb_a}\bda\cr
\Kb^{D}(({}^{w_\flat}\Zc^0_{\leq w_\flat})^a)_a&\times
&\Kb^{D}(\Scu_\phi^a)_a&
\lra&\Kb^{D}(\Scu_\phi^a)_a\cr
{\ss{}^\sharp p_a}\bda&&{\ss p^*_a}\bda&&{\ss p^*_a}\bda\cr
\Kb^{D}({}^{w_\flat}\Zc_{\leq w_\flat}^a)_a&\times
&\Kb^{D}(\Sc_{\phi,1}^a)_a
&\lra&\Kb^{D}(\Sc_{\phi,1}^a)_a.
\endmatrix$$
Note that $\g_{1,a}$ equals $p_a$.
The intermediate horizontal maps are the $\star$-products
relative to ${}^{w_\flat}\Pc^{\flat,k}\times_\Gu\dot\Ncu^2$ and
$((\Tc^0)^a)^2$.

We have $\Gamma(E)\star F=E\star F$ by (A.4.2),
$\cb_a(E)\star\tilde\cb_a(F)=\tilde\cb_a(E\star F)$ by (A.4.3),
and ${}^\sharp p_a(E)\star p_a^*(F)=p_a^*(E\star F)$ by (A.3.2).
Hence the diagram is commutative.
\qed

\subhead 5.7\endsubhead
Fix $y\in W^\flat$ and $i\in I$ with $s_iy<y$.
Hence $s_iy\in W^\flat$.
Let 
$$\sigma_i\,:\,\Kb^{D}(\Sc^a_{\phi,\leq s_iy})_a
\to\Kb^{D}(\Sc^a_{\phi,\leq y})_a
\to\Kb^{D}(\Sc^a_{\phi,\leq y})_a
\to\Kb^{D}(\Sc^a_{\phi,y})_a$$
be the composition of the direct image, 
of the map $E\mapsto D_{i,a}\star E$,
and of the restriction.

\proclaim{Lemma}
$\sigma_i$ is surjective.
\endproclaim

\noindent{\sl Proof:}
We have $\Kb^D(\Sc^a_{\phi,<s_iy})_a\subset\Ker(\sigma_i)$,
because $\Psi_a(t_i)\star\Kb^D(\Sc^a_{\phi,\leq y'})_a\subset
\Kb^D(\Sc^a_{\phi,\leq y'})_a$ whenever $s_iy'w_\flat<y'w_\flat$.
The map $\sigma_i$ factorizes uniquely through a map
$\sigma'_i\,:\,\Kb^{D}(\Sc^a_{\phi,s_iy})_a\to\Kb^{D}(\Sc^a_{\phi,y})_a$
by Lemma 5.5.$(iii)$.
We have the map
$$(\g_{y,a})^*\,:\,
\Kb^D(\Scu_\phi^a)_a\to\Kb^D(\Sc_{\phi,y}^a)_a.$$
It suffices to prove that $\sigma'_i\circ(\g_{s_iy,a})^*$ and $(\g_{y,a})^*$
have the same image, because
$(\g_{s_iy,a})^*$ and $(\g_{y,a})^*$ are invertible by Corollary 5.3.
Fix $w,w'$ with $w\geq y^{-1}s_i$ and $w'$ large enough.
Hence $\Sc_{\phi,\leq s_iy}\subset{}^w\Tc$ and 
$\Sc_{\phi,\leq y}\subset{}^{w'}\Tc$.
Fix an open subset $\Uc\subset{}^{w'}\Tc\times{}^w\Tc$ such that 
$$\matrix
\Uc\cap\Zc_{s_i}=\Uc\cap\bar\Zc_{s_i},\hfill\cr
{}^w\Zc_{s_i}\cap(\Tc\times\Sc_{\phi,s_iy})\subset\Uc\cap\Zc_{s_i},\hfill\cr
\Uc\cap(\Tc\times\Sc_{\phi,\leq s_iy})=
\Uc\cap(\Tc\times\Sc_{\phi,s_iy}).\hfill
\endmatrix$$
The map $\pen\mapsto\pen^!$ is $B_\flat$-equivariant and takes
$\ben_w$ to $\ben$ for all $w\in {}^\flat W$.
Hence $\pen_i^!=\pen^!$ whenever 
$\pen_i=\ad_g(\ben_{y^{-1}s_i})$,
$\pen=\ad_g(\ben_{y^{-1}})$, and $g\in B_\flat$.
Thus we have the commutative diagram 
$$\matrix
{}^w\bar\Zc_{s_i}^a\cap(\Tc^a\times\Sc_{\phi,\leq s_iy}^a)
&{\buildrel\iota\over\lla}&
{}^w\Zc^a_{s_i}\cap(\Tc^a\times\Sc^a_{\phi,s_iy})
&{\buildrel q_2\over\lra}&\Sc^a_{\phi,s_iy}\cr
\bda&&{\ss q_1}\bda&&\bda{\ss\g_{s_iy}}\cr
\Sc^a_{\phi,\leq y}&\lla&\Sc^a_{\phi,y}&
{\buildrel\g_y\over\lra}&\Scu^a_\phi.
\endmatrix\leqno(5.7.1)$$
Given an equivariant sheaf $E$ on $\Scu_\phi^a,$
pick $F_i\in\Kb^{D}(\Sc^a_{\phi,\leq s_iy})_a$ 
whose restriction to $\Sc^a_{\phi,s_iy}$ is $(\g_{s_iy,a})^*(E)$.
Let $F\in\Kb^{D}(\Sc^a_{\phi,\leq y})_a$ be the direct image of $F_i$.
As in 4.9 we fix invertible elements 
$R_i\in\Kb_{D}({}^{w'}\Tc^a)_a$ and  
$S_i\in\Kb_{D}({}^{w}\Tc^a)_a$ such that
the restriction of $D_{i,a}$ to ${}^{w}\Zc_{s_i}^a$ is
$(R_i\boxtimes S_i)|_{{}^{w}\Zc_{s_i}^a}$.
Since the left square in (5.7.1) is Cartesian we have 
$$\matrix
\sigma_i(S_i^{-1}\otimes F_i)
&=q_{1*}\iota^*
(D_{i,a}\otimes(\CC\boxtimes(S_i^{-1}\otimes F)))\hfill\cr
&=q_{1*}((R_i\boxtimes S_i)\otimes_{\Uc^a}
(\CC\boxtimes(S_i^{-1}\otimes(\g_{s_iy,a})^*(E)))).
\hfill\cr
\endmatrix$$
The map $q_1$ in (5.7.1) is an isomorphism.
The varieties 
$\Uc\cap\Zc^a_{s_i}$ and 
$\Uc\cap(\Tc\times\Sc_{\phi,s_iy})^a$
are smooth and intersect transversally along $\Sc^a_{\phi,y}$.
Hence, 
$\sigma_i(S_i^{-1}\otimes F_i)=
q_{1*}(R_i\boxtimes(\g_{s_iy,a})^*(E))|_{\Sc^a_{\phi,y}}$.
We are done.
\qed

\subhead 5.8\endsubhead
Let $\Upsilon\,:\,\Hb_{\tau,\zeta}\otimes_{\Hbu_\zeta}
\Kb^{D}(\Scu_\phi)_a\to\Kb^{D}(\Sc^a_\phi)_a$
be the linear map such that
$x\otimes y\mapsto\Psi_a(x)\star\tilde\Theta_\flat(y).$

\proclaim{Theorem}
(i)
$\Upsilon$ is invertible.

(ii)
The natural map 
$\Kb^{D}(\Sc^a_\phi)_a\to\Kb^{\la a\ra}_\top(\Sc^a_\phi)_a$ 
is an isomorphism.
The same holds for $\Scu_\phi^a$.
\endproclaim

\noindent{\sl Proof:}
The map $\Upsilon$ is a $\Hb$-homomorphism by Lemma 5.6.
Let us prove that it is invertible. 
For any $y\in W^\flat$ the subspace 
$\Hb_{\leq yw_\flat}\subset\Hb_{\tau,\zeta}$
is a free right $\Hbu_\zeta$-submodule of rank $\ell(y)$, see 3.1.
We claim that $\Upsilon$ restricts to a surjective map
$$\Upsilon_y\,:\,\Hb_{\leq yw_\flat}\otimes_{\Hbu_\zeta}
\Kb^{D}(\Scu_\phi)_a
\to\Kb^{D}(\Sc^a_{\phi,\leq y})_a.$$
The existence and invertibility of $\Upsilon_1$
follows from Lemma 5.6.
Assume that $\ell(y)>0$ and that $\Upsilon_{y'}$ is defined
and surjective for all $y'\in W^\flat$ such that $y'< y$. 
We can find $i\in I$ such that $\ell(s_iy)=\ell(y)-1$.
Then $s_iy\in W^\flat$.
We get
$$\Upsilon_y\bigl(
\Hb_{\leq yw_\flat}\otimes_{\Hbu_\zeta}\Kb^{D}(\Scu_\phi)_a\bigr)
=\{1,\Psi_a(t_i)\}\star\Upsilon_{s_iy}
\bigl(\Hb_{\leq s_iyw_\flat}
\otimes_{\Hbu_\zeta}\Kb^{D}(\Scu_\phi)_a\bigr)=$$
$$=\{1,\Psi_a(t_i)\}\star\Kb^{D}(\Sc^a_{\phi,\leq s_iy})_a
=\Kb^{D}(\Sc^a_{\phi,\leq y})_a,$$
where the second equality comes from the surjectivity of $\Upsilon_{s_iy}$, 
and the third from Lemma 5.7.
The claim is proved. 

The map $\Upsilon_y$ is invertible because
$\Hb_{\leq yw_\flat}\otimes_{\Hbu_\zeta}\Kb^{D}(\Scu_\phi)_a$
and $\Kb^{D}(\Sc^a_{\phi,\leq y})_a$ are finite dimensional vector spaces
of the same dimension by Lemma 5.5.$(iii)$ and Corollary 5.3.$(iii)$.
Claim $(i)$ is proved.

Claim $(ii)$ follows from Lemma 5.5.$(ii)$. 
\qed

\head 6. The regular case\endhead

In 6.1 we give a ring homomorphism 
$\Kb_a\to\Extb(L_a,L_a).$
It is injective with a dense image.
The standard modules are introduced in 6.2.
In the regular case, 
the simple modules of $\Extb(L_a,L_a)$ are given in 6.3.
In 6.4-7 we compute the Jordan-H\"older 
multiplicities of the induced modules.

\subhead 6.1\endsubhead
The schemes $\dot\Nc^a_i$ in 2.13 are smooth of finite type by Lemma 2.13.
They are preserved by the $G^s$-action.
The Springer map restricts to a $G^{s}$-equivariant proper map
$p\,:\,\dot\Nc^a_i\to\Nc^a$.
Set $L_{a,i}=p_!\CC_{\dot\Nc^a_i}\in\Db_{G^s}(\Nc^a)$.
By the Beilinson-Deligne-Gabber decomposition theorem we have 
$L_{a,i}\dot{=}\bigoplus_{\chi}\Lb_{a,\chi,i}\otimes S_{a,\chi},$
where $S_{a,\chi}$ are 
simple $G^s$-equivariant perverse sheaves on $\Nc^a$
and $\Lb_{a,\chi,i}$ are finite dimensional $\ZZ$-graded vector spaces.
Set $\Lb_{a,\chi}=\bigoplus_i\Lb_{a,\chi,i}$, and
$L_a\dot{=}\bigoplus_{\chi}\Lb_{a,\chi}\otimes S_{a,\chi}.$
Note that $L_a$ is not well-defined as an object of
$\Db_{G^s}(\Nc^a)$ because $\Xi_s$ is infinite.
We only use it as a convenient notation.
We shall write $\{\chi\}$ for the set of labels $\chi$ parametrizing
the non-zero $\ZZ$-graded vector spaces $\Lb_{a,\chi}$.
Observe that $\Lb_{a,\chi,i}$ may be $\{0\}$ for some $i$,
although $\Lb_{a,\chi}\neq\{0\}$ by convention. 

Consider the topological ring
$$\Extb(L_a,L_a):=\prod_i\bigoplus_j\Extb_{\Db(\Nc^a)}(L_{a,i},L_{a,j}).$$
The product is opposite to the Yoneda product,
and the topology is the product topology.
Recall that we have
an injective continuous ring homomorphism 
$$\Kb_a\to\prod_j\bigoplus_i\Kb^A(\Zc^a_{ij})_a,\leqno(6.1.1)$$
see Remark 4.4.
Its image contains $\bigoplus_{i,j\in S}\Kb^A(\Zc^a_{ij})_a$
for all finite subset $S\subset\Xi_s$.
For each $S$ we have the map 
$$\Phi_S\,:\,\bigoplus_{i,j\in S}\Kb^A(\Zc^a_{ij})_a\to
\bigoplus_{i,j\in S}\Extb_{\Db(\Nc^a)}(L_{a,i},L_{a,j})$$
as in A.2, because
$\Zc^a_{ij}=\{(x,\pen;x',\pen')\in\dot\Nc^a_i\times\dot\Nc^a_j\,;\,x=x'\}$
by Proposition 4.3.
Composing $\Phi_S$ with the direct image by the flip 
$\Zc_{ij}^a\to\Zc_{ji}^a$ and taking the limit, 
we get a continuous ring homomorphism
$\Phi_a\,:\,\Kb_a\to\Extb(L_a,L_a).$

The ring $\Extb(L_a,L_a)$ is the sum 
$nilp\oplus\bigoplus_\chi\Endb(\Lb_{a,\chi})$,
where $nilp$ is an ideal consisting of nilpotent elements.
The topology on $\bigoplus_\chi\Endb(\Lb_{a,\chi})$
is the finite topology, because for each $i$ there is a finite
number of parameters $\chi$ such that $\Lb_{a,\chi,i}\neq\{0\}.$ 
In particular $\Lb_{a,\chi}$ is endowed
with a smooth action of $\Extb(L_a,L_a)$.

\proclaim{Lemma} 
(i)
$\Phi_a$ is injective, and $\Phi_a(\Kb_a)$  
contains $\bigoplus_{i,j\in S}\Extb_{\Db(\Nc^a)}(L_{a,i},L_{a,j})$
for all $S$.

(ii)
The set of smooth and simple $\Kb_a$-modules is
$\{\Phi_a^\bullet(\Lb_{a,\chi})\}.$
\endproclaim

\noindent{\sl Proof:}
The map (6.1.1) is injective.
Hence to prove $(i)$ it is sufficient to prove that 
$\Phi_S$ is invertible for all $S$.
For any $y\in W$ the second projection 
$\Zc^a_y\cap\Zc^a_{ij}\to\Xc^s_j$ 
is an vector bundle by Lemma 4.1.$(iv)$ and Lemma 4.2.$(i)$.
Moreover the Chern character $\Kb_0(\Xc^s_j)\to\Hb_*(\Xc^s_j)$ is an 
isomorphism because $\Xc^s_j$ is isomorphic 
to the flag variety of the group $G^s$.
Hence the bivariant Riemann-Roch map
$\R\R\,:\,\Kb_0(\Zc^a_{ij})\to\Hb_*(\Zc^a_{ij})$ 
relative to $\dot\Nc_i^a\times\dot\Nc_j^a$ is invertible by 
\cite{CG, Theorem 5.9.19}.
The rest is standard, see A.2.

The set of smooth and simple $\Extb(L_a,L_a)$-modules is
$\{\Lb_{a,\chi}\},$ see \cite{J} and Appendix B.
Thus $(ii)$ follows from $(i)$ as in the proof of Theorem 4.9.$(iv)$.
\qed

\subhead 6.2\endsubhead
Given $e\in\Nc^a$ we set $\Bc^s_e=\Bc_e\cap\Bc^s$.
Note that $\Bc^s_e\neq\emptyset$ because $e$ is a nil-element, and
the $\CC^\times$-fixed point subset of a projective variety is non-empty. 
Thus $\Bc^s_e$ is a disjoint union of (reduced) projective 
schemes of finite type.
For any $y\in W$ we set 
$\Bc^s_{e,\leq y}=\Bc^s_e\cap\Xc_{\leq y^{-1}}.$
It is a projective scheme. 
Put $\Hb_*(\Bc^s_e)=\ind_y\Hb_*(\Bc^s_{e,\leq y})$.
As in 5.4 we a have the $\star$-product
$$\Hb_*({}^{w}\Zc^a_{\leq y})\times\Hb_*(\Bc^s_{e,\leq y'})\to
\Hb_*(\Bc^s_{e,\leq y''})$$
relative to $(\Tc^a)^2$ if $w\geq{y'}^{-1}$.
The Riemann-Roch map
$RR\,:\,\Kb^A({}^w\Zc_{\leq y}^a)_a\to\Hb_*({}^w\Zc^a_{\leq y})$
yields a representation of $\Kb_a$ on $\Hb_*(\Bc^s_e)$.

Let $G(s,e)\subset G^s$ be the centralizer of $e$,
and $\Pi(s,e)$ be the group of components of $G(s,e)$.
It is finite because $G^s$ is reductive.
The natural $\Pi(s,e)$-action on $\Kb_a$ is trivial, 
because the group $G^s$ acts on $\Zc^a_{ij}$ and
$G^s$ is connected by Lemma 2.13.$(i)$.
Hence the natural $\Pi(s,e)$-action on $\Hb_*(\Bc_e^s)$ 
commutes with the $\Hb_{\tau,\zeta}$-action.
For any irreducible representation $\chi$ of $\Pi(s,e)$
we have the {\sl standard module} 
$$\Nb_{a,e,\chi}=\Homb_{\Pi(s,e)}(\chi,\Hb_*(\Bc_e^s)).$$ 

Let $\Pi(s,e)^\vee$ be the set of irreducible representations of
$\Pi(s,e)$ occuring in $\Hb_*(\Bc_e^s).$ 
Let $\iota$ be the inclusion $\{e\}\subset\Nc^a$.
We have
$$\Nb_{a,e,\chi}\dot{=}\Oplus_{\chi'}\Lb_{a,\chi'}\otimes
\Homb_{\Pi(s,e)}(\chi,\Hb^*(\iota^!S_{a,\chi'}))\leqno(6.2.1)$$
because
$\Hb_*(\Bc_e^s)
\dot{=}\Oplus_i\Hb^*(\iota^!L_{a,i}).$
Hence $\Nb_{a,e,\chi}\neq\{0\}$ if and only if $\chi\in\Pi(s,e)^\vee$.

\subhead 6.3\endsubhead
In the rest of Section 6
we assume that $(\tau,\zeta)$ is regular.
By Proposition 2.14 there is a finite number of $G^s$-orbits in $\Nc^a$.
Hence the complex $S_{a,\chi}$ is the intersection cohomology complex
of an irreducible $G^s$-equivariant local system
on a $G^s$-orbit $\Oc\subset\Nc^a$.
Therefore, in the regular case, we may identify the set $\{\chi\}$ in 6.1
with a set of irreducible $G^s$-equivariant local systems on $\Nc^a$.

Given any $e\in\Oc$, a local system as above may also be viewed as
a representation of the group $\Pi(s,e)$.
Let $\Pi(a)^\vee$ be the set of all $G^s$-equivariant local systems
corresponding to a representation in $\Pi(s,e)^\vee$
for some $e\in\Nc^a$.
%We will set $\Lb_{a,\chi}=\{0\}$ if $\chi$ is a local system
%which does not lie in $\Pi(a)^\vee$.

\proclaim{\bf Proposition}
Assume that $(\tau,\zeta)$ is regular.
The set $\{\chi\}$ is finite and coincides with $\Pi(a)^\vee$.
\endproclaim

\noindent{\sl Proof:}
The corresponding statement for affine Hecke algebras is due to \cite{KL1}.
Another proof due to Grojnowski is given in \cite{CG, 8.8}.
Our proof of the proposition is very similar
to {\it loc. cit.}, to which we refer for some technical result.
Fix a $G^s$-orbit $\Oc\subset\Nc^a$ and $e\in\Oc$.
Let $\bar\Oc\subset\Nc^a$ be the Zariski closure.

\vskip2mm

\noindent{\sl Claim 1.
There is a group homomorphism $v\,:\,(\CC^\times,\times)\to(\QQ,+)$
such that $v(\zeta)>0$ and $v(\tau)<0$.}

\vskip2mm

If $\zeta=\tau^{k/m}$ with $k<0$, $m>0$, the existence follows from
\cite{CG, Lemma 8.8.12}. Else, fix $t,z\in\CC$ such that
$\exp(2i\pi t)=\tau$, $\exp(2i\pi z)=\zeta$. 
The elements $1,t,z\in\CC$ are $\QQ$-independent because
$(\tau,\zeta)$ is regular.
Hence there is a $\QQ$-linear map $\CC\to\QQ$
such that $1\mapsto 0$, $t\mapsto 1$, $z\mapsto -1$.
This map induces a group homomorphism 
$\CC^\times\to\QQ$ such that
$\tau\mapsto -1$, $\zeta\mapsto 1$.

\vskip2mm

\noindent{\sl Claim 2.
Up to conjugating $e$ by an element of $G^s$,
there is a $\slen_2$-triplet $\phi\subset\gen_{KM}$ containing $e$
such that $\CC^\times_\phi\subset H$.}

\vskip2mm

By Proposition 2.14 there are 
$S\subset\Hu\times\CC_q^\times$, 
$d\in\ZZ_{>0}$, $r\in(1/d)\ZZ_{\leq 0}$,
and $\g_d\in\Gu(K_d)$ such that $S$ is finite and the conjugation 
by $\g_d$ takes $s$ into $H$, $e$ into $\Ncu^S\otimes\eps^r$,
and $\Gu(K)^s$ onto a group containing $\Gu^S$.

Given a $\sen\len_2$-triple $\phiu=\{\eu,\fu,\hu\}\subset\genu$ 
such that $\eu\otimes\eps^r=\ad_{\g_d}e$, we set 
$$\matrix
\Mu&=\{(g,z)\in\Gu\times\CC_q^\times\,;\,\ad_g\eu=z\eu\}\hfill\cr
\Mu_\phiu&=\{(g,z)\in\Mu\,;\,\ad_gx=\ad_{\phiu(z^{1/2})}x,
\,\forall x\in G_\phiu\}\hfill
\endmatrix$$
(the choice of a square root of $z$ is irrelevant).
The subgroup $\Mu_\phiu\subset\Mu$ is maximal reductive by \cite{BV, 2.4}, 
see also \cite{CG, 8.8}.
Therefore we can choose $\phiu$ such that $S\subset\Mu_\phiu$,
because $S$ generates a reductive subgroup of $\Mu$.
Then $\CC_{\phiu}^\times\subset\Gu^S$.

Consider the $\sen\len_2$-triple
$\phiu^r=\{\eu\otimes\eps^r,\fu\otimes\eps^{-r},\hu\}\subset\genu\otimes K_d$.
Set $\phi'=\{e,f,h'\}=\ad_{\g_d}^{-1}(\phiu^r).$ 
Note that $h'\in(\genu\otimes K)^s$ because $\hu\in\genu^S$.
We have also $e\in\genu\otimes\CC[\eps^{\pm 1}]$.
Hence $\phi'\subset\genu\otimes\CC[\eps^{\pm 1}]$.
Thus there is a unique $h\in(\gen_{KM})^s$ such that
$\phi=\{e,f,h\}$  is a $\sen\len_2$-triple in $\gen_{KM}$.
We have $\CC^\times_\phi\subset G^s$ because 
$\CC_{\phi'}^\times\subset\Gu(K)^s$.
There is an element $g\in G^s$ such that $\ad_g(\CC^\times_\phi)\subset H$ 
because $H\subset G^s$ is a maximal torus.

\vskip2mm

\noindent{\sl Claim 3.
There is a non-empty subset
$\hat\Qc\subset\Bc_e^s$
such that $\ad_{G^s}(\nen^a)=\bar\Oc$
for any $\pen\in\hat\Qc$,
where $\nen\subset\pen$ is the pro-nilpotent radical.}

\vskip2mm

Fix $\zeta^{1/2}$, $s_\phi$ as in 5.1, and put
$$\gen_{t,i}=\{x\in\gen\,;\,\ad_{s_\phi}x=t\, x,\,
\ad_{\phi(z)}x=z^ix,\,\forall z\in\CC^\times\}.$$
%Hence 
%$$\gen=\bigoplus_{t\in\CC^\times}\bigoplus_{i\in\ZZ}\gen_{t,i},\quad
%\gen^s=\bigoplus_{i\in\ZZ}\gen_{\zeta^{-i/2},i},\quad 
%\gen^a=\bigoplus_{i\in\ZZ}\gen_{\zeta^{-i/2},i+2}.$$
Set also
$\qen=\bigoplus_{v(t)\leq 0}\bigoplus_{i\in\ZZ}\gen_{t,i}.$
If $k\ll 0$ then $\bigoplus_{v(t)<k}\bigoplus_{i\in\ZZ}\gen_{t,i}$
is a finite codimensional Lie subalgebra of $\uen$ which is preserved 
by the adjoint action of $\qen$ and $H$.
Hence there are $J\subsetneq I$ and $g\in G$ such that 
$\ad_g\ben_J$ is also preserved by $\qen$ and $H$
by \cite{KcW, Proposition 2.8}.
Thus $\qen\subset\ad_g\ben_J$ because
$\ad_g\ben_J$ is equal to its own normalizer in $\gen$ 
(see \cite{KcW, Lemma 1.17}).
Then $\ad_g\uen_J\subset\qen$,
and there is $g$ such that $\ad_g\ben\subset\qen$.
Hence we may choose $J$ and $g$ such that $\ad_g\ben_J=\qen$ 
by \cite{KcW, Lemma 1.5}.
Then $\ad_g\gen_J=\len$ with
$\len=\bigoplus_{v(t)=0}\bigoplus_{i\in\ZZ}\gen_{t,i}.$
Set $Q=\ad_gB_J$ and $L=\ad_gG_J$.

The set $\{\pen\in\Bc\,;\,\pen\subset\qen\}$ is non-empty and projective.
It is preserved by the action of $s,\exp(e)\in G$,
because $\ad_s\qen=\qen$ and $[e,\qen]\subset\qen$ (since $e\in\gen_{1,2}$).
Hence the simultaneous fixed point set of $s$ and $\exp(e)$ is non-empty. 
Set $\hat\Qc=\{\pen\in\Bc_e^s\,;\,\pen\subset\qen\}$.

Given $\pen\in\hat\Qc$, and $\nen\subset\pen$ its pro-nilpotent radical, we have
$$\Oc=\ad_{G^s}e\subset\ad_{G^s}(\nen^a)
\subset\ad_{G^s}(\qen^a).$$
Let $P=\ad_gB$ with $g\in G$ such that $\pen=\ad_g\ben$.
The group $G^s$ is reductive and $P^s,Q^s\subset G^s$ are parabolic.
The subset $\ad_{G^s}(\nen^a)\subset\Nc^a$ is closed because
$\ad_{P^s}(\nen^a)\subset\nen^a$ and $G^s/P^s$ is projective.
Similarly $\ad_{G^s}(\qen^a)\subset\Nc^a$ is closed because
$\ad_{Q^s}(\qen^a)=\qen^a$ and $G^s/Q^s$ is projective. 
To prove that $\ad_{G^s}(\nen^a)=\bar\Oc$
it is sufficient to prove that 
$\Oc\subset\ad_{G^s}(\qen^a)$ is dense,
or that $\ad_{Q^s}e\subset\qen^a$ is dense.
This follows from the equality $[e,\qen^s]=\qen^a$,
which is a consequence of
the theory of representations of $\slen_2$ 
and the identities
$$e\in\gen_{1,2},\quad
\qen^s=\bigoplus_{i\geq 0}\gen_{\zeta^{-i/2},i},\quad
\qen^a=\bigoplus_{i\geq 0}\gen_{\zeta^{-i/2},i+2}.$$
Claim 3 is proved.

Let $\hat\Nc^a\subset\dot\Nc^a$ be the union of the connected components 
intersecting $\hat\Qc$. 
Set also $\hat\Bc_e^{s}=\hat\Nc^a\cap\Bc_e^{s}$.
Since $\hat\Nc^a\subset\dot\Nc^a$ is open, closed, and $G^{s}$-stable,
we conclude that $\hat\Bc_e^{s}\subset\Bc_e^s$ is open, closed,
and $G(s,e)$-stable.
Therefore, $\Hb_*(\hat\Bc_e^s)$
is a $\Pi(s,e)$-submodule of $\Hb_*(\Bc_e^s)$.

\vskip2mm

\noindent{\sl Claim 4.
Any representation $\chi\in\Pi(s,e)^\vee$
occurs in $\Hb_*(\hat\Bc_e^s)$.}

\vskip2mm

Set $L(s,e)=G(s,e)\cap L$, 
and $Z^\circ_L\subset L$ equal to the connected center.
We have $e\in\gen_{1,2}\subset\len$.
Similarly $s\in L$ because 
$s\in G^{s_\phi}$, $G^{s_\phi}$ is connected by Lemma 2.13.$(i)$,
and $\gen^{s_\phi}\subset\len$.
Therefore $Z^\circ_L\subset L(s,e)$.
Thus the $G(s,e)$-action on $\Bc^{s}_e$ restricts to a $Z^\circ_L$-action.
The chain of inclusions
$\hat\Qc\subset\hat\Bc_e^{s}\subset\Bc_e^{s}$
gives a chain of inclusions of $L(s,e)$-varieties
$\hat\Qc\subset(\hat\Bc_e^{s})^{Z^\circ_L}\subset(\Bc_e^{s})^{Z^\circ_L}$.
The fixed point set $\Bc^{Z^\circ_L}$ is a disjoint union of pieces,
each $L$-equivariantly isomorphic to $\{\pen\in\Bc\,;\,\pen\subset\qen\}$.
Hence $(\Bc_e^{s})^{Z^\circ_L}$ is a disjoint union of pieces,
each $L(s,e)$-equivariantly isomorphic to $\hat\Qc$.
Hence any simple $L(s,e)$-module occuring in
$\Hb_*((\Bc_e^s)^{Z^\circ_L})$ occurs in $\Hb_*(\hat\Qc)$.
The set $\hat\Qc\subset(\hat\Bc_e^{s})^{Z^\circ_L}$ is open and closed,
because it is open and closed in $(\Bc_e^{s})^{Z^\circ_L}$.
Hence the $L(s,e)$-module $\Hb_*(\hat\Qc)$ is a direct summand in 
$\Hb_*((\hat\Bc_e^s)^{Z^\circ_L})$.
Thus any simple $L(s,e)$-module occuring in
$\Hb_*((\Bc_e^s)^{Z^\circ_L})$ occurs 
also in $\Hb_*((\hat\Bc_e^s)^{Z^\circ_L})$.
By \cite{CG, Proposition 2.5.1}
the $L(s,e)$-modules $\Hb_*((\Bc_e^s)^{Z^\circ_L})$ 
and $\Hb_*(\Bc_e^s)$ have the same class 
in the Grothendieck group.
The same holds for the $L(s,e)$-modules $\Hb_*((\hat\Bc_e^s)^{Z^\circ_L})$ 
and $\Hb_*(\hat\Bc_e^s)$.
Hence any irreducible $L(s,e)$-module that occurs in
$\Hb_*(\Bc_e^s)$ occurs also in $\Hb_*(\hat\Bc_e^s)$.
Therefore, to prove Claim 4 it is
enough to prove that the projection $G(s,e)\to\Pi(s,e)$ 
restricts to a surjective group homomorphism
$L(s,e)\to\Pi(s,e)$. 

In order to prove this, we first check that $G(s,e)\subset Q^s$.
Let $G(s_\phi,\phi)$ be the simultaneous centralizer of 
$G_\phi$ and $s_\phi$.
It is a reductive group because $G_\phi\subset G^{s_\phi}$,
and $G^{s_\phi}$, $G_\phi$ are reductive.
Set
$$\matrix
M^s&=\{(g,z)\in G^{s}\times\CC^\times_q\,;\,\ad_ge=ze\}\hfill\cr
M^s_\phi&=\{(g,z)\in M^s\,;\,\ad_gx=\ad_{\phi(z^{1/2})}x,\,
\forall x\in G_\phi\}.\hfill
\endmatrix$$
The group $M^s$ is linear. 
Let $M^s_u\subset M^s$ be the unipotent radical.
We have $M^s_u\subset G(s,e)\times\{1\}$.
Hence $\Lie(M^s_u)\subset\gen(s,e)$,
where $\gen(s,e)$ is the Lie algebra of $G(s,e)$. 
Thus $M^s_u\subset Q^s\times\{1\}$,
because $M^s_u$ is unipotent, $\gen(s,e)\subset\qen^s$
by the representation theory of $\slen_2$, 
and $Q^s=N_{G^s}(\qen^s)$.
Similarly, if $(g,z)\in M^s_\phi$ then $g$
commutes with $s_\phi$ and $\CC^\times_\phi$.
Hence $\ad_g(\qen^s)\subset\qen^s$.
Thus $M^s_\phi\subset Q^s\times\CC^\times_q$,
because $Q^s=N_{G^s}(\qen^s)$.
The group $M^s_\phi$ is reductive because it is the central extension of
$\CC^\times_q$ by $G(s_\phi,\phi)$.
We claim that $M^s_\phi\subset M^s$ is maximal reductive,
i.e. $M^s=M^s_\phi\cdot M^s_u$.
Then $G(s,e)\times\{1\}\subset M^s\subset Q^s\times\CC^\times_q$,
hence $G(s,e)\subset Q^s$.

To prove that $M^s_\phi\subset M^s$ is maximal reductive
it is sufficient to prove that $G(s_\phi,\phi)\subset G(s,e)$ is
maximal reductive, because, then, any element
$(g,z)\in M^s$ has the form $u\cdot r$ where $r$ is $(\phi(z^{1/2}),z)$ times
the component in $G(s_\phi,\phi)$ of the element
$\phi(z^{-1/2})g\in G(s,e)$, and 
$u$ belongs to the unipotent radical of $G(s,e)$.
By the representation theory of $\slen_2$, 
we have $\ad_hx\in\ZZ_{\geq 0}\cdot x$ whenever $x\in\gen(s,e).$
Let $\ven\subset\gen(s,e)$ be the Lie subalgebra linearly spanned by
$\{x\,;\,\ad_hx\in\ZZ_{>0}\cdot x\},$
and $V\subset G(s,e)$ be the corresponding unipotent normal subgroup.
The adjoint $V$-action on the affine space
$h+\ven$ is transitive, see \cite{CG, 3.7}.
The group $G(s,e)$ acts naturally on $h+\ven$,
and any reductive subgroup in $G(s,e)$ has a fixed point in $h+\ven$,
hence is conjugated to a subgroup fixing $h$, hence $\phi$.
We are done.

We can now complete the proof of Claim 4 by proving that
the projection $L(s,e)\to\Pi(s,e)$ is surjective. 
Let $g\in G(s,e)$ be a representative of
an element $\bar g\in\Pi(s,e)$.
The group $G^s$ is reductive, 
$Q^s\subset G^s$ is parabolic, and
$L^s\subset Q^s$ is a Levi subgroup.
Set $g=ru$ with $r\in L^s$ and $u$ in the unipotent radical of $Q^s$.
Then $r\in L(s,e)$ because $e\in\len$.
Moreover $r$ maps to $\bar g$ because
$\Pi(s,e)$ is a finite group.
This yields the surjectivity.

Using Claim 3 and 4 we can now prove the proposition.
Fix $\hat\Xi_s\subset\Xi_s$ such that
$\hat\Nc^a=\bigcup_{i\in\hat\Xi_s}\dot\Nc_i^a$.
We have $e\in\Oc$, and $p(\dot\Nc_i^a)=\bar\Oc$ 
for any $i\in\hat\Xi_s$ by Claim 3. 
Therefore a representation $\chi$ of $\Pi(s,e)$ occurs in
$\Hb_*(\hat\Bc_e^s)$ if and only if
the complex $IC_\chi$ occurs in $L_{a,i}$ for some $i\in\hat\Xi_s$.
Thus, for any $\chi\in\Pi(s,e)^\vee$ the complex $IC_\chi$
occurs in $L_{a,i}$ for some $i\in\Xi_s$ by Claim 4.
\qed

\subhead 6.4\endsubhead
Fix $\eu\in\Ncu^a$. 
Fix a $\slen_2$-triple $\phiu\subset\genu$ containing $\eu$.
Up to conjugating $\eu$ by an element in $\Gu^s$ we may assume that
$\CC^\times_\phiu\subset\Hu$.
Let
$$\star\,:\,\Hb_*({}^w\Zc^a_{\leq y})\times\Hb_*(\Sc^a_{\phiu,\leq y'})\to
\Hb_*(\Sc^a_{\phiu,\leq y''})$$
be as in 6.2.
It yields a representation of $\Kb_a$ on 
$\Hb_*(\Sc^a_\phiu)=\ind_y\Hb_*(\Sc^a_{\phiu,\leq y})$.
The subspace 
$\Hb_*(\Scu^a_\phiu)\subset\Hb_*(\Sc^a_\phiu)$
is preserved by $\Psi_a(\Hbu)$.

\subhead 6.5\endsubhead
Let $\Piu(\su,\eu)$ be the group of components of 
$\Gu(\su,\eu)=\Gu\cap G(\su,\eu)$.
Let $\Piu(\su,\eu)^\vee\subset\Piu(\su,\eu)$ 
be the set of irreducible representations 
occuring in $\Hb_*(\Bcu^\su_\eu).$ 
Let $\Piu(a)^\vee$ be the set
of $\Gu^\su$-equivariant local systems on $\Ncu^a$
corresponding to a representation in $\Piu(\su,\eu)^\vee$ 
for some $\eu$.

The simple representations of $\Hbu_\zeta$ 
are classified in \cite{KL1, Theorem 7.12}.
They are labelled by conjugacy classes of triples $\{s',\eu,\chiu\}$
where $s'\in\Gu\times\CC^\times_{\o_0}$ is semisimple, 
$\eu$ is as above, 
and $\chiu\in\Piu(\su,\eu)^\vee.$ 
Let $\Lbu_{a,\chiu}$ be the simple module
labelled by $\{s',\eu,\chiu\}$.

Let $\Gu(\su,\phiu)\subset\Gu$ be the 
simultaneous centralizer of $G_\phiu$ and $\su$.
The group of connected components of $\Gu(\su,\phiu)$ is $\Piu(\su,\eu)$. 
The $\Gu(\su,\phiu)$-action on $\Scu^a_\phiu$ gives a
$\Piu(\su,\eu)$-action on $\Hb_*(\Scu_\phiu^a)$ which 
commutes with the $\Hbu_\zeta$-action.
Set
$$\Mbu_{a,\phiu,\chiu}=
\Homb_{\Piu(\su,\eu)}\bigl(\chiu,\Hb_*(\Scu_\phiu^a)\bigr),\quad
{\underline n}_{a,\chiu,\chiu'}=
\dim\Homb_{\Piu(\su,\eu)}\bigl(\chiu',\Hb^*(\iota^!IC_{\chiu})\bigr),$$
where $\iota$ is the inclusion $\{\eu\}\subset\Ncu^a.$

\proclaim{Lemma}
(i)
The Grothendieck group of finite dimensional $\Hbu_\zeta$-modules
is spanned by $\{\Mbu_{a,\phiu,\chiu}\}.$

(ii)
$[\Mbu_{a,\phiu,\chiu}:\Lbu_{a,\chiu'}]={\underline n}_{a,\chiu,\chiu'}$.
\endproclaim

\noindent{\sl Proof:}
This is standard, see \cite{CG} for instance. 
For $(ii)$ use the identity
$${\underline n}_{a,\chiu,\chiu'}=\dim\Homb_{\Piu(\su,\eu)}
\bigl(\chiu',\Hb^*(\senu^a_{\phiu},\eps^!IC_\chiu)\bigr),\leqno(6.5.1)$$
where $\eps$ is the natural inclusion $\senu^a_\phiu\subset\Ncu^a$.
Since $IC_\chiu$ is $C_\phiu$-equivariant (6.5.1) follows from 
the basic result below. 

Let $V$ be a finite dimensional vector space with a 
$\CC^\times$-action contracting to 0. 
Let $Z\subset V$ be a closed subset preserved by $\CC^\times$,
$\iota\,:\,\{0\}\to Z$ and $p\,:\,Z\to\{0\}$ 
be the obvious inclusion and projection.
Then the canonical map $\iota^!\Ec\to p_!\Ec$ is an isomorphism for
any $\Ec\in\Db_{\CC^\times}(Z)$.
See \cite{KS, Proposition 3.7.5} for instance.
\qed

\subhead 6.6\endsubhead
The $\Piu(\su,\eu)$-action on $\Hb_*(\Sc_\phiu^a)$
induced by the $\Gu(\su,\phiu)$-action on $\Sc_\phiu^a$
commutes with $\Hb_{\tau,\zeta}$ because $\Gu^\su$ is connected.
Set
$\Mb_{a,\phiu,\chiu}=
\Homb_{\Piu(\su,\eu)}\bigl(\chiu,\Hb_*(\Sc^a_\phiu)\bigr)$
for each $\chiu\in\Piu(\su,\eu)^\vee$.

\proclaim{Lemma}
$\Hb_{\tau,\zeta}\otimes_{\Hbu_\zeta}\Mbu_{a,\phiu,\chiu}\simeq
\Mb_{a,\phiu,\chiu}$.
\endproclaim

\noindent{\sl Proof:}
We have $a\in A_e$.
Taking $D\subset A_e$ to be the closed subgroup generated by 
and $(A_e)^\circ\cup\{a\}$, 
Proposition 5.8 yields an isomorphism of $\Hb_{\tau,\zeta}$-modules
$$\Hb_{\tau,\zeta}\otimes_{\Hbu_\zeta}\Hb_*(\Scu^a_\phiu)\simeq
\Hb_*(\Sc_\phiu^a).$$
This isomorphism is $\Piu(\su,\eu)$-equivariant, because
the vector bundle $\g$ is $\Gu(\su,\phiu)$-equivariant. 
\qed

\vskip3mm

Using Lemma 6.5 and Lemma 6.6, 
it is enough to compute
the Jordan-H\"older multiplicities of $\Mb_{a,\phiu,\chiu}$ 
to get the Jordan-H\"older multiplicities of the induced module
$\Hb_{\tau,\zeta}\otimes_{\Hbu_\zeta}\Lb$ for 
any simple $\Hbu_\zeta$-module $\Lb$.

\subhead 6.7\endsubhead
For any $\chi\in\Pi(a)^\vee$, $\chiu\in\Piu(a)^\vee$ we set
$n_{a,\chiu,\chi}=
\dim\Homb_{\Piu(\su,\eu)}\bigl(\chiu,\Hb^*(\sen_\phiu^a,\eps^!IC_{\chi})\bigr),$
where $\eps$ is the natural inclusion $\sen^a_\phiu\subset\Nc^a$.

\proclaim{Lemma}
Assume that $(\tau,\zeta)$ is regular.
Then $\Mb_{a,\phiu,\chiu}$ has a finite length, and
$[\Mb_{a,\phiu,\chiu}:\Lb_{a,\chi}]=n_{a,\chiu,\chi}$.
\endproclaim

\noindent{\sl Proof:}
Fix $e\in\Nc^a$ such that $\chi\in\Pi(s,e)^\vee$.
The Springer resolution restricts to a map
$\Sc_\phiu^a\to\sen_\phiu^a$,
yielding a Cartesian square
$$\matrix
\Sc_\phiu^a&\lra&\dot\Nc^a&&\cr
\bda&&\bda{\ss p}&&\cr
\sen_\phiu^a&{\buildrel\eps\over\lra}&\Nc^a
&{\buildrel\iota\over\lla}&\{\eu\}.
\endmatrix$$
Thus
$\Hb_*(\Sc^a_\phiu)\dot{=}
\bigoplus_i\Hb^*(\sen_\phiu^a,\eps^!L_{a,i}),$
yielding the lemma.
\qed

\vskip3mm

Note that 
$n_{a,\chiu,\chi}=
\dim\Homb_{\Piu(\su,\eu)}\bigl(\chiu,\Hb^*(\iota^!IC_{\chi})\bigr)$
because $IC_\chi$ is $C_\phiu$-equivariant.
In particular $n_{a,\chiu,\chi}$
does not depend on the choice of $\phiu$.

\subhead 6.8\endsubhead
This subsection is given for the sake of completness.
It is not used in the rest of the paper.

\proclaim{Lemma}
If $e\in\Nc^a$ is nilpotent  
then $\Bc^s_e$ has no odd rational homology.
\endproclaim

\noindent{\sl Proof:}
The proof is identical to the proof of \cite{KL1, Theorem 4.1}.
By Claim 2 in 6.3 there is a $\slen_2$-triple $\phi\subset\gen_{KM}$ 
containing $e$ with $C_\phi\subset A$.
Then the element $s_\phi$ in Claim 3 belongs to $H'\times\{\tau\}$ 
and $e\in\gen^{s_\phi}$. 
In particular, $\Bc^{s_\phi}_e$ is locally of finite type
and each connected component is isomorphic to the variety of Borel
subgroups of $G^{s_\phi}$ containing $e$.
Hence $\Bc^{s_\phi}_e$ has no odd rational homology.
Therefore, to conclude it is sufficient to observe that
$(\Bc^{s_\phi}_e)^{C_\phi}=(\Bc^s_e)^{C_\phi}$, that
the $C_\phi$-action on $\Sc^a_\phi$ is a contraction to $\Bc_e^s$,
and that $\Sc^a_\phi$ is locally of finite type and smooth.
\qed

\head 7. Induction of sheaves and Fourier transform\endhead
\subhead 7.1\endsubhead
Assume that $\tau$, $\zeta$ are not roots of unity.
Fix $J\subsetneq I$ and $\qen\in(\Bc^J)^{s}$. 
Let $\nen$ be the pro-nilpotent radical of $\qen$, 
and $\len\subset\qen$ a Levi Lie subalgebra. 
Let $Q,N,L$ be the corresponding subgroups of $G$.
We have the obvious projection 
$\qen\to\len$, $x\mapsto x+\nen$,
and the natural inclusion $\qen\subset\gen$.
Taking the group $L$ instead of $\Gu$ in 2.9 we get the varieties
$\Bc_L$, $\Nc_L$, $\dot\Nc_L$ and $\dot\len$.
We may identify $\Bc_L$ with the subset 
$\{\ben\in\Bc\,;\,\ben\subset\qen\}.$
It yields an inclusion $\dot\len\subset\dot\gen.$

Set $\gen_\pm=\{x\in\gen\,;\,\ad_sx=\zeta^{\mp 1}x\}$. 
Idem for $\qen_\pm$, $\len_\pm$, $\nen_\pm$.
The groups $Q^{s}$, $N^{s}$ act on 
$\qen_\pm$ in the obvious way.
Put $\Ec_\pm=G^{s}\times_{N^{s}}\qen_\pm$
and $\Fc_\pm=G^{s}\times_{Q^{s}}\qen_\pm$.
We have the map 
$p_1\,:\,\Ec_\pm\to\len_\pm$, $[g:x]\mapsto q_1(x)$, 
the map $p_2\,:\,\Ec_\pm\to\Fc_\pm$,  $[g:x]\mapsto [g:x]$, 
and the map $p_3\,:\,\Fc_\pm\to\gen_\pm$, $[g:x]\mapsto\ad_gx$,
yielding the diagram
$$\len_\pm
{\buildrel p_1\over\lla}
\Ec_\pm
{\buildrel p_2\over\lra}
\Fc_\pm
{\buildrel p_3\over\lra}
\gen_\pm.$$
The induction functor
$\Indb^{\gen_\pm}_{\qen_\pm}\,:\,
\Db_{L^{s}}(\len_\pm)\to\Db_{G^{s}}(\gen_\pm)$
takes $E$ to $p_{3!}(p_2^*)^{-1}p_1^*(E).$

Fix connected components 
$\dot\gen_\pm\subset\{x\in\dot\gen\,;\,\ad_{(s,\zeta^{\pm 1})}x=x\}$, 
and 
$\dot\len_\pm\subset\{x\in\dot\len\,;\,\ad_{(s,\zeta^{\pm 1})}x=x\}$.
Set $\dot\Nc_\pm=\dot\Nc\cap\dot\gen_\pm$ and
$\dot\Nc_{L,\pm}=\dot\Nc_L\cap\dot\len_\pm$.
Recall the projection $p\,:\,\dot\gen\to\gen$ from 2.8.
Then $\dot\gen_\pm$ is a connected component in $p^{-1}(\gen_\pm).$ 
Let $p$ denote also the restriction of $p$ to $\dot\gen_\pm$, $\dot\Nc_\pm$.
Idem for $\dot\len_\pm$, $\dot\Nc_{L,\pm}$.

\proclaim{7.2. Lemma} 
If $\dot\len_\pm\subset\dot\gen_\pm$ there are canonical isomorphisms
$\Indb_{\qen_\pm}^{\gen_\pm}(p_!\CC_{\dot\len_\pm})=p_!\CC_{\dot\gen_\pm}$,
and
$\Indb_{\qen_\pm}^{\gen_\pm}(p_!\CC_{\dot\Nc_{L,\pm}})=p_!\CC_{\dot\Nc_\pm}$.
\endproclaim

\noindent{\sl Proof:}
Fix an Iwahori Lie algebra $\pen\subset\qen$ such that
$(0,\pen/\nen)\in\dot\len_\pm$ and $(0,\pen)\in\dot\gen_\pm$.
Set $P\subset G$ equal to the Iwahori subgroup associated to $\pen$,
and $\pen_\pm=\pen\cap\gen_\pm$.
We have
$\dot\gen_\pm=G^{s}\times_{P^{s}}\pen_\pm$ 
and
$\dot\len_\pm=Q^{s}\times_{P^{s}}(\pen_\pm/\nen_\pm).$
Set $\dot\qen_\pm=Q^{s}\times_{P^{s}}\pen_\pm$,
$\dot\Ec_\pm=G^{s}\times_{N^{s}}\dot\qen_\pm$,
and
$\dot\Fc_\pm=G^{s}\times_{Q^{s}}\dot\qen_\pm$.
Note that $\dot\Fc_\pm\simeq\dot\gen_\pm$.
We get the commutative diagram
$$\matrix
\dot\len_\pm\hfill&
{\buildrel\dot p_1\over\lla}&
\dot\Ec_\pm&
{\buildrel\over\lra}&
\dot\Fc_\pm&=&\dot\gen_\pm\cr
{\ss p}\bda&&{\ss p}\bda&&{\ss p}\bda&&{\ss p}\bda\cr
\len_\pm&
{\buildrel p_1\over\lla}&
\Ec_\pm&
{\buildrel p_2\over\lra}&
\Fc_\pm&
{\buildrel p_3\over\lra}&
\gen_\pm,
\endmatrix
$$
where $\dot p_1([g_1:g_2:x])=[g_2:x+\nen_\pm]$ 
for all $g_1\in G^s$, $g_2\in Q^s$.
Both left squares are Cartesian.
By base change we get 
$p_2^*p_!(\CC_{\dot\Fc_\pm})=p_1^*p_!(\CC_{\dot\len_\pm}).$
Hence
$$\Indb^{\gen_\pm}_{\qen_\pm}(p_!\CC_{\dot\len_\pm})=
(p_3p)_!(\CC_{\dot\Fc_\pm})=p_!\CC_{\dot\gen_\pm}.\leqno(7.2.1)$$
Put
$\Ec^{nil}_\pm=G^{s}\times_{N^{s}}(\qen_\pm\cap\Nc)$ and
$\Fc^{nil}_\pm=G^{s}\times_{Q^{s}}(\qen_\pm\cap\Nc)$.
The diagram
$$\matrix
\len_\pm&
{\buildrel p_1\over\lla}&
\Ec_\pm&
{\buildrel p_2\over\lra}&
\Fc_\pm\hfill&
{\buildrel p_3\over\lra}&
\gen_\pm\hfill\cr
\bua&&\bua&&\bua&&\bua\cr
\len_\pm\cap\Nc_L&
{\buildrel\over\lla}&
\Ec^{nil}_\pm&
{\buildrel\over\lra}&
\Fc^{nil}_\pm&
{\buildrel\over\lra}&
\gen_\pm\cap\Nc\hfill\cr
\endmatrix
$$
is formed by Cartesian squares because
$q_1^{-1}(\len_\pm\cap\Nc_L)=\qen_\pm\cap\Nc$.
Hence the lemma follows from (7.2.1) and base change.
\qed

\subhead 7.3\endsubhead
Fix a $G$-invariant nondegenerate pairing $\la\,:\,\ra$ on $\gen$.
It restricts to a non-degenerate pairing
$\gen_+\times\gen_-\to\CC$.
The complexes $p_!\CC_{\dot\Nc_\pm}$, $p_!\CC_{\dot\gen_\pm}$ are
$\CC^\times_{\gen_\pm}$-equivariant.
Let $\rho$ be as in 2.5.

\proclaim{\bf Lemma}
(i)
If $E\in\Db_{L^s}(\len_\pm)$ is $\CC^\times_{\len_\pm}$-equivariant then 
$\Indb_{\qen_\pm}^{\gen_\pm}(E)$ is $\CC^\times_{\gen_\pm}$-equivariant, and
$\Fb_{\gen_\pm}\circ\Indb_{\qen_\pm}^{\gen_\pm}(E)\dot{=}
\Indb_{\qen_{\mp}}^{\gen_{\mp}}\circ\Fb_{\len_\pm}(E)$.

(ii) 
If $\rho(\dot\gen_\pm)=\rho(\dot\gen_\mp)$
then $\Fb_{\gen_\pm}(p_!\CC_{\dot\Nc_\pm})\dot{=}p_!\CC_{\dot\gen_\mp}$
and
$\Fb_{\gen_\mp}(p_!\CC_{\dot\gen_\mp})\dot{=}p_!\CC_{\dot\Nc_\pm}$.
\endproclaim

\noindent{\sl Proof:}
Claim $(i)$ is analogous to \cite{L1, Corollary 10.5}.
The first part of $(ii)$ follows from $(i)$ and Lemma 7.2 
with $\qen\in\rho(\dot\gen_\pm)$, because
$$\Fb_{\gen_\pm}(p_!\CC_{\dot\Nc_\pm})=
\Fb_{\gen_\pm}\circ\Indb_{\qen_\pm}^{\gen_\pm}(\CC_{\{0\}})\dot{=}
\Indb_{\qen_\mp}^{\gen_\mp}\circ\Fb_{\len_\pm}(\CC_{\{0\}})=
\Indb_{\qen_\mp}^{\gen_\mp}(\CC_{\len_\mp})=
p_!\CC_{\dot\gen_\mp}.$$
The second part follows from the first one because the complexes 
$p_!\CC_{\dot\gen_\pm}$, $p_!\CC_{\dot\Nc_\pm}$ 
are $\CC^\times_{\gen_\pm}$-equivariant.
\qed

\subhead 7.4\endsubhead
For any pair $(\tau,\zeta)$, possibly non regular,
we set $\Pi(a)^\vee=\{\chi\}$, see 6.1.
Hence $\Pi(a)^\vee$ may be viewed as a set consisting of
simple $G^s$-equivariant perverse sheaves on $\Nc^a$.
This notation is coherent with Proposition 6.3 if $(\tau,\zeta)$ is regular.

Assume that $(\tau^{-1},\zeta)$ is singular.
Hence $(\tau,\zeta)$ is regular.
Recall that $a=(s,\zeta)$, see 2.12.
Set $\check a=(s^{-1},\zeta)$.
Fix $\chi\in\Pi(a)^\vee$.
Note that $\Psi_a^\bullet\Phi_a^\bullet(\Lb_{a,\chi})$ 
is a simple integrable $\Hb_{\tau,\zeta}$-module by 
Lemma 6.1.$(ii)$ and Claim 1 in 7.6 below.
Hence $\IM^\bullet\Psi_a^\bullet\Phi_a^\bullet(\Lb_{a,\chi})$ is again 
is a simple integrable $\Hb_{\tau^{-1},\zeta}$-module.
Therefore, by Lemma 6.1.$(ii)$ and Claim 2 in 7.6
there is a unique element $\check\chi\in\Pi(\check a)^\vee$ such that
the $\Hb_{\tau^{-1},\zeta}$-modules
$\Psi_{\check a}^\bullet\Phi_{\check a}^\bullet(\Lb_{\check a,\check\chi})$ 
and  
$\IM^\bullet\Psi_a^\bullet\Phi_a^\bullet(\Lb_{a,\chi})$ are isomorphic. 
The map $\Pi(a)^\vee\to\Pi(\check a)^\vee$, $\chi\mapsto\check\chi$
is a bijection.
In particular $\Pi(\check a)^\vee$ is a finite set by Proposition 6.3.

We claim that the Fourier transform yields another bijection 
$\Pi(a)^\vee\to\Pi(\check a)^\vee$. 
More precisely, set $\bar a=(s,\zeta^{-1})$.
We have $\Xi_s=\Xi_{s^{-1}}$ and $L_{\bar a,i}=L_{\check a,i}$,
hence $\Pi(\bar a)^\vee=\Pi(\check a)^\vee$, because the fixed point sets
$\gen^{\bar a}$, $\gen^{\check a}$ coincide.
We have also
$L_{\bar a,i}\dot{=}\Fb_{\gen^a}(L_{a,i})$
by Lemma 7.3.$(ii)$, and $S_{a,\chi}=IC_\chi$ for all $\chi\in\Pi(a)^\vee$. 
The complex $IC_\chi$ is $\CC^\times_{\gen^a}$-equivariant, and
$\Fb_{\gen^a}(IC_\chi)$ is a simple perverse sheaf.
Hence we get a bijection $\Pi(a)^\vee\to\Pi(\bar a)^\vee$, 
$\chi\mapsto\bar\chi$, such that
$S_{\bar a,\bar\chi}=\Fb_{\gen^a}(IC_\chi)$.
The composed map 
$\Pi(a)^\vee\to\Pi(\bar a)^\vee=\Pi(\check a)^\vee$  
is still denoted by $\chi\mapsto\bar\chi$.

\proclaim{Conjecture}
$\bar\chi=\check\chi$ for each $\chi\in\Pi(a)^\vee$ .
\endproclaim

\noindent
The corresponding statement for $\Hbu$ is proved in \cite{EM}.
Probably, this proof generalizes to the case of $\Hb$.

\subhead 7.5\endsubhead
For any $\chiu\in\Piu(a)^\vee$, $\chi\in\Pi(a)^\vee$ let 
$m_{a,\chiu,\chi}$ be the Jordan-H\"older multiplicity of $\Lb_{a,\chi}$ in 
$\Hb_{\tau,\zeta}\otimes_{\Hbu_\zeta}\Lbu_{a,\chiu}$.
If $(\tau,\zeta)$ is regular then $m_{a,\chiu,\chi}$ 
is known by 6.5-7.
If $(\tau,\zeta)$ is singular then 
$m_{a,\chiu,\chi}=m_{\check a,\check\chiu,\check\chi}$,
thus it is known because $(\tau^{-1},\zeta)$ is regular.

\subhead 7.6\endsubhead
An $\Hb_{\tau,\zeta}$-module is said to be {\sl integrable}
if the action of the subring $\Rb_{\tau,\zeta}$ is locally finite.
Let $\Modb^{int}_{\Hb_{\tau,\zeta}}\subset\Modb_{\Hb_{\tau,\zeta}}$
be the full subcategory consisting of integrable modules. 
The following theorem proves a refined version of \cite{BEG1, Conjecture 6.5}.
Note that this conjecture was known, 
at least to the author, before {\it loc. cit.} 

\proclaim{Theorem}
Assume that $\tau,\zeta$ are not roots of unity.

(i) If $(\tau,\zeta)$ is regular,
the simple integrable $\Hb_{\tau,\zeta}$-modules are classified 
by the $G$-conjugacy classes of triples $\{s,e,\chi\}$ with 
$s\in G'\times\{\tau\}$ semisimple, $e\in\Nc$, $\chi\in\Pi(s,e)^\vee$,
and $\ad_se=\zeta e$.

(ii) If $(\tau,\zeta)$ is not regular,
the simple integrable $\Hb_{\tau,\zeta}$-modules are 
in one-to-one correspondence with the $\Hb_{\tau^{-1},\zeta}$-modules.
\endproclaim

\noindent{\sl Proof:}
Note that part $(ii)$ follows from part $(i)$ by 7.4.
The proof consists in three steps.

\vskip3mm

\noindent{\sl Claim 1. 
The pull-back by $\Psi_a$ takes smooth and simple $\Kb_a$-modules
to simple integrable $\Hb_{\tau,\zeta}$-modules supported on
$W(s)\times\{\zeta\}\subset A$. 
}

\vskip3mm
The restriction of the map
$\Psi_a\,:\,\Hb\to\Kb_a$ to $\Rb_A$ factorizes through the obvious map
$\Rb_A\to\prod_{b\in W(a)}\Rb_A/\Ib_b$ by Lemma 4.8.$(ii)$.
Therefore, given a smooth $\Kb_a$-module $\Mb$, 
the $\Rb_A$-action on $\Psi_a^\bullet(\Mb)$ factorizes through
$\bigoplus_{b\in W(a)}\Rb_A/\Ib_b$.
Hence it is locally finite 
because $\Rb_A/\Ib_b$ is a finite dimensional $\CC$-algebra,
and the restriction of the 
$\Hb_{\tau,\zeta}$-action to $\Rb_{\tau,\zeta}$ 
takes $\Psi_a^\bullet(\Mb)$ to a quasi-coherent sheaf on 
$H'\times\{\tau\}$ supported on $W(s)$. 
Thus $\Psi_a^\bullet(\Mb)$ is integrable and 
the claim follows from Theorem 4.9.$(iv)$.

\vskip3mm

\noindent{\sl Claim 2. 
For any simple integrable $\Hb_{\tau,\zeta}$-module $\Lb$
there is an element $s\in H'\times\{\tau\}$
and a smooth simple $\Kb_a$-module $\Lb'$ such that $\Lb=\Psi_a^\bullet(\Lb')$.}

\vskip3mm

Taking $e=0$ in Theorem 5.8 we get
$\Kb^{A}(\Sc_\phi^a)_a=
\Hb_{\tau,\zeta}\otimes_{\Hbu_\zeta}\Kb^{A}(\dot\Ncu)_a$,
with $\Sc_\phi=\dot\Nc\cap(\ben_\flat\times\Bc)$.
By \cite{L3, Proposition 8.5} we have an isomorphism of $\Hbu_\zeta$-modules
$\Kb^{A}(\dot\Ncu)_a=\Hbu_\zeta\otimes_{\Rb_A}\CC_a$.
Hence 
$\Hb_{\tau,\zeta}\otimes_{\Rb_A}\CC_a=\Psi_a^\bullet(\Kb^{A}(\Sc_\phi^a)_a)$
where $\Kb_a$ acts on $\Kb^A(\Sc_\phi)_a$ as in 5.4.
If $\Lb$ is an integrable $\Hb_{\tau,\zeta}$-module there is
a non-zero element $x\in\Lb$ which is an $\Rb_A$-eigenvector.
Fix $s\in H'\times\{\tau\}$ such that
the $\Rb_A$-module spanned by $x$ 
is isomorphic to $\CC_a$.
If $\Lb$ is simple the linear map
$\rho_x\,:\,\Hb_{\tau,\zeta}\otimes_{\Rb_A}\CC_a\to\Lb$ 
such that $\rho_x(y\otimes 1)=y(x)$ is surjective.
The subspace $\Ker(\rho_x)\subset\Hb_{\tau,\zeta}\otimes_{\Rb_A}\CC_a$ 
is preserved by $\Hb_{\tau,\zeta}$, hence by $\Kb_a$
because $\Psi_a(\Hb)\subset\Kb_a$ is dense and
$\Kb^A(\Sc_\phi)_a$ is a smooth $\Kb_a$-module. 
We are done.

\vskip3mm

\noindent{\sl Claim 3. 
$\Psi_{a_1}^\bullet(\Lb_{a_1,\chi_1})\simeq
\Psi_{a_2}^\bullet(\Lb_{a_2,\chi_2})$ if and only if
the corresponding triples are $G$-conjugated.}

\vskip3mm

By Claim 2 a simple integrable $\Hb_{\tau,\zeta}$-module $\Lb$
is isomorphic to $\Psi_a^\bullet(\Lb')$ for an element 
$s\in H'\times\{\tau\}$ and a simple smooth $\Kb_a$-module $\Lb'$.
The support of $\Lb$, viewed as a $\Rb_A$-module,
is contained in a $W$-orbit because $\Lb$ is simple.
By Claim 1 this orbit is precisely $W(s)\times\{\zeta\}.$ 
Hence $a$ is uniquely determined, up to conjugation by $W$.
Given $a$, there is a unique module $\Lb'$  such that $\Lb=\Psi_a^\bullet(\Lb')$
by Theorem 4.9.$(iv)$.
Assume that $(\tau,\zeta)$ is regular.
Then, by Lemma 6.1 and Proposition 6.3, the set of simple $\Kb_a$-modules 
is $\{\Lb_{a,\chi}\,;\,\chi\in\Pi(a)^\vee\}$.
Hence the theorem follows from Proposition 4.4.$(ii)$.
\qed

\vskip3mm

Given a simple integrable $\Hb_{\tau,\zeta}$-module $\Lb$,
the corresponding $G$-conjugacy class of triples $\{s,e,\chi\}$,
or, abusively, any element in this conjugacy class, are 
called the {\sl Langlands parameters} of $\Lb$.

\head 8. The type A case\endhead

In Sections 8.1 to 8.3
we compare our classification of the simple modules with the 
classification announced by Cherednik, in \cite{C2}, in type $A$.
In Sections 8.4 to 8.5 we make a link with representations of the cyclic quiver
and with Kazhdan-Lusztig polynomials.

\proclaim{8.1. Lemma}
Assume that $\Gu=\SL_d$ and $(\tau,\zeta)$ is regular.
Then $\Pi(s,e)^\vee$ is reduced to the trival representation of
$\Pi(s,e)$ for all $s,e$.
\endproclaim

\noindent{\sl Proof:}
Up to conjugation by an element of $G$
there is a maximal proper subset $J\subset I$ such that
$s\in G_J$ and $G^s=G_J^s$, see the proof of Lemma 2.13. 
Let $Z_J,G'_J\subset G_J$ be the center and the derived subgroup of $G_J$.
We write $Z_J^\circ$ for $(Z_J)^\circ$.
We have $G_J=Z^\circ_J\cdot G'_J$, $G'_J\simeq\SL_d$,
and $G^s_J=Z^\circ_J\cdot(G'_J)^s.$
The restriction of the $G_J^s$-action on $\dot\Nc^a$ to $Z_J$ is trivial.
Hence the $(G'_J)^s$-action on $\dot\Nc^a$
extends to an action of $(\GL_d)^s$.
Thus the $G_J^s$-equivariant local systems on $\Nc^a$
which belong to $\Pi(a)^\vee$ are $(\GL_d)^s$-equivariant.
Hence they are trivial. 
\qed

\vskip3mm

\noindent{\bf Remark.}
The group $\Pi(s,e)$ may be non trivial, 
see Remark 9.3.$(ii)$ for instance.

\subhead 8.2\endsubhead
A multi-segment is a familly
$\sigma=(\sigma_{ij})$ with $\sigma_{ij}\in\ZZ_{\geq 0}$,
$i,j\in\ZZ$, $i\leq j$.
Given $\sigma$ and $n\in\ZZ$ let $\sigma[n]$ be the multi-segment
such that $(\sigma[n])_{ij}=\sigma_{n+i,n+j}$.
Fix $d,m\in\ZZ_{>0}$.
The multi-segment $\sigma$ is $m$-periodic of order $d$ if
$\sigma[m]=\sigma$ and $\sum_i\sum_{1\leq j\leq m}(j+1-i)\sigma_{ij}=d.$

A $m$-periodic pair of order $d$ is a pair $\{(\sigma_a), (z_a)\}$,
where $(\sigma_a)$ is a familly of $m$-periodic multisegments 
such that the sum of the orders of each of them is $d$, 
and $(z_a)$ is a familly of complex numbers
such that $z_a/z_b\notin\tau^{(1/m)\ZZ}$ for all $a\neq b$.
Let $\Sc$ be the set of $m$-periodic pairs of order $d$. 
Two pairs are equivalent
if they can be obtained from each other by a simultaneous permutation 
of $(\sigma_a)$ and $(z_a)$, 
combined with a translation 
$\{(\sigma_a),(z_a)\}\to\{(\sigma_a),(z_a\tau^{m_a})\}$,
and with a simultaneous translation 
$\{(\sigma_a),(z_a)\}\to\{(\sigma_a[-m_a]),(z_a\zeta^{m_a})\}$ 
with $m_1,...,m_r\in\ZZ$.
Let $\Sc/\!\!\sim$ 
be the set of equivalence classes in $\Sc$. 

From now on we fix $k\in\ZZ_{<0}$ with $(m,k)=1$, 
a $m$-th root of $\tau$, and $\zeta=\tau^{k/m}$.
Hence $\tau^\ZZ\zeta^\ZZ=\tau^{(1/m)\ZZ}.$
Set $\Gu=\GL_d$. 
Let $\Cc$ be the set of pairs $(s,e)$,
with $s=(s',\tau)\in G_1$ semisimple and $e\in\genu\otimes K$ 
nilpotent such that $\ad_se=\zeta e$.
Let $\Cc/\!\!\sim$ be the set of equivalence classes in $\Cc$,
where $(s_1,e_1)\sim(s_2,e_2)$ if they are $G_1$-conjugated.

\proclaim{Proposition}
There is a natural bijection $\Cc/\!\!\sim\to\Sc/\!\!\sim$.
\endproclaim

\noindent{\sl Proof:}
The group $G_1=L\Gu\rtimes\CC^\times_\delta$ acts naturally
on $\CC^d[\eps,\eps^{-1}]$, so that
$(g,z)\cdot(v\otimes\eps^j)=z^jg(v)\otimes\eps^j$.
For any semisimple element $s\in G_1$ 
we fix a familly of complex numbers $({}^sz_a)$ such that
$\Spec(s)=\{{}^sz_a\}\,\mod\,\tau^{(1/m)\ZZ}$ 
and ${}^sz_a/{}^sz_b\notin\tau^{(1/m)\ZZ}$ for all $a\neq b$,
and we set 
${}^sV_{ai}=\{v\in\CC^d[\eps,\eps^{-1}]\,;\,s(v)={}^sz_a\tau^{i/m}v\}$,
${}^sV_{a}=\bigoplus_i{}^sV_{ai}$.
If $s\in\Gu(K)\times\{\tau\}$ then 
${}^sV_{a}$ is a $\ZZ$-graded $\CC[\eps,\eps^{-1}]$-module
such that ${}^sV_{ai}$ is the degree $i$ subspace.
Set ${}^sd_a=\sum_i(\dim{}^sV_{ai})t^i$.
Let $\Dc$ be the set of pairs $\{(d_a),(z_a)\}$
where $(z_a)$ is as above, $d_a\in\ZZ_{\geq 0}[[t, t^{-1}]]$,
$\sum_ad_a=p(t)\sum_{i\in\ZZ}t^{im}$ with $p(t)\in\ZZ_{\geq 0}[t,t^{-1}]$
such that $p(1)=d$.
Two pairs in $\Dc$ are equivalent if they can be obtained from each other
by a simultaneous permutation of $(d_a)$ and $(z_a)$, 
combined with a translation 
$\{(d_a),(z_a)\}\mapsto\{(d_at^{-mm_a}),(z_a\tau^{m_a})\}$
or
$\{(d_a),(z_a)\}\mapsto\{(d_at^{-km_a}),(z_a\zeta^{m_a})\}$
with $m_1,...,m_r\in\ZZ$.
The map $s\mapsto\{({}^sd_a), ({}^sz_a)\}$ 
gives a bijection from the set of $G_1$-conjugacy 
classes of semisimple elements in $\Gu(K)\times\{\tau\}$ 
to the set of equivalence classes in $\Dc$,
because $(m,k)=1$.

For any pair $\{s,e\}\in\Cc$ 
the element $e$ acts naturally on $\CC^d[\eps,\eps^{-1}]$.
Let $e_a\in\End({}^sV_a)$ be the restriction of $e$.
The operator $e_a$ is a nilpotent endomorphism 
of degree $k$ commuting with $\eps$.
Let ${}^e\sigma_{aij}$ be the number of Jordan blocks of $e_a$
of multi-degree $\{ki,ki+k,...,kj\}$.
Hence ${}^e\sigma_a=({}^e\sigma_{aij})$ is a 
$m$-periodic multi-segment of order the rank of the
$\CC[\eps,\eps^{-1}]$-module ${}^sV_a$, such that
${}^sd_a=\sum_{i,j}{}^e\sigma_{aij}(t^{ki}+t^{k(i+1)}+\cdots+t^{kj})$.

The map $\Cc\to\Sc/\!\!\sim$ taking
the conjugacy class of $\{s,e\}$ to the class of
$\{({}^e\sigma_a),({}^sz_a)\}$ gives a map 
$\Cc/\!\!\sim\to\Sc/\!\!\sim$.
It is a bijection.
\qed

\subhead 8.3\endsubhead
In this subsection we compare the simple modules of different versions
of the double affine Hecke algebra when $\Gu$ is of type $A_{d-1}$, 
with $d\geq 2$ : 

$(i)$ $\Hb^\ch_{\tau,\zeta}$ is the $\CC$-algebra generated by 
$x_1^{\pm 1}$,...,$x_d^{\pm 1}$,$\pi^{\pm 1}$,
$t_1$,...,$t_{d-1}$, modulo the relations
$$\matrix
x_ix_j=x_jx_i,\ \pi x_i=x_{i+1}\pi,\ \pi^d x_i=\tau^{-1}x_i\pi^d,\hfill&\cr
(t_i-\zeta)(t_i+1)=0,\ \pi t_i=t_{i+1}\pi,\ \pi^d t_i=t_i\pi^d,\hfill&\cr
t_it_jt_i\cdots=t_jt_it_j\cdots\ \text{if}\ i\neq j\ 
\text{($m_{ij}$\ factors\ in\ both\ products)},\hfill\cr
t_ix_it_i=\zeta x_{i+1},\ t_ix_j=x_jt_i
\ \text{if}\ j\neq i,i+1.\hfill\cr
\endmatrix$$

$(ii)$ $\Hb'_{\tau,\zeta}\subset\Hb^\ch_{\tau,\zeta}$ 
is the subalgebra generated by 
$x_1^{\pm 1}$,...,$x_d^{\pm 1}$,
$t_1$,...,$t_{d-1}$, and the element $t_0=\pi^{-1}t_1\pi$.

$(iii)$ $\Hb_{\tau,\zeta}$ is obtained by taking $\Gu=\SL_d$ in 3.1.
In particular the subring spanned by $\{x_\l\,;\,\l\in X'\}$ is isomorphic to
$\Rb_{H'}$.

Since the algebras $\Hb_{\tau,\zeta}$ and $\Hb'_{\tau,\zeta}$ 
differ only in the polynomial subalgebras $\Rb_{H'}$ and 
$\CC[x_1^{\pm 1},...x_d^{\pm 1}]$,
a simple integrable $\Hb'_{\tau,\zeta}$-module 
may be viewed as a pair $\{\Lb,\a\}$
where $\Lb$ is a simple integrable $\Hb_{\tau,\zeta}$-module, and
$\a$ is a locally finite representation of
$\CC[x_1^{\pm 1},...x_d^{\pm 1}]$ in $\Lb$ 
compatible with the $\Rb_{H'}$-action in the obvious way.
Two pairs $\{\Lb,\a\}$, $\{\Lb,\a'\}$ 
can only differ by a scalar $z\in\CC^\times$ 
such that $\a'(x_1)=z\a(x_1)$, 
because $\Lb$ is simple and integrable.
The $\Hb'_{\tau,\zeta}$-modules $\{\Lb,\a\}$ and $\{\Lb',\a'\}$
are isomorphic if and only if there is a linear isomorphism $\Lb\to\Lb'$ 
commuting to $\Hb_{\tau,\zeta}$ up to $x_{\o_0}$,
and taking $\a$ to $\a'$.

Let $c_\pi\,:\,\Hb'_{\tau,\zeta}\to\Hb'_{\tau,\zeta}$ 
be the conjugation by $\pi$.
For any $\Hb'_{\tau,\zeta}$-module $\Mb$ set 
$\Mb^\pi=\Oplus_{i=0}^{\ell-1}(c_\pi^i)^\bullet(\Mb)$, where
$\ell\in\ZZ_{>0}$ is minimal such that $(c_\pi^\ell)^\bullet(\Mb)=\Mb.$
A simple integrable $\Hb^\ch_{\tau,\zeta}$-module may be viewed as a pair
$\{\Mb^\pi,\b\}$, where $\Mb=\{\Lb,\a\}$ is simple and integrable,
and $\b$ is a representation of $\CC[\pi^{\pm 1}]$ in $\Mb^\pi$
compatible with the $\Hb'_{\tau,\zeta}$-action in the obvious way.
If $\b$, $\b'$ are two such representations then 
$\b'(\pi)\circ \b(\pi)^{-1}$ 
is a scalar because $\Mb$ is simple and integrable.

The simple integrable $\Hb_{\tau,\zeta}$-modules are classified by Theorem 7.6. 
Thus the simple integrable $\Hb^\ch_{\tau,\zeta}$-modules, 
modulo the automorphisms $\pi\mapsto z\pi$ for any $z\in\CC^\times$,
are classified by $\Sc/\!\!\sim$ by Proposition 8.2.
While we were writing this paper Cherednik announced a classification  
of the simple integrable $\Hb^\ch$-modules for $\Gu$ of type $A$
and $\zeta^m=\tau^k$, similar to the Bernstein-Rogawski-Zelevinsky 
classification for affine Hecke algebras of type $A$, see \cite{C2}.
The set $\Sc/\!\!\sim$ is precisely the set of new pairs of generalized
$m$-periodic infinite skew-diagrams of total order $d$
entering into Cherednik's classification. 
Cherednik gives also a classification of the simple $\Hb^\ch$-modules
such that the action of the subring $\Rb_A$ is semisimple. 
This last result can not be deduced from our work.

\subhead 8.4\endsubhead
Let $\Cc_d\subset\Cc/\!\!\sim$ be the subset 
consisting of the conjugacy classes
of pairs $(s,e)$ such that $\Spec(s)\subset\tau^{(1/m)\ZZ}$.
The bijection in 8.2 identifies $\Cc_d$ with the subset
$\Sc_d\subset\Sc/\!\!\sim$ consisting of equivalence classes of 
$m$-periodic pairs of order $d$ such that 
$(z_a)$ consists of a single element of $\tau^{(1/m)\ZZ}$.

Let $\Qb_m$ be the quiver of type $A^{(1)}_{m-1}$ with the cyclic orientation,
i.e. the set of vertices is $\ZZ/m\ZZ$
and there is one arrow $i\to i+1$ for each $i\in\ZZ/m\ZZ$.
There is a unique bijection from $\Sc_d$ to the set $Q_d/\!\!\sim$ 
of isomorphism classes of nilpotent representations of $\Qb_m$ of dimension $d$
taking the $m$-periodic multisegment $\sigma$ to
the representation with $\sigma_{ij}$ Jordan blocks
of multi-degree $\{i\,\mod m,i+1\,\mod m,...,j\,\mod m\}$.

Assume that $\Spec(s)\subset\tau^{(1/m)\ZZ}$.
Set ${}^sd_i=\dim\{v\in\CC^d[\eps,\eps^{-1}]\,;\,s(v)=\tau^{i/m}v\}$.
Let $G_1^s\subset G_1$ be the centraliser of $s$.
Note that ${}^sd_{i+m}={}^sd_i$ for all $i$.
Let $\Rep_{{}^sd}(\Qb_m)$ 
be the space of ${}^sd$-dimensional representation of $\Qb_m$.
The group $\GL_{{}^sd}=\prod_{i\in\ZZ/m\ZZ}\GL_{{}^sd_i}$
acts on $\Rep_{{}^sd}(\Qb_m)$.
The bijection $\Cc_d\to Q_d/\!\!\sim$ lifts to
an equivariant isomorphism of varieties
$\{e\,;\,(s,e)\in\Cc_d\}\to\Rep_{{}^sd}(\Qb_m)$.
In particular the closure of the $G^s$-orbit of $e$ in $\Nc^a$
is isomorphic to the closure of 
the corresponding $\GL_{{}^sd}$-orbit in $\Rep_{{}^sd}(\Qb_m)$.
Note that the intersection cohomology of the latter
is expressed by Kazhdan-Lusztig polynomials of type $A^{(1)}$ 
by \cite{L4, Corollary 11.6}. 
The appearence of Kazhdan-Lusztig polynomials
was already mentioned in a less precise way in \cite{AST}.

\subhead 8.5\endsubhead
Let $\Ub_m$ be the specialization at $q=1$ of the generic Hall algebra 
of the category of nilpotent finite dimensional representations of $\Qb_m$.
The $\CC$-algebra $\Ub_m$ admits two natural $\CC$-bases. 
We use the conventions in \cite{VV, Section 3}.
Let $\Sc_\infty=\bigcup_{d>0}\Sc_d.$
The first basis, $\Bb=(\bb_x\,;\,x\in\Sc_\infty)$, 
is given by the Euler characteristic of the stalks of the
cohomology sheaves of the intersection cohomology complex
of the constant local systems on the isomorphism classes. 
The second one, $\Fb=(\fb_x\,;\,x\in\Sc_\infty)$,
is formed by the characteristic functions of the isomorphism 
classes of representations of $\Qb_m$.
For any $m$-uple of non-negative integers $\db$,
let $\fb_\db$ be the element in $\Fb$ labelled by the 
zero representation with graded dimension $\db$.
The algebra $\Ub_m$ is generated by the elements $\fb_\db$,
see \cite{VV, Proposition 3.5}.

Let $\Ub_\infty$ be the negative part of
the enveloping algebra of type $A_\infty$, 
i.e. $\Ub_\infty$ is the $\CC$-algebra
generated by elements $\fbu_i$ with $i\in\ZZ$ 
satisfying the Serre relations.
Let $\Qb_\infty$ be the quiver of type $A_\infty$.
The algebra $\Ub_\infty$ is identified with the specialization at $q=1$
of the generic Hall algebra of the category of finite-dimensional 
representations of $\Qb_\infty$.
It admits two natural bases analoguous to $\Bb$ and $\Fb$.
Let us denote them by $\Bbu=(\bbu_x\,;\,x\in\Scu_\infty)$ and 
$\Fbu=(\fbu_x\,;\,x\in\Scu_\infty)$ respectively.
For any sequence $\dbu$ of non-negative integers, let
$\fb_\dbu$ be the element in $\Fbu$ labelled by the 
zero representation with graded dimension $\dbu$.

Set $\Hb_d=\Hb^\ch_{\tau,\zeta}$ if $d\geq 2$,
and let $\Hb_1$ be the $\CC$-algebra generated 
by $x_1^{\pm 1}$ and $\pi^{\pm 1}$,
modulo the relation $\pi x_i=\tau^{-1}x_i\pi.$
Let $\Gb_d$ be the complexified Grothendieck group 
of the category of integrable $\Hb_d$-modules of finite length
whose Jordan-H\"older factors are labelled by $\Sc_d$.
Set $\Gb_\infty=\bigoplus_{d>0}\Gb_d$,
and $\Tb_\infty=(\Gb_\infty)^*$ (=the restricted dual).
For any $x\in\Sc_d$ let $\Lb_x$ denote the corresponding simple $\Hb_d$-module.
Then $(\Lb_x\,;\,x\in\Sc_\infty)$ is a basis of $\Gb_\infty$.
Let $(\Lb^x\,;\,x\in\Sc_\infty)$ be the dual basis of $\Tb_\infty$.
A conjugacy class $x\in\Cc_d$ contains a pair $(s,e)$ with $e\in\genu$. 
Given a $\slen_2$-triple $\phi\subset\genu$ containing $e$,
let $\Mb_x$ be the class of the module $\Mb_{a,\phi,1}$ in $\Gb_d$.
It depends only on $x$, see 6.7.
Then $(\Mb_x\,;\,x\in\Sc_\infty)$ is a basis of $\Gb_\infty$.
Let $(\Mb^x\,;\,x\in\Sc_\infty)$ be the dual basis.

For any $\ell\in\ZZ$
let $\Ccu_d[\ell]$ be the set of conjugacy classes of pairs $(\su,\eu)$
such that $\eu\in\slen_d$ is nilpotent, $\su\in\GL_d$ is semisimple
with $\Spec(\su)\subset\tau^{\ell/m}\zeta^\ZZ$, and $\ad_\su\eu=\zeta\eu$.
We write $\Ccu_d$ for $\Ccu_d[0]$. 
The set $\Scu_d$ of isomorphism
classes of $d$-dimensional representations of $\Qb_\infty$ 
is naturally identified with $\Scu_d$.
Set $\Scu_\infty=\bigcup_{d>0}\Scu_d.$
Let $\Hbu_d\subset\Hb_d$ be the subalgebra generated by
$x_1^{\pm 1}$,...,$x_d^{\pm 1}$,$t_1$,...,$t_{d-1}$.
Let $\Gbu_d$ be the complexified Grothendieck group 
of the category of finite dimensional $\Hbu_d$-modules 
whose Jordan-H\"older factors are labelled by $\Scu_d$.
Set $\Gbu_\infty=\bigoplus_{d>0}\Gbu_d$,
and $\Tbu_\infty=(\Gbu_\infty)^*.$ 
For any $x\in\Scu_d$ let $\Lbu_x$ 
denote the corresponding simple $\Hbu_d$-module.
Let $(\Lbu^x\,;\,x\in\Scu_\infty)$ be the dual basis of $\Tbu_\infty$.
Given a pair $(\su,\eu)$ in a conjugacy class in $\Ccu_d$
and a $\slen_2$-triple $\phiu\subset\genu$ containing $\eu$,
let $\Mbu_x$ be the class in $\Gbu_d$
of the module $\Mbu_{(\su,\zeta),\phiu,1}$.
Then $(\Mbu_x\,;\,x\in\Scu_\infty)$ is a basis of $\Gbu_\infty$.
Let $(\Mbu^x\,;\,x\in\Scu_\infty)$ be the dual basis.

There is a unique $\CC$-linear isomorphism
$\Ub_m\to\Tb_\infty$ 
such that $\Lb^x\mapsto\bb_x$ and $\Mb^x\mapsto\fb_x$,
and a unique $\CC$-linear isomorphism $\Ub_\infty\to\Tbu_\infty$ 
such that $\Lbu^x\mapsto\bbu_x$ and $\Mbu^x\mapsto\fbu_x$.
There is a unique algebra homomorphism
$$\Delta\,:\,\Ub_m\to\hat\Ub_\infty\ 
\roman{such\ that\ }\fb_\db\mapsto\Sum_{\dbu}\fbu_\dbu,$$
where $\db=(d_j)$ and the sum is over all $\dbu=(\du_j)$ such that
$d_i=\sum_{j=i\,\mod m}\du_j$, see \cite{VV, Remark 6.1}.
Here $\hat\Ub_\infty$ is a completion of $\Ub_\infty$.
Given $x\in\Sc_d$ and $\xu\in\Scu_d$ let
$m_{x,\xu}$ be the coefficient of $\bbu_\xu$ in the infinite sum
$\Delta(\bb_x)$.

\proclaim{Theorem}
(i)
If $\Lbu$ is a simple $\Hbu_d$-module and 
$\Lb_x$ is a Jordan-H\"older factor in $\Hb_d\otimes_{\Hbu_d}\Lbu$, 
then the Langlands parameter of 
$\Lbu$ belongs to $\Ccu_d[\ell]$ for some $\ell$.

(ii)
$m_{x,\xu}$ is the Jordan-H\"older multiplicity of
$\Lb_x$ in $\Hb_d\otimes_{\Hbu_d}\Lbu_\xu.$ 
\endproclaim

\noindent{\sl Proof:}
The induction from $\Hbu_d$-modules to $\Hb_d$-modules is an exact functor. 
It yields a linear map $\Gamma\,:\,\Gbu_\infty\to\Gb_\infty$.
Let $\g\,:\,\Scu_\infty\to\Sc_\infty$ be the map
taking the $\Gu$-conjugacy class of $(\su,\eu)$
to the $G_1$-conjugacy class of $((\su,\tau),\eu)$.
We have $\Gamma(\Mbu_\xu)=\Mb_{\g(\xu)}$ by Lemma 6.6.
Let $\delta\,:\,\Tb_\infty\to\hat\Tbu_\infty$ be the map dual to $\Gamma$,
where $\hat\Tbu_\infty$ is a completion of $\Tbu_\infty$.
Hence $\delta(\fb_x)=\sum_{\g(\xu)=x}\fbu_\xu$.
The square 
$$\matrix
\Tb_\infty&{\buildrel\delta\over\lra}&\hat\Tbu_\infty\cr
\Vert&&\Vert\cr
\Ub_m&{\buildrel\Delta\over\lra}&\hat\Ub_\infty
\endmatrix$$
is commutative, because $\Delta(\fb_x)=\sum_{\g(\xu)=x}\fbu_\xu$
by \cite{VV, Lemma 6.4}.
Therefore $m_{x,\xu}$ is the Jordan-H\"older multiplicity of
the simple $\Hb_d$-module $\Lb_x$ in
the induced module $\Gamma(\Lbu_\xu).$ 
Let $(\su,\eu)$ be the Langlands parameters of $\Lbu_\xu$.
Claim $(i)$ is obvious because the Langlands parameters 
$(s,e)$ of a Jordan-H\"older factor of $\Hb_d\otimes_{\Hbu_d}\Lbu_\xu$ 
are such that $\Spec(\su)=\Spec(s)\,\mod\tau^{\ZZ}$.
\qed

\vskip3mm

\noindent{\bf Remark.}
It would be interesting to give a representation theoretic
construction of the product in $\Ub_m$ via the Grothendieck group
$\Gb_\infty$ as for the affine case, see \cite{A} for instance.

\vskip3mm

%\noindent{\bf Examples.}
%$(i)$
%If $d=1$ the map $\Sc_1\to\tau^{(1/m)\ZZ}/\tau^\ZZ$,
%$\{\sigma,z\}\mapsto z\tau^\ZZ$ is a bijection.
%The simple $\Hb_1$-module corresponding to $z\tau^\ZZ$, $z\in\tau^{(1/m)\ZZ}$,
%is the unique simple integrable module such that the spectrum of $x_1$ is
%$z\tau^\ZZ.$ 
%
%$(ii)$
%If $m=d$ and $k=\pm 1$ there is a unique one-dimensional 
%$\Hb'_{\tau,\zeta}$-module
%such that $t_i\mapsto -1$, $x_i\mapsto\zeta^{-i}$,
%and $t_i\mapsto\zeta$, $x_i\mapsto\zeta^i$ respectively.
%Let $\Lb_\pm$ be the corresponding simple $\Hb_d$-modules.
%For any $j\in\ZZ$ let $\Lbu_\pm[j]$ be the unique one-dimensional 
%$\Hbu_d$-modules such that $t_i\mapsto -1$, $x_i\mapsto\zeta^{j-i}$,
%and $t_i\mapsto\zeta$, $x_i\mapsto\zeta^{j+i}$ respectively. 
%The restriction of $\Lb_\pm$ to $\Hbu_d$ is 
%$\bigoplus_{j\in\ZZ}\Lbu_\pm[j]$.

\head 9. Other examples\endhead
Let $\Modb_a=\Psi_a^\bullet(\Modb_{\Kb_a})\subset\Modb_{\Hb_{\zeta,\tau}}$.
Let $\Modbu_a$ be the category of finite-dimensional 
$\Hbu_\zeta$-modules whose central character belongs to 
the orbit $W(s')$.

\subhead 9.1. The generic case\endsubhead
Assume that 
$1,\zeta\notin\{\tau^n\su^\a\,;\,\a\in\Phiu,n\in\ZZ\setminus\{0\}\}$.
Then the induction functor 
$\Lb\to\Hb^\ch_{\zeta,\tau}\otimes_{\Hbu_\zeta}\Lb$
takes simple modules in $\Modbu_a$ to simple modules in $\Modb_a$
by \cite{C1, Proposition 6.6}. 
The proof of Cherednik uses intertwining operators.
A geometric proof of the same result is as follows : 
we have $\Pi(a)^\vee=\Piu(a)^\vee$ and
$n_{a,\chi,\chi'}={\underline n}_{a,\chi,\chi'}$,
because $G^s=\Gu^s$ and $\Nc^a=\Ncu^a$;
hence the claim follows from Lemmas 6.5, 6.6 and 6.7.

\subhead 9.2. The case $s'=1$\endsubhead
Set $s'=1$.
The Fourier transform gives an equivalence of categories
$\Modb_a\to\Modb_{\bar a}$ by 7.3.
Hence we may assume that $(\tau,\zeta)$ is regular.
Then $\Nc^a=\{0\}$ unless $\zeta=\tau^k$ with $k\in\ZZ_{<0}$,
and we get $\Nc^a=\Ncu\otimes\eps^k$, $G^{s}=H\cdot\Gu$.
Therefore the map $x\mapsto x\otimes\eps^k$ gives a bijection from
the set of $\Gu$-orbits in $\Ncu$ to
the set of $G^s$-orbits in $\Nc^a$.
Conjecturally, $\Pi(a)^\vee$ is equal to 
the set of irreducible $\Gu$-equivariant local systems 
on $\Ncu$ which occur in the homology of the Springer fibers in $\Bcu$.
Then, if $\Gu$ is a classical group, there would be a natural bijection
between the simple $\Wu$-modules
and the simple objects in $\Modb_a$ by \cite{S, Theorem 2.5}.
However the category $\Modb_a$ is not semisimple.
Probably $\Modb_a$ is equivalent to 
the category of $W$-modules with trivial central character.
Note also that the induction functor 
$\Modbu_a\to\Modb_a$
is trivial because $\Modbu_a$ is trivial.

\subhead 9.3. Finite dimensional modules\endsubhead
Fix $k,m\in\ZZ_{>0}$ with $(m,k)=1$.
Assume that $\zeta^m=\tau^k$.
Fix $e\in\Nc^a$. 
The nil-element $e$ may be not nilpotent.
Recall that $e$ is regular semisimple {\sl nil-elliptic} if
it is regular semisimple and
the group of cocharacters of the centralizer of $e$ in $\Gu(K)$ is trivial. 
Then $\Bc_e$ is an algebraic variety by\cite{KL2, Corollary 3.2}.
Thus $\Nb_{a,e,\chi}$ is finite dimensional for any $\chi$.

Let $\h\t\,:\,\Phiu\to\ZZ$ be the height, and $h=\h\t(\theta)+1$.
Let $r$ be the rank of $\genu$, and $\rho^\vee=\sum_{i\in\Iu}\o_i^\vee$.
Hence $2\rho^\vee\in\Yu^\vee$.
Fix $a,b\in\ZZ_{\geq 0}$ with $b<h$ and $k=ah+b$.
Set $\Pi_k=\{(\a,\ell)\in\Phiu\times\ZZ\,;\,(\h\t(\a),\ell)=(b,a),(b-h,1+a)\}.$
For any $\a\in\Phiu$ we fix a non-zero element $e_\a\in\genu_\a$.
Assume that $m=h$.
Fix a $2h$-th root of $\tau$ such that $\zeta=\tau^{k/h}$. Put 
$$e_k=\sum_{(\a,\ell)\in\Pi_k}e_\a\otimes\eps^\ell,
\quad
\su=(2\rho^\vee)\otimes\tau^{1/2h},
\quad
s=(\su,\tau),
\quad
a_k=(s,\zeta).$$

\proclaim{Lemma}
(i)
$\gen^{a_k}$ is the linear span of 
$\{e_\a\otimes\eps^\ell\,;\,(\a,\ell)\in\Pi_k\}$,
and $G^s=H$.

(ii)
There is a finite number of $G^s$-orbits in $\gen^{a_k}$.
The orbit containing $e_k$ is dense.
\endproclaim

\noindent{\sl Proof:}
The first part of $(i)$ is immediate because
$$\ad_s(e_\a\otimes\eps^\ell)=\tau^{\ell+\h\t(\a)/h}e_\a\otimes\eps^\ell$$
for all $\a$, $\ell$.
For the second part it is sufficient to observe that
$G^s$ is generated by $H$ and the $U_\b$'s with $\gen^s_\b\neq\{0\}$
because it is connected by Lemma 2.13, and that
$\ad_s(e_\a\otimes\eps^\ell)\neq e_\a\otimes\eps^\ell$ for all
$(\a,\ell)\in\Phiu\times\ZZ$ because $\tau,\zeta$ are not roots of unity.

Let us prove $(ii)$.
Assume first that $k=1$.
Recall that $H=\Hu\times\CC^\times_\delta\times\CC^\times_{\o_0}$ with
$\CC^\times_{\o_0}$ acting trivially on $\gen$.
We have $\Pi_1=\{\a_i\,;\,i\in I\}$.
Thus we must count the orbits of a $(r+1)$-torus acting 
on a $(r+1)$-vector space. 
The set of weights of $\gen^{a_1}$ is $\Pi_1$,
and $\Lie(\Hu\times\CC^\times_\delta)$ is linearly spanned by 
$\{\o_i^\vee\,;\,i\in I\}$.
Thus there are $2^{r+1}$ orbits.
Let $R_k\subset\Phiu$ be the image of $\Pi_k$ by the first projection.
By \cite{F, Lemma 2.2} there is an element $w\in\Wu$ such that $w(R_k)=R_1$.
Set $\rho^\vee=\sum_{i\in\Iu}\o_i^\vee\in\Yu^\vee$, and
$g=\rho^\vee\otimes\eps^{1/h}\in\Gu(K_h)$.
Then
$\dot w\circ\ad_g(\gen^{a_k}),
\ad_g(\gen^{a_1})\subset\genu\otimes K_h$
are linearly spanned by
$\{e_\a\otimes\eps^{k/h}\,;\,\a\in R_1\}$,
$\{e_\a\otimes\eps^{1/h}\,;\,\a\in R_1\}$
respectively.
Thus there is the same number of $H$-orbits
in $\dot w\circ\ad_g(\gen^{a_k})$ and $\ad_g(\gen^{a_1})$.
The lemma follows. 
\qed

\vskip3mm

We have $\Nc^{a_k}=\gen^{a_k}$ by Proposition 2.14.$(iii)$.
Since there are a finite number of $G^s$-orbits in $\Nc^{a_k}$
the discussion before Proposition 6.3 still holds.
In particular the simple perverse sheaves $S_{a_k,\chi}$ 
occuring in $L_{a_k}$ are intersection cohomology complexes of irreducible
$G^s$-equivariant local systems on $\Nc^{a_k}$,
and the corresponding representations of $\Pi(s,e)$ 
must occur in $\Hb_*(\Bc_e^s)$ for any $e\in\Nc^{a_k}$.
It can be shown that $\Bc^s=\{\ben_w\,;\,w\in W\}$
by using the Bruhat decomposition for $G'$.
Hence $G^s$ acts trivially on $\Bc^s$,
because it is connected. 
Thus $\Pi(s,e)$ acts trivially on $\Hb_*(\Bc_e^s)$.
Therefore the local systems above are all trivial.
Moreover
the $\Hb_{\tau,\zeta}$-module $\Hb_*(\Bc_{e_k}^s)$
is simple because the orbit of $e_k$ is dense.
Obviously $\Bc^s_{e_k}$ has no odd cohomology. 
Thus $\dim\Hb_*(\Bc^s_{e_k})=k^r$ by \cite{F, Proposition 1}. 
Therefore we get the following.

\proclaim{Proposition}
$\Hb_*(\Bc^s_{e_k})$ is a simple $\Hb_{\tau,\zeta}$-module
of dimension $k^r$.
\endproclaim

\noindent{\bf Remarks.} 
$(i)$
Probably $\Hb_*(\Bc^s_{e_k})$ and
$\IM^\bullet(\Hb_*(\Bc^s_{e_k}))$, with $k>0$, $m=h$, and $(k,h)=1$,
are the only finite dimensional simple $\Hb$-modules
if $\tau$, $\zeta$ are not roots of unity and $\Gu$ is of type $A$.
This is false in general : for instance, if $\Gu=\SP(2r)$,
$e=\sum_{(\a,\ell)\in\Pi_1}e_\a\otimes\eps^{r\ell}$,
$\su=(2\rho^\vee)\otimes\tau^{1/4}$,
and $\zeta=\tau^{1/2}$, then $\ad_se=\zeta e$
and the element $e$ is regular semisimple nil-elliptic.
This yields a finite dimensional $\Hb_{\tau,\zeta}$-representation on
$\Hb_*(\Bc^s_e)$ which can not be of the previous type because $h=2r$.

$(ii)$
A direct computation gives $\Pi(s,e_k)=\ZZ/k\ZZ$ if $\Gu=\SL_{r+1}$.

\subhead 9.4. Comparison with \cite{L5}\endsubhead
Let us compare our construction with
Lusztig's $W$-action on the homology of the affine Springer fibers.
This comparison is not used in the rest of the paper.
Proofs are not given.
The $W$-action in \cite{L5} is described via a familly of compatible
$W_J$-actions, for each $J\subsetneq I$. 

Let $\Hb'_J$ be the graded Hecke algebra associated to $\Hb_J$, 
see \cite{EM} for instance.
There is an obvious surjective $\CC$-algebra homomorphism $\Hb'_J\to\CC W_J$.

Let $\Nc_J$ the nilpotent cone of $G_J$,
and $p_J\,:\,\dot\Nc_J\to\Nc_J$ be the Springer map.
There is a commutative square of algebra homomorphisms
$$\matrix
\Hb'_J&\simto&\Extb_{D_{G_J\times\CC^\times_q}(\Nc_J)}
\bigl((p_J)_!\CC,(p_J)_!\CC\bigr)
\hfill\cr
\bda&&\bda\cr
\CC W_J&\simto&\Homb_{D_{G_J\times\CC^\times_q}(\Nc_J)}
\bigl((p_J)_!\CC,(p_J)_!\CC\bigr).\hfill\cr
\endmatrix$$
Set 
$${}^w\Tc_J^{k\ell}=
(U^-_k\setminus{}^w\Gc)\times_{B_J}(\Nc_J+\uen_J/\uen_\ell),\quad
{}^w\Tc_J={}^w\Gc\times_{B_J}(\Nc_J+\uen_J),$$
with $\uen_\ell=\eps^\ell\cdot\uen_J$, $\ell,k\geq 0$, and $w\in W^Jw_J$.
The group $A$ acts on ${}^w\Tc_J^{k\ell}$, ${}^w\Tc^{k\ell}$.
The obvious projection
$p\,:\,{}^w\Tc^{k\ell}\to{}^w\Tc_J^{k\ell}$ 
is proper, and commutes with $A$-action.
Fix a non empty, separated, open subset $\Uc_J\subset{}^w\Tc_J^{k\ell}$ 
which is preserved by the $A$-action.
Set $\Uc=p^{-1}\Uc_J.$ 
Let $p_\Uc$ be the restriction of $p$ to $\Uc$.
From \cite{L5, 5.3} we get a graded algebra homomorphism
$$\Extb_{D_{G_J\times\CC^\times_q}(\Nc_J)}
\bigl((p_J)_!\CC,(p_J)_!\CC\bigr)
\to\Extb_{D_{A}(\Uc_J)}\bigl((p_\Uc)_!\CC,(p_\Uc)_!\CC\bigr).$$
The right hand side is $\Hb_*^A(\Zc_\Uc)$, where
$\Zc_\Uc={}^w\Zc^{k\ell}_{\leq w_J}\cap\Uc^2$.
We get the commutative square of algebra homomorphisms
$$\matrix
\Hb'_J&\lra&
\Hb_*^A(\Zc_\Uc)
\hfill\cr
\bda&&\bda\cr
\CC W_J&\lra&
\Hb_{top}(\Zc_\Uc),
\hfill\cr
\endmatrix\leqno(9.4.1)$$
where the subscript `top' means that we only consider 
the homology groups of maximal degrees.
The right vertical map is evaluation at 0 in $\Lie(A)$.

Fix $e\in\Nc$.
Set $\Bc_{e,\leq w}=\Bc_e\cap\Bc_{\leq w}$.
Fix $w$ such that $\Bc_{e,\leq w}\neq\emptyset$.
Fix $k$ such that $\Bc_{e,\leq w}\subset{}^w\Tc^{k\ell}$.
Fix $\Uc_J$ such that $\Bc_{e,\leq w}\subset\Uc$.
The ring $(\Hb_{top}(\Zc_\Uc),\star)$ acts on 
$\Hb_*(\Bc_{e,\leq w})$,
yielding a $W_J$-action on
$\Hb_*(\Bc_e)=\ind_w\Hb_*(\Bc_{e,\leq w})$.
This is Lusztig's action.

Let $a=(s,\zeta)$ be as in 2.12.
Fix $k,\ell$ and $\Uc_J$ such that 
$$\Uc^a=({}^w\Tc)^a=({}^w\Tc^{k\ell})^a,$$
see 2.13.
Note that $\Uc$ is a smooth scheme of finite type.
Composing the upper horizontal map in (9.4.1) with the concentration map
relative to $\Uc^2$ 
we get an algebra homomorphism
$$\psi\,:\,\Hb'_J\to\Hb_*(\Zc_{\leq w_J}^a).$$
Hence $\Hb'_J$ acts on $\Hb_*(\Bc^s_e)$ by convolution whenever $e\in\Nc^a$.
In 6.2 we defined a $\Hb_J$-action on $\Hb_*(\Bc_e^s).$

\proclaim{Claim}
There are increasing $\ZZ_{\geq 0}$-filtrations on $\Hb_J$ and 
$\Hb_*(\Bc_e^s)$ such that the graded ring associated to
$\Hb_J$ is $\Hb'_J$, and the representations of
$\Hb_J$ and $\Hb'_J$ on $\Hb_*(\Bc_e^s)$ coincide modulo 
lower terms in the filtration. 
\endproclaim

\noindent{\bf Remark.}
Clearly the representations of
$\Hb'_J$ on $\Hb_*(\Bc^s_e)$ are compatible, when $J$ varies,
yielding a representation of the degenerated double affine Hecke algebra.

\head A. Appendix : K-theory\endhead
\subhead A.1. Generalities\endsubhead
Let $G$ be a linear group.
In A.1 all schemes are Noetherian.
Given a $G$-scheme $\Xc$, we say that $(G,\Xc)$ has the resolution property
if any $G$-equivariant coherent sheaf on $\Xc$ is the quotient of
a $G$-equivariant locally free sheaf,
and the quotient map commutes with the $G$-action.
It is known that $(G,\Xc)$ has the resolution property
if $\Xc$ is smooth and separated, see \cite{T1}.
If the scheme $\Xc$ is smooth and not separated 
the pair $(G,\Xc)$ has not the resolution property in general,
see \cite{TT, Exercice 8.6}.
If $(G,\Xc)$ has the resolution property and $\Xc$ is smooth
then the natural group homomorphism
$\Kb^i_G(\Xc)\to\Kb^G_i(\Xc)$ is an isomorphism.
In the following we give a reminder on equivariant K-theory
for non separated schemes. 

\vskip2mm

\noindent(A.1.1)
Let $f\,:\,\Pc\to\Xc$ be a $G$-equivariant map 
of $G$-schemes which is a $H$-torsor, 
with $H\subset G$ a closed normal subgroup.
Then $Lf^*$ is an isomorphism 
$\Kb^{G/H}_i(\Xc)\to\Kb^{G}_i(\Pc)$, see \cite{T3}.
In particular if $G=G_1\times G_2$ and $\Yc$ 
is a $G_1$-scheme,
then $\Pc\times_{G_1}\Yc$ is a 
$G_2$-scheme and there is the map
$$\Kb_i^{G_1}(\Yc)\to\Kb_i^{G}(\Pc\times\Yc)\to
\Kb_i^{G_2}(\Pc\times_{G_1}\Yc),$$
such that the first arrow takes $E$ to $\CC\boxtimes E$, 
and the second arrow is the inverse of the pull-back by the projection
$\Pc\times\Yc\to\Pc\times_{G_1}\Yc$.
The homotopy invariance and
the existence of localization long exact sequences
are proved in \cite{T3}.

\vskip2mm

\noindent(A.1.2) The following lemma and its proof 
are taken from the proof of \cite{T2, Lemma 4.3}.

\proclaim{Lemma}
$\Xc$ contains an open dense subset which is $G$-invariant and separated.
\endproclaim

\noindent{\sl Proof:}
Fix a dense open affine subset $\Uc\subset\Xc$.
The map $G\times\Xc\to\Xc$, 
$(g,x)\mapsto g\cdot x$ is flat of finite type.
Hence it is open. 
Thus $G\cdot\Uc\subset\Xc$ is an open dense subset.
The open dense subset $G^\circ\cdot\Uc\subset\Xc$ is separated 
by the following argument taken from \cite{R, Lemma V.3.11}
that we reproduce for the comfort of the reader. 
Assume that $S$ is the spectrum of a discrete valuation field
whith generic point $x$ and closed point $y$,
and assume that $f_1,f_2$ are maps $S\to G^\circ\cdot\Uc$ such that
$f_1(x)=f_2(x)$.
Set $U_i=\{g\in G^\circ\,;\,g\cdot f_i(y)\in\Uc\}$.
They are non-empty open subsets of $G^\circ$.
Hence $U_1\cap U_2\neq\emptyset$.
Fix $g\in U_1\cap U_2$.
Then $f_1(y), f_2(y)\in g^{-1}\cdot\Uc$.
Thus $f_1=f_2$ because $g^{-1}\cdot\Uc$ is separated.
Since $G/G^\circ$ is finite, $G\cdot(\Xc\setminus G^\circ\cdot\Uc)$ is a finite union
of nowhere dense closed subspaces of $\Xc$.
Hence $\Xc\setminus G\cdot(\Xc\setminus G^\circ\cdot\Uc)$ is a $G$-invariant open
dense subset. It is separated because it is contained in $G^\circ\cdot\Uc$.
\qed

\vskip2mm

\noindent(A.1.3)
Assume that the group $G$ is diagonalizable and $g\in G$.
Recall that $\Kb_j^G(\Xc)_g$ is the specialization of the $\Rb_G$-module
$\Kb_j^G(\Xc)$ at the maximal ideal associated to $g$, see 1.4.

\proclaim{Theorem}
If the natural inclusion
$f\,:\,\Xc^g\to\Xc$ is a closed immersion, then the direct image
$Rf_*\,:\,\Kb_i^G(\Xc^g)_g\to\Kb_i^G(\Xc)_g$ is an isomorphism.
\endproclaim

\noindent{\sl Proof:}
The proof of \cite{T4, Theorem 2.1} extends easily to the case
of non separated schemes. 
Since $f$ is a closed immersion
and localization long exact sequences do exist, 
we are reduced to check that $\Kb_i^G(\Xc)_g=0$ for any $G$-scheme
such that $\Xc^g=\emptyset$.
To do so, it is sufficient to construct a non empty open subset 
$\Uc\subset\Xc$ which is preserved by the $G$-action
and such that $\Kb_i^G(\Uc)_g=0$ because, then,
a Noetherian induction on closed $G$-subsets of $\Xc$ implies
that $\Kb_i^G(\Xc\setminus\Uc)_g=0$, and thus
$\Kb_i^G(\Xc)_g=0$ by the localization long exact sequence.
By Lemma A.1.2, $\Xc$ contains an open dense subset $\Uc$ which is
$G$-invariant and separated.
Hence $\Kb_i^G(\Uc)_g=0$ by the concentration theorem for 
separated $G$-schemes.
We are done.
\qed

\vskip2mm

Assume that $\Xc$ is smooth and $f$ is a closed immersion.
Then $\Xc^g$ is smooth, possibly non separated,
because smoothness is a local condition and
$\Xc$ can be covered by $G$-invariant open and separated subschemes.
Let $N$ be the conormal bundle of $\Xc^g$ in $\Xc$.
The element $\l(N)=\sum_{i\geq 0}(-1)^i\bigwedge^i(N)\in\Kb_G(\Xc^g)_g$
is invertible.
For any $G$-equivariant vector bundle $E$ on $\Xc$ we have
$f_*(\l(N)^{-1}\otimes(E|_{\Xc^g}))=E$ in $\Kb^G(\Xc)_g$.

\vskip2mm

\noindent(A.1.4)
The tor-product 
$\otimes^L\,:\,\Kb^G(\Xc)\times\Kb^G(\Xc)\to\Kb^G(\Xc)$
is well-defined if $\Xc$ is smooth, possibly non separated.
A naive proof is that the local Tor exist because
$\Xc$ can be covered by $G$-invariant separated open subsets, see above.
More rigorously, let $\Qcohb^G(\Xc)$ be the Abelian category of 
$G$-equivariant quasi-coherent sheaves on $\Xc$. 
Let $\Db(\Qcohb^G(\Xc))$ be its derived category.
Let $\Db^b(\Qcohb^G(\Xc))\subset\Db(\Qcohb^G(\Xc))$
be the full subcategory whose objects are the complexes 
with bounded cohomology.
Let also $\Db^G(\Xc)_{pf}\subset\Db^b(\Qcohb^G(\Xc))$
be the full subcategory whose objects are perfect complexes
(i.e. complexes locally quasi-isomorphic to a bounded complex
of locally free sheaves of finite type), and let
$\Db^G(\Xc)_{coh}\subset\Db^b(\Qcohb^G(\Xc))$
be the full subcategory whose objects are the complexes 
with coherent cohomology.
Then $\Kb^G(\Xc)=\Kb(\Db^G(\Xc)_{coh})$ because $\Xc$ is Noetherian, 
and $\Kb(\Db^G(\Xc)_{coh})=\Kb(\Db^G(\Xc)_{pf})$ because
$\Xc$ is smooth and Noetherian, see \cite{SGA6, Ch IV, \S 2.4-5}.
Note that \cite{SGA6} considers only the non-equivariant case. 
The general case follows from the resolution property 
of smooth separated $G$-schemes 
and the fact that $\Xc$ can be covered
by $G$-invariant separated open subsets by Lemma A.1.2.
Since the derived tensor product 
is well-defined on $\Db^G(\Xc)_{pf}$, 
we have a derived tensor product on $\Kb^G(\Xc)$.

\vskip2mm

\noindent(A.1.5)
Let $\Yc_1,\Xc_1$ be $G$-schemes
with $\Yc_1\subset\Xc_1$ closed.
Let $f\,:\,\Xc_2\to\Xc_1$ be a closed immersion of $G$-schemes,
and $\Yc_2=f^{-1}(\Yc_1)\subset\Xc_2$.
If $\Xc_1$ is smooth then the pull-back 
$Lf^*\,:\,\Kb^G(\Yc_1)\to\Kb^G(\Yc_2),\,
\Ec\mapsto\Oc_{\Xc_2}\otimes^L\Ec$,
respectively to the tor-product on $\Xc_1$,
is well-defined. 

\vskip2mm

\noindent{\bf Other conventions.}
We may identify a sheaf on $\Xc$ and its class in $\Kb^G(\Xc)$. 
We write $\otimes$, $f_*$, $f^*$ instead of $\otimes^L$, $Rf_*$, $Lf^*$.
To avoid confusion we may also write $\otimes_\Xc$ for $\otimes$.

\subhead A.2. The $\star$-product\endsubhead
Let $\Sc$ be a smooth Noetherian $G$-scheme.
Fix closed $G$-subschemes 
${}^{12}\tilde\Zc, {}^{13}\tilde\Zc, {}^{23}\tilde\Zc\subset\Sc$
such that ${}^{12}\tilde\Zc\cap{}^{23}\tilde\Zc\subset{}^{13}\tilde\Zc$.
Let $g_{ij}\,:\,{}^{ij}\tilde\Zc\to{}^{ij}\Zc$ be smooth $G$-maps
such that the restriction of $g_{13}$ to 
${}^{12}\tilde\Zc\cap{}^{23}\tilde\Zc$ is proper.
The $\star$-product relative to $\Sc$ is
$$\Kb^G({}^{12}\Zc)\times\Kb^G({}^{23}\Zc)\to\Kb^G({}^{13}\Zc),\ 
(E, F)\mapsto g_{13*}(g_{12}^*E\otimes_\Sc g_{23}^*F).$$

The following situation is standard : 
$\Tc$ is a smooth Noetherian $G$-scheme,
${}^{12}\Zc,{}^{13}\Zc,{}^{23}\Zc$ 
are closed subsets in $\Tc^2$
preserved by the diagonal $G$-action,
$\Sc=\Tc^3$,
and $q_{ij}\,:\,\Sc\to\Tc^2$
is the projection along the factor not named such that 
$q_{13}$ restricts to a proper map
$q_{12}^{-1}({}^{12}\Zc)\cap q_{23}^{-1}({}^{23}\Zc)\to{}^{13}\Zc$.
Then the $\star$-product 
relative to $\Sc$ takes 
$(E, F)$ to $q_{13*}(q_{12}^*E\otimes_\Sc q_{23}^*F).$
If ${}^{12}\Zc={}^{23}\Zc={}^{13}\Zc=\Zc$ then
$(\Kb^G(\Zc),\star)$ is a $\Rb_G$-algebra.
This example adapts easily to the case where
$\Sc$  is the product of three distinct smooth schemes.
This generalization is left to the reader.

If $G=\{1\}$, 
$p\,:\,\Tc\to\Nc$ is proper with $\Tc$, $\Nc$ separated, 
and  $\Zc=\Tc\times_\Nc\Tc$, 
there is a linear map 
$\Phi\,:\,\Kb_0(\Zc)\to\Extb_{\Db(\Nc)}(p_*\CC_\Tc,p_*\CC_\Tc)$ taking
the $\star$-product to the Yoneda product.
It is the composition of the bivariant Riemann-Roch map
$\R\R\,:\,\Kb_0(\Zc)\to\Hb_*(\Zc)$ relative to $\Tc^2$, 
see \cite{CG, 5.11}, and of the chain of isomorphisms
$$\Hb_*(\Zc)\dot{=}\Hb^*(\Zc,\DD_\Zc)\dot{=}
\Extb_{\Db(\Zc)}(q_1^*\CC_\Tc,q_2^!\CC_\Tc)=$$
$$=\Extb_{\Db(\Tc)}(q_{2!}q_1^*\CC_\Tc,\CC_\Tc)=
\Extb_{\Db(\Tc)}(p^*p_!\CC_\Tc,\CC_\Tc)=
\Extb_{\Db(\Nc)}(p_*\CC_\Tc,p_*\CC_\Tc).$$
Moreover $\Phi$ is invertible if $\R\R$ is invertible.

\subhead A.3. Functoriality of $\star$\endsubhead
For each $k=1,2$ we fix
$\Tc_k$, $\Zc_k\subset(\Tc_k)^2$, and
$\Sc_k=(\Tc_k)^3$ as in A.2.
Let $f\,:\,\Tc_2\to\Tc_1$ be a flat $G$-map
such that $(f\times f)(\Zc_2)\subset\Zc_1$.

\vskip2mm

\noindent(A.3.1)
Assume that $f\times\Id$ restricts to a proper map 
$\Zc_2\to(\Id\times f)^{-1}(\Zc_1)$.
Then the maps 
$(f\times\Id)_*\,:\,\Kb^G(\Zc_2)\to
\Kb^G((\Id\times f)^{-1}(\Zc_1))$
and
$(\Id\times f)^*\,:\,\Kb^G(\Zc_1)\to
\Kb^G((\Id\times f)^{-1}(\Zc_1))$
are well-defined.
If 
$(f\times\Id)_*(E_2)=(\Id\times f)^*(E_1),$
and
$(f\times\Id)_*(F_2)=(\Id\times f)^*(F_1),$
then
$$\matrix
(f\times\Id)_*(E_2\star F_2)
&=(f\times\Id)_*q_{13*}
((E_2\boxtimes\CC)\otimes(\CC\boxtimes F_2)),\hfill\cr
&=q_{13*}(f\times f\times\Id)_*((E_2\boxtimes\CC)\otimes
(f\times\Id\times\Id)^*(\CC\boxtimes F_2)),\hfill\cr
&=q_{13*}(\Id\times f\times\Id)_*
((\Id\times f\times\Id)^*(E_1\boxtimes\CC)
\otimes(\CC\boxtimes F_2)),\hfill\cr
&=q_{13*}((E_1\boxtimes\CC)\otimes
(\Id\times\Id\times f)^*(\CC\boxtimes F_1)),\hfill\cr
&=q_{13*}(\Id\times\Id\times f)^*
((E_1\boxtimes\CC)\otimes(\CC\boxtimes F_1)),\hfill\cr
&=(\Id\times f)^*q_{13*}
((E_1\boxtimes\CC)\otimes(\CC\boxtimes F_1)),\hfill\cr
&=(\Id\times f)^*(E_1\star F_1).\hfill
\endmatrix$$
The 3-rd and 4-th equalities follow from the projection formula,
the 6-th equality follows from base change.

If $(f\times\Id)_*$ is invertible we set
$f^\sharp=(f\times\Id)_*^{-1}(\Id\times f)^*\,:\,
\Kb^G(\Zc_1)\to\Kb^G(\Zc_2)$,
if $(\Id\times f)^*$ is invertible we set 
${}_\sharp f=
(\Id\times f)^{*-1}(f\times\Id)_*\,:\,
\Kb^G(\Zc_2)\to\Kb^G(\Zc_1)$.
Both maps are $\star$-homomorphisms.

\vskip2mm

\noindent{\bf Examples.}
$(i)$
If $\Tc_1=\Tc_2$, $f=\Id$,
and $\Zc_2$ is closed in $\Zc_1$ 
then ${}_\sharp f$ makes sense.

$(ii)$
If $f$ is open immersion such that
$\Zc_2=(\Id\times f)^{-1}(\Zc_1)$
then $f^\sharp$ makes sense.

\vskip2mm

\noindent(A.3.2)
If $\Id\times f$ restricts to a proper map 
$\Zc_2\to(f\times\Id)^{-1}(\Zc_1)$ and 
$(\Id\times f)_*(E_2)=(f\times\Id)^*(E_1),$
$(\Id\times f)_*(F_2)=(f\times\Id)^*(F_1),$
then $E_1\star F_1$, $E_2\star F_2$ satisfy similar relations. 

If $(\Id\times f)_*$ is invertible we set 
${}^\sharp f=(\Id\times f)_*^{-1}(f\times\Id)^*\,:
\,\Kb^G(\Zc_1)\to\Kb^G(\Zc_2)$,
if $(f\times\Id)^*$ is invertible we set 
$f_\sharp=(f\times\Id)^{*-1}(\Id\times f)_*\,:
\,\Kb^G(\Zc_2)\to\Kb^G(\Zc_1)$.
Both maps are $\star$-homomorphisms.

\vskip2mm

\noindent(A.3.3)
Fiven $g\in G$, let $i\,:\,\Tc^g\to\Tc$ be the natural inclusion.
Assume that $G$ is a diagonalizable group and $i$ is a closed immersion
(for instance if $\Tc$ is separated).
The map $(i\times i)_*\,:\,\Kb^G(\Zc^g)_g\to\Kb^G(\Zc)_g$
is invertible by the Thomason concentration theorem.
Let $N$ be the conormal bundle of $\Tc^g$ in $\Tc$.
The element $\l(N)\in\Kb_G(\Tc^g)_g$ is invertible, see (A.1.3). 
The concentration map relative to $\Tc^2$ is 
$$\cb_g=(1\boxtimes\l(N))\otimes(i\times i)_*^{-1}
\,:\,\Kb^G(\Zc)_g\to\Kb^G(\Zc^g)_g,$$
see \cite{CG, 5.11}.
It commutes with $\star$ 
(indeed $\cb_g=i^\sharp$, see (A.3.1)).

\subhead A.4. Induction\endsubhead
Let $G,$ $\Tc$, $\Zc\subset \Tc^2$ be as in A.2.
Let also $\Sc=\Tc^3$, $\Uc=\Tc^2$,
and $q_{ij}\,:\,\Sc\to\Tc^2$ be the obvious projection.
We will assume that $\Tc$ is separated.
Fix another linear group $H$.
Let $\Pc$ be a $H$-equivariant $G$-torsor over a smooth $H$-scheme $\Xc$.
Set $\Tc_\Pc=\Pc\times_{G}\Tc$.
Idem for $\Sc_\Pc$, $\Zc_\Pc$, $\Uc_\Pc$.
The $H$-schemes $\Uc_\Pc$, $\Tc_\Pc$, $\Sc_\Pc$ are smooth because
the $G$-action on $\Pc$ is locally free and $\Xc$ is smooth.
Moreover $\Zc_\Pc$ is closed in $\Uc_\Pc$.

\vskip2mm

\noindent(A.4.1)
The obvious projection $g_{ij}\,:\,\Sc_\Pc\to\Uc_\Pc$ is smooth,
because the base change of a smooth map is again smooth.
Set ${}^{ij}\tilde\Zc_\Pc=g_{ij}^{-1}(\Zc_\Pc).$
The $\star$-product
$\Kb^{H}(\Zc_\Pc)\times\Kb^{H}(\Zc_\Pc)\to
\Kb^{H}(\Zc_\Pc)$
relative to $\Sc_\Pc$ is well-defined.

Assume that $\Xc$ is a separated scheme.
Then the map $\delta\,:\,\Zc_\Pc\to(\Tc_\Pc)^2$ composed of 
the inclusion $\Zc\subset\Tc^2$
and the diagonal $\Pc\to\Pc^2$
is a closed immersion.
Set $\Zc^\delta_\Pc=\delta(\Zc_\Pc)$.
There is a $\star$-product
$\Kb^{H}(\Zc^\delta_\Pc)\times
\Kb^{H}(\Zc^\delta_\Pc)\to
\Kb^{H}(\Zc^\delta_\Pc)$
relative to $(\Tc_\Pc)^3$.
The base change and the projection formula imply that 
$\delta_*$ is a ring isomomorphism
$\Kb^{H}(\Zc_\Pc)\to\Kb^{H}(\Zc_\Pc^\delta)$.

\vskip2mm

\noindent(A.4.2)
Let $\varphi\,:\,\Pc\times\Zc\to\Zc$
and $\psi\,:\,\Pc\times\Zc\to\Zc_\Pc$
be the natural projections.
The composition of 
the obvious map $\Kb^{G}(\Zc)\to\Kb^{G\times H}(\Zc)$
and of $(\psi^*)^{-1}\circ\varphi^*$
is a $\star$-homomorphism 
$\Kb^{G}(\Zc)\to\Kb^{H}(\Zc_\Pc)$.

\vskip2mm

\noindent(A.4.3)
Assume that $H$ is diagonalisable. 
Fix $h\in H$, and write
$\Zc^h_\Pc$, $\Tc^h_\Pc$, $\Uc_\Pc^h$ 
for $(\Zc_\Pc)^h$, $(\Tc_\Pc)^h$, $(\Uc_\Pc)^h$. 
Assume that $\Xc^h$ is a closed separated subset of $\Xc$. 
Then $\delta$ restricts to a closed immersion 
$\Zc_\Pc^h\to(\Tc_\Pc^h)^2$. 
We endow $\Kb^{H}(\Zc_\Pc^h)$
with the $\star$-product relative to $(\Tc_\Pc^h)^3$.
Let $i\,:\,\Uc_P^h\to\Uc_\Pc$ be the natural inclusion.
It is a closed immersion.
The concentration map relative to $\Uc_\Pc$ is 
$$\cb_h=\l(N)\otimes i_*^{-1}
\,:\,\Kb^H(\Zc_\Pc)_h\to\Kb^H(\Zc_\Pc^h)_h,$$
where $N$ is the pull-back by the 2-nd projection
$\Uc_\Pc\to\Tc_\Pc$ conormal bundle to $\Tc^h_\Pc$ in $\Tc_\Pc$.
The following are easily checked :

(i) $\cb_h$ commutes with $\star$,

(ii) if $\Xc$ is separated then 
$\delta_*\circ\cb_h\circ(\delta_*)^{-1}$ 
is the concentration relative to $(\Tc_\Pc)^2$.

\vskip2mm

\noindent(A.4.4)
Let $\Xc'\to\Xc$ be a smooth $H$-map.
Set $\Pc'=\Xc'\times_\Xc\Pc$.
The 1st projection $\Pc'\to\Xc'$ is a $H$-equivariant $G$-torsor.
The 2nd projection $q\,:\,\Pc'\to\Pc$ is smooth. 
Set $\Tc_{\Pc'}=\Pc'\times_{G}\Tc$, etc. 
Set $q_\Tc=q\times\Id\,:\,\Tc_{\Pc'}\to\Tc_\Pc$,
and
$q_\Zc=q\times\Id\,:\,\Zc_{\Pc'}\to\Zc_\Pc$.
The pull-back
$(q_\Zc)^*\,:\,\Kb^{H}(\Zc_\Pc)\to\Kb^{H}(\Zc_{\Pc'})$
makes sense because $q_\Zc$ is smooth.
The following are easily checked :

(i) $(q_\Zc)^*$ is a $\star$-homomorphism,

(ii) if $\Xc^h$, $(\Xc')^h$ are closed and separated then
$\cb_h\circ(q_\Zc)^*=(q_h)^\sharp\circ\cb_h$,
where $q_h\,:\,\Tc_{\Pc'}^h\to\Tc_\Pc^h$ is the restriction of $q_\Tc$,

(iii) if $\Xc$ is separated then
$\delta_*\circ(q_\Zc)^*\circ(\delta_*)^{-1}=(q_\Tc)^\sharp$.

%
%\proclaim{Lemma}
%$\pi^*$ is a $\star$-homomorphism.
%We have ${}^\sharp p\circ\delta_*=\delta_*\circ\pi^*$.
%\endproclaim
%
%\noindent{\sl Proof:}
%The first claim is immediate.
%For the second one we must prove that
%$$(q_2\times\Id_{\Pc\times\Tc})^*\circ(\delta_{\Pc}\times\Id_\Tc)_*(E)=
%(\Id_{\Pc'}\times q_2\times\Id_\Tc)_*\circ(\delta_{\Pc'}\times\Id_\Tc)_*\circ
%(q_2\times\Id_\Tc)^*(E)$$ 
%where $\Tc=\Tc_1\times\Tc_1$ and $E\in\Kb^G(\Pc\times\Tc)$.
%This is obvious by base change.
%\qed
%
%\vskip2mm

\subhead A.5. Ind-objects and Pro-objects.
Definition of $\Kb_G(\Tc)$\endsubhead
From now on
let $\Tc^k_\ell$ be Noetherian $G$-schemes for $k,\ell\geq 0$,
with commutative diagrams of affine $G$-maps 
$$\matrix
\Tc^{k_2}_{\ell_2}
&{\buildrel p^{k_2k_1}_{\ell_2}\over\lra}&
\Tc^{k_1}_{\ell_2}\cr
{\ss p^{k_2}_{\ell_2\ell_1}}\bda\qquad&
&\qquad\bda{\ss p^{k_1}_{\ell_2\ell_1}}\cr
\Tc^{k_2}_{\ell_1}
&{\buildrel p^{k_2k_1}_{\ell_1}\over\lra}&
\Tc^{k_1}_{\ell_1}
\endmatrix
\leqno(A.5.1)$$
satisfying the obvious composition rules,
for any $k_2\geq k_1$, $\ell_2\geq\ell_1$.
We write $p^{k_2k_1}$ for $p^{k_2k_1}_{\ell}$, and  
$p_{\ell_2\ell_1}$ for $p^{k}_{\ell_2\ell_1}$.
We write also
$p^{k_2k_1*}$,
$p_{\ell_2\ell_1}^*$,
$p^{k_2k_1}_*$,
for
$(p^{k_2k_1})^*$,
$(p_{\ell_2\ell_1})^*$,
$(p^{k_2k_1})_*$.
Let $\Tc$ be the scheme representing the projective
system $(\Tc^k_\ell, p^{k_2k_1}, p_{\ell_2\ell_1})$.
Let $p^k_\ell\,:\,\Tc\to\Tc^k_\ell$ be the canonical map.
We set 
$$\Kb_G(\Tc)=
\ind_{k,\ell}\bigl(\Kb_G(\Tc^k_\ell),p^{k_2k_1*},p_{\ell_2\ell_1}^*\bigr).$$
The tensor product by vector bundles endows 
$\Kb_G(\Tc)$ with the structure of a $\Rb_G$-algebra.
If $p_{\ell_2\ell_1}$, $p^{k_2k_1}$ 
are torsors over a vector bundle
then $\Kb_G(\Tc)\simeq\Kb_G(\Tc^0_0).$
Given two pro-$G$-schemes $({}^1\Tc^k_\ell),$ $({}^2\Tc^k_\ell)$ as above
and $G$-maps $f^k_\ell\,:\,{}^2\Tc^k_\ell\to{}^1\Tc^k_\ell$ 
commuting to $p_{\ell_2\ell_1}$ and $p^{k_2k_1},$
we have the $\Rb_G$-algebra homomorphism 
$f^*\,:\,\Kb_G({}^1\Tc)\to\Kb_G({}^2\Tc)$
such that
$(E^k_\ell)\mapsto((f^k_\ell)^*(E^k_\ell))$.

If $\Cb$ is a category, let $\Cb^\circ$ be the opposite category
and $\Indb(\Cb)$, $\Prob(\Cb)$ 
the categories of ind- and pro-objects in $\Cb$, see \cite{GV} for background.
Objects in $\Indb(\Cb)$ are denoted by $``\ind_i"C_i$ where $(C_i)$ 
is a filtering inductive system over $\Cb$,
while objects in $\Prob(\Cb)$ are denoted by $``\pro_i"C_i$ where $(C_i)$ 
is a filtering projective system over $\Cb$.
Recall that $\Indb(\Cb)$ can be considered as a full subcategory
in the category of functors $\Cb^\circ\to\Setsb$, while
$\Prob(\Cb)=\Indb(\Cb^\circ)^\circ.$ 
Two pro-objects (resp. ind-objects) are isomorphic if they are
isomorphic in $\Prob(\Cb)$ (resp. $\Indb(\Cb)$).
See \cite{KT1, Section 1} for more details on pro-schemes, 
including pro-smoothness.
Let us quote the following fact from 
\cite{EGAIV, Corollary 8.13.2}.

\proclaim{Proposition}
For any scheme $X$,
the category of pro-objects in the category of affine schemes of finite type
over $X$ is equivalent to the category of affine schemes over $X$ 
via the functor $``\pro_i"C_i\mapsto\pro_iC_i$.
\endproclaim

\noindent
In particular, two pro-objects 
$(\Tc^k_\ell, p^{k_2k_1}, p_{\ell_2\ell_1})$ and
$('\Tc^k_\ell, '\!\!p^{k_2k_1}, '\!\!p_{\ell_2\ell_1})$
as above representing the same scheme $\Tc$ are necessarily isomorphic.

\subhead A.6. Definition of $\Kb^G(\Zc)$\endsubhead
From now on we assume that the square (A.5.1) is Cartesian, 
$\Tc^k_\ell$ is smooth, and $p^{k_2k_1}$, $p_{\ell_2\ell_1}$
are smooth locally trivial fibrations.

Given a $G$-subscheme $\Zc\subset\Tc^2$ we 
put $\Zc^k_\ell=(p^k_\ell\times p^{k}_{\ell})(\Zc)$.
Assume that $\Zc^k_\ell$
is closed in $(\Tc^k_\ell)^2$,
that $\Id\times p_{\ell_2\ell_1}$ restricts to an isomorphism
$\Zc^k_{\ell_2}\to(p_{\ell_2\ell_1}\times\Id)^{-1}(\Zc^k_{\ell_1})$,
and that $p^{k_2k_1}\times\Id$ restricts to an isomorphism
$\Zc^{k_2}_\ell\to(\Id\times p^{k_2k_1})^{-1}(\Zc^{k_1}_\ell).$
The map $(p^{k_2k_1})^\sharp$ commutes with ${}^\sharp p_{\ell_2\ell_1}$
by base change.
We set 
$$\Kb^G(\Zc)=\ind_{k,\ell}(\Kb^G(\Zc^k_\ell),
(p^{k_2k_1})^\sharp,{}^\sharp p_{\ell_2\ell_1}).$$

If $p_{\ell_1\ell_2}$, $p^{k_1k_2}$ are torsors over vector bundles
then $\Kb^G(\Zc)$ is isomorphic to $\Kb^G(\Zc^0_0)$.
%Note that we do not impose that
%the subset $\Zc\subset{}^1\Tc\times{}^2\Tc$ is closed. 

\subhead A.7. The $\star$-product on $\Kb^G(\Zc)$\endsubhead
Given ${}^{12}\Zc$, ${}^{23}\Zc$, ${}^{13}\Zc$ as above, 
assume that each
${}^{ij}\Zc^k_\ell$ is closed in $(\Tc^k_\ell)^2$,
and that $q_{13}$ restricts to a proper map
$q_{12}^{-1}({}^{12}\Zc^k_\ell)\cap
q_{23}^{-1}({}^{23}\Zc^k_\ell)\to{}^{13}\Zc^k_\ell.$
If $(E^k_\ell)\in\Kb^G({}^{12}\Zc)$ and
$(F^k_\ell)\in\Kb^G({}^{23}\Zc)$,
then $(E^k_\ell\star F^k_\ell)$ lies in $\Kb^G({}^{13}\Zc)$.
It yields the $\star$-product
$\Kb^G({}^{12}\Zc)\times\Kb^G({}^{23}\Zc)\to\Kb^G({}^{13}\Zc)$
relative to $\Tc^3$.

\vskip2mm

\head B. Appendix : infinite dimensional vector spaces.
\endhead

Let $\kb$ be a field, $V$ be a $\kb$-vector space.
Let $\Lb(V)$ be the $\kb$-algebra consisting of linear endomorphisms of $V$.
It is a topological ring for the {\sl finite topology}, i.e. such that
a basis of neighborhoods of the element $\phi$ is formed by the subsets
$$\{\phi'\,;\,\phi'(x)=\phi(x),\,\forall x\in S\}$$
where $S\subset V$ is finite.

\proclaim{Proposition}
A smooth and simple representation of the topological 
ring $\Lb(V)$ is isomorphic to $V$.
\endproclaim

\noindent{\sl Proof:}
Let $W$ be a smooth representation of $\Lb(V)$.
It gives a continuous ring homomorphism $\rho\,:\,\Lb(V)\to\Lb(W)$.
If $\Ker(\rho)\neq\{0\}$ then it contains the two-sided ideal
$\Lb(V)_f\subset\Lb(V)$ consisting of endomorphisms of finite rank,
hence $\rho=0$ because $\Ker(\rho)$ is closed and $\Lb(V)_f$ is dense.
Therefore, if $W$ is simple
$\rho$ is an isomorphism of $\Lb(V)$ onto a dense subring in $\Lb(W)$
by the Jacobson density theorem. 
Then, by \cite{J, Chap. IX, \S 11, Theorem 6 and 7} there is an isomorphism
$\psi\,:\,V\to W$ such that $\rho(x)=\psi\circ x\circ\psi^{-1}$ for all $x$.
We are done.
\qed

\vskip3mm

\head Index of notations\endhead

\noindent{\bf 1.4.} 
$G^\circ$, $\Kb_G(\Xc)$, $\Kb^G(\Xc)$, $\Kb_i^G(\Xc)$, 
$\Kb^i_G(\Xc)$, $\Rb_G$, $\Jb_g$, $\CC_g$,
$\Mb_g$, $\delta_\Xc$, $q_{ij}$, $q_i$, $\CC^ \times_\rho$, 
$\la g\ra$, $N_G(H)$,
$\Hb_*(\Xc)$, $\Kb_{i,\top}^G(\Xc)$, $\Kb^G_\top(\Xc)$, $\Db(\Xc)$,
$\Db_G(\Xc)$, $\CC_\Xc$, $\DD_\Xc$, $IC_\Xc$, $V\dot= W$,
$\Fb_V$, $\Db_\con(V)$, $\CC^\times_V$, $\Modb_\Ab$, $c^\bullet$

\vskip2mm
\noindent{\bf 2.0.} 
$\Gu$, $\Hu$, $\genu$, $\genu_\a$, $\henu$, $\benu$, $\Phiu$, $\Xu$, $\Yu$,
$\Xu^\vee$, $\Yu^\vee$, $\a_i$, $\o_i$, $\a_i^\vee$, $\o_i^\vee$,
$\theta$, $\theta^\vee$, $W$, $\Wu$, $\ell(w)$, $I$, $\Iu$, $X$, $Y$,
$X^\vee$, $Y^\vee$, $\delta$, $\a_0$,
$\a_0^\vee$, $c$, $\o_0$, $\o_0^\vee$, $(\ :\ )$, $m_{ij}$, $\O$

\vskip2mm
\noindent{\bf 2.1.}
$K$, $K_d$, $\eps$, $\gen$, $\hen$, $\Phi$, $\Phi^+$, $\Phi^-$,
$\gen_\a$, $\Theta$, $\uen(\pm\Theta)$, $\Theta_\ell$, $\uen_\ell$,
$\uen^-_\ell$, $\uen$, $\uen^-$ 

\vskip2mm
\noindent{\bf 2.2.}
$G'$, $G_{KM}$, $G$, $G_d$, $\CC^\times_\delta$, $\CC^\times_{\o_0}$,
$H'$, $H$, $X'$, $W_s$, $U$, $U^-$, $U(\Theta)$, $U(-\Theta)$,
$U_\ell$, $U^-_\ell$, $U_{\pm\a}$, $B$, $B^-$

\vskip2mm
\noindent{\bf 2.3.}
$\Xc$, $\Gc$, $\pi$, $\dot W$, $1_\Gc$, $1_\Xc$, $\Xc_y$, $\Xc^w$, $\Uc^w$,
$\Xc_{\leq y}$, ${}^w\Xc$, ${}^w\Gc$, $y''\succeq(y,y')$

\vskip2mm
\noindent{\bf 2.4.}
${}^w\Xc^k$, $p^{k_2k_1}$, $p^k$

\vskip2mm
\noindent{\bf 2.5.}
$\Tc$, $\rho\,:\,\Tc\to\Xc$, $\Vc^w$, ${}^w\Tc$, $[g:x]$, $\Tc_y$

\vskip2mm
\noindent{\bf 2.6.}
${}^w\Tc^{k\ell}$, $\Tc^\ell$, $p_\ell$,
$p^{k_2k_1}_\ell$, $p^k_{\ell_2\ell_1}$,
$p^k_\ell$

\vskip2mm
\noindent{\bf 2.7.}
$\ben$, $\Bc$, $\ben_w$ 

\vskip2mm
\noindent{\bf 2.8.}
$\dot\gen$, $p\,:\,\dot\gen\to\gen$, $\Bc_x$, $\Nc$, $\dot\Nc$, 

\vskip2mm
\noindent{\bf 2.9.}
$\Phi_J$, $\Phi_J^\pm$, $U_J$, $\uen_J$, $G_J$, $B_J$, $\ben_J$, $W_J$,
$w_J$, ${}^JW$, $W^J$, $\Xc^J$, ${}^w\Xc^J$, ${}^w\Xc^{J,k}$,
$\Bc^J$, $\pi_J$, $\Bcu$, $\dot\genu$, $p\,:\,\dot\genu\to\genu$,
$\Ncu$, $\dot\Ncu$, $\Bcu_x$, $\Zcu$, $\Bc_J$, $\dot\Nc_J$, $\Zc_J$,
$\Oc_{J,y}$, $\Zc_{J,y}$, $\Zc_{J,\leq y}$, ${}^J\Bc_y$, ${}^J\Bc$,
$\flat$, $!\,:\,\Bc\to\Bcu$, $\pen^!$

\vskip2mm
\noindent{\bf 2.10.}
$\Pc^J$, ${}^w\Pc^{J,k}$

\vskip2mm
\noindent{\bf 2.11.}
$\ad_{(g,\zeta)}x$, $\CC^\times_q$, $A$

\vskip2mm
\noindent{\bf 2.12.}
$\su$, $s'$, $s$, $a$, $G^s$, $\Bc^s$, $\gen^s$, $\Tc^a$ 

\vskip2mm
\noindent{\bf 2.13.}
${}^w\Xc^s$, ${}^w\Tc^a$, $\Xc_i^s$, $\Xi_s$, $\dot\Nc_i^a$

\vskip2mm
\noindent{\bf 3.1.}
$\Hb$, $t_i$, $x_\l$, $\Rb_H$, $\Rb_A$, $t_y$, $IM$, $\Rb_{\tau,\zeta}$,
$\Hb_{\tau,\zeta}$, $\Hb_{\leq w}$, $\Hb^\ch$

\vskip2mm
\noindent{\bf 3.2.}
$\Hb_J$, $\Hbu$, $\Hbu_\zeta$

\vskip2mm
\noindent{\bf 4.1.}
$\Oc_y$, $\Oc_{\leq y}$, ${}^w\Oc_y$, ${}^w\Oc_y^k$

\vskip2mm
\noindent{\bf 4.2.}
$\Zc_y$, $\Zc_{\leq y}$, $\Zc$, ${}^w\Zc_y$, ${}^w\Zc_y^{k\ell}$, $\Zc_y^\ell$

\vskip2mm
\noindent{\bf 4.3.}
$\Zc^a_{ij}$

\vskip2mm
\noindent{\bf 4.4.}
${}^w\Zc_{\leq y}^a$, $D\subset A$, $h^{w'w}$, $h_{yy'}$,
$\Kb^D(\Zc^a_{\leq y})$, $\Kb^D(\Zc^a)$, $\Kb_a$, $\Kb_{a,\leq y}$

\vskip2mm
\noindent{\bf 4.5.}
$\Kb^A({}^w\Zc^0_{\leq w_J})$, $\cb_a$

\vskip2mm
\noindent{\bf 4.6.}
${}^wL_\l^k$, ${}^wL_\l^{k\ell}$, $\Kb_H(\Xc)$, $\Kb_A(\Tc)$,
$\Kb_H({}^w\Xc)$, $\Kb_A({}^w\Tc)$, ${}^wL_\l$, $L_\l$,
$\Kb_H(\Xc^s)$, $\Kb_A(\Tc^a)$, $\Lc_\l$

\vskip2mm
\noindent{\bf 4.7.}
$D_\l$, ${}^wD^{k\ell}$, ${}^wD^0_\l$, $D_i$, $\bar\Zc_{s_i}$,
$({}^w\Zc_{s_i})^{k\ell}$, ${}^wD^{k\ell}_i$, ${}^wD_i^0$, $p_a$,
$D_{\l,a}$, $D_{i,a}$, ${}^\sharp p_a$

\vskip2mm
\noindent{\bf 4.8.}
${}^w\beta_s$, $I_{w,s}$, $\Ib_a$, $\hat\Rb_a$, $W(a)$

\vskip2mm
\noindent{\bf 4.9.}
$\Psi_a$, $\Theta_J$, ${}^w\Theta_J$, $\phi_k$, $\psi_k$, $\Gamma_k$, $\Gamma$,
$\theta_i$, ${}^w\Psi_a$, $\hat\Psi_a$, $\hat\Hb_a$

\vskip2mm
\noindent{\bf 5.1.}
$\gen_{KM}$, $\phi=\{e,f,h\}$, $\gen_\phi$, $G_\phi$, $Z_\phi(G)$,
$A_e$, $C_\phi$, $s_\phi$

\vskip2mm
\noindent{\bf 5.3.}
$\zen(f)$, $\senu_\phi$, $\Scu_\phi$, $\sen_\phi$, $\Sc_{\phi,y}$,
$\Sc_{\phi,\le y}$, $\Sc_{\phi,<y}$, $\Sc_\phi$, $\Sc_{\phi,y}^{k,\ell}$,
$\g_y^{k\ell}$, $\g_y$, $\g_{y,a}$

\vskip2mm
\noindent{\bf 5.4.}
$\Kb^D(\Sc^a_{\phi,\leq y})$

\vskip2mm
\noindent{\bf 5.5.}
$(S_\QQ)$, $\Bcu_e^s$, $(T_H)$

\vskip2mm
\noindent{\bf 5.6.}
$\tilde\cb_a$, $\tilde\Theta_\flat$

\vskip2mm
\noindent{\bf 5.7.}
$\sigma_i$, $\sigma'_i$ 

\vskip2mm
\noindent{\bf 5.8.}
$\Upsilon$, $\Upsilon_y$

\vskip2mm
\noindent{\bf 6.1.}
$L_{a,i}$, $S_{a,\chi}$, $\Lb_{a,\chi,i}$, $\Lb_{a,\chi}$, $\{\chi\}$, $L_a$,
$\Extb(L_a,L_a)$, $\Phi_S$, $\Phi_a$, $nilp$, RR

\vskip2mm
\noindent{\bf 6.2.}
$\Bc_e^s$, $\Bc_{e,\le y}^s$, $\Hb_*(\Bc_e^s)$, $G(s,e)$, $\Pi(s,e)$,
$\Nb_{a,e,\chi}$, $\Pi(s,e)^\vee$, $\iota$

\vskip2mm
\noindent{\bf 6.3.}
$\Pi(a)^\vee$, $v\,:\,\CC^\times\to\QQ$, $\phiu=\{\eu,\fu,\hu\}$,
$\Mu$, $\Mu_\phiu$, $\gen_{t,i}$, $Q$, $L$, $\qen$, $\len$,
$\hat\Qc$, $\hat\Nc^a$, $\hat\Bc^s_e$, 
$L(s,e)$, $Z^\circ_L$, $M^s$, $M^s_\phi$, $M^s_u$, $\hat\Xi_s$

\vskip2mm
\noindent{\bf 6.5.}
$\Piu(\su,\eu)$, $\Gu(\su,\eu)$, $\Pi(\su,\eu)^\vee$, $\Piu(a)^\vee$,
$\Lbu_{a,\chiu}$, $\Gu(\su,\phiu)$, $\Mbu_{a,\phiu,\chiu}$,
$\underline n_{a,\chiu,\chiu'}$

\vskip2mm
\noindent{\bf 6.7.}
$n_{a,\chiu,\chi}$, $\eps$

\vskip2mm
\noindent{\bf 7.1.}
$\Bc_L$, $\Nc_L$, $\dot\Nc_L$, $\dot\len$, $\gen_\pm$, $\len_\pm$,
$\dot\gen_\pm$, $\dot\len_\pm$,
$\Ec_\pm$, $\Fc_\pm$, $\Indb_{\qen_\pm}^{\gen_\pm}$, $\dot\Nc_{\pm}$,
$\dot\Nc_{L,\pm}$ 

\vskip2mm
\noindent{\bf 7.3.}
$\la\ :\ \ra$

\vskip2mm
\noindent{\bf 7.4.}
$\Pi(a)^\vee$, $\check a$, $\Lb_{\check a,\check\chi}$, $\bar a$ 

\vskip2mm
\noindent{\bf 7.5.}
$m_{a,\chiu,\chi}$

\vskip2mm
\noindent{\bf 7.6.}
$\Modb^{int}_{\Hb_{\tau,\zeta}}$

\vskip2mm
\noindent{\bf 8.2.}
$\sigma=(\sigma_{ij})$, $\sigma[n]$, $\{(\sigma_a),(z_a)\}$, $\Sc$,
$\Sc/\sim$, $\Cc$, $\Cc/\sim$

\vskip2mm
\noindent{\bf 8.3.}
$\Hb'_{\tau,\zeta}$, $\Hb^\ch_{\tau,\zeta}$, $\{\Lb,\a\}$, $\Mb^\pi$,
$\{\Mb^\pi,\beta\}$

\vskip2mm
\noindent{\bf 8.4.}
$\Cc_d$, $\Sc_d$, $\Qb_m$, $Q_d/\sim$, ${}^sd_i$,
$\Rep_{{}^sd}(\Qb_m)$

\vskip2mm
\noindent{\bf 8.5.}
$\Ub_m$, $\Sc_\infty$, $\Bb$, $\bb_x$, $\Fb$, $\fb_x$, $\db$, $\fb_\db$,
$\Ub_\infty$, $\fbu_i$, $\Qb_\infty$, $\Bbu$, $\Fbu$, $\bbu_x$, $\fbu_x$,
$\dbu$, $\fb_\dbu$, $\Hb_d$, $\Gb_d$, $\Gb_\infty$, $\Tb_\infty$, $\Lb_x$,
$\Lb^x$, $\Mb_x$, $\Mb^x$, $\Ccu_d[\ell]$, $\Ccu_d$, $\Scu_d$, $\Scu_\infty$,
$\Hbu_d$, $\Gbu_d$, $\Gbu_\infty$, $\Tbu_\infty$, $\Lbu_x$, $\Lbu^x$,
$\Mbu_x$, $\Mbu^x$, $\Delta$, $\hat\Ub_\infty$, $m_{x,\xu}$

\vskip2mm
\noindent{\bf 9.}
$\Modb_a$, $\Modbu_a$

\vskip2mm
\noindent{\bf 9.3.}
ht, $h$, $r$, $\rho^\vee$, $\Pi_k$, $e_k$, $a_k$

\vskip2mm
\noindent{\bf 9.4.}
$\Hb'_J$, $\Nc_J$, $p_J\,:\,\dot\Nc_J\to\Nc_J$, 
${}^w\Tc_J^{k\ell}$, ${}^w\Tc_J$, 
$p\,:\,{}^w\Tc^{k\ell}\to{}^w\Tc_J^{k\ell}$,
$\Uc_J$, $\Uc$, $\Zc_\Uc$, $\Bc_{e,\leq w}$,
$\Hb_*(\Bc_e)$, $\psi$

\enddocument